%% file: main.tex
\documentclass{article}

\input{packages}
\input{commands}

\begin{document}

\title{Randomized  orthogonalization and Krylov subspace methods: principles and algorithms}

\author{Jean-Guillaume de Damas\thanks{Sorbonne University, Paris, France. Work done while the author was at INRIA, Paris, France, \texttt{jean-guillaume.de-damas@inria.fr}} \and
Laura Grigori\thanks{Institute of Mathematics, EPFL, Lausanne, and PSI Center for Scientific Computing, Theory and Data, Villigen PSI, Switzerland, \texttt{laura.grigori@epfl.ch}} \and
Igor Simunec\thanks{Institute of Mathematics, EPFL, Lausanne, Switzerland, \texttt{igor.simunec@epfl.ch}} \and
Edouard Timsit\thanks{Éducation Nationale, Académie de Versailles, France. Part of this work was done while the author was at INRIA, Paris, France. \texttt{edouard-gaston.timsit@ac-versailles.fr}}
}

\date{}

\maketitle
\abstract{
We present an overview of randomized orthogonalization techniques that construct a well-conditioned basis whose sketch is orthonormal. Randomized orthogonalization has recently emerged as a powerful paradigm for reducing the computational and communication cost of state-of-the-art orthogonalization procedures on parallel architectures, while preserving, and in some cases improving, their numerical stability.    
 This approach can be employed within Krylov subspace methods to mitigate the cost of orthogonalization, yielding a randomized Arnoldi relation. We review the main variants of the randomized Gram--Schmidt and Householder QR algorithms, and discuss their application to Krylov methods for the solution of large-scale linear algebra problems, such as linear systems of equations, eigenvalue problems, the evaluation of matrix functions, and matrix equations. 
}

\section{Introduction}
\input{introduction}

\input{sketching}
\input{randqr}

\input{rgs}

\input{householder}

\input{randorth-extra}

\input{krylov}

\input{linearsystems}

\input{eigenvalues}

\input{matrix-functions}

\input{matrix-equations}

\section{Conclusions}
\label{sec:conclusions}

The randomized orthogonalization framework provides efficient ways to construct sketch-orthogonal bases that are very well-conditioned. The randomized Gram--Schmidt and Householder QR algorithms have excellent numerical stability properties and significantly lower communication costs on parallel architectures. These techniques can be employed within Krylov subspace methods to lower the cost of orthogonalization, and the resulting randomized Arnoldi relation can be used to construct approximate solutions to linear systems of equations, eigenvalue problems, to evaluate matrix functions, and to solve matrix equations. Two main approaches can be identified in this context. On one hand, we can construct a sketch-orthogonal basis explicitly, for instance via a randomized Gram--Schmidt process, which requires roughly half the number of flops with respect to the deterministic Gram--Schmidt process, has communication costs which are comparable to those of CGS, and numerical stability which is comparable to that of MGS. On the other hand, we can construct a non-orthogonal basis with a cheap procedure, such as the $k$-truncated Arnoldi process, and then obtain a sketch-orthogonal basis implicitly through whitening, by computing a QR factorization of the sketched basis. This second approach is asymptotically cheaper than randomized Gram--Schmidt, as the sketch-orthogonal basis is never formed explicitly, but its main issue is that the basis constructed with $k$-truncated Arnoldi quickly becomes ill-conditioned, and this may have a negative impact on the convergence of the approximate solutions extracted from the Krylov subspace. Nevertheless, although less robust than the first, the approach employing whitening is computationally efficient, and it very often performs well in practice, although this behavior is still not completely understood theoretically. 
Several research directions remain still open. For instance, it is not known whether it is possible to obtain a sketch-orthogonal basis of a Krylov subspace with an algorithm that has the same numerical stability as randomized Gram--Schmidt, and the same computational efficiency as whitening. Moreover, when the matrix $A$ is symmetric, sketch-orthogonalization in the Arnoldi process does not yield a short-term recurrence in general, destroying the symmetry of the projected matrix. This phenomenon is undesirable especially when computing eigenvalues of a symmetric matrix, since the eigenvalues of the projected matrix obtained with the randomized Arnoldi process are not guaranteed to be real. However, it is unclear if the randomized orthogonalization process can be adapted in order to preserve the symmetry of the small projected matrix. Lastly, the development of a standard library for randomized orthogonalization routines would enable their use in real applications, thus allowing to more easily benchmark and gain feedback on the numerical behavior of these algorithms in large-scale applications.

\section*{Acknowledgment}
The first, second, and fourth authors acknowledge funding from the European Research Council (ERC) under the European Union’s Horizon 2020 research and innovation program, grant agreement No 810367.

\bibliographystyle{siam}
\bibliography{references}

\end{document}

%% file: packages.tex
\usepackage[T1]{fontenc}
\usepackage[utf8]{inputenc}
\usepackage[english]{babel}
\usepackage{geometry}

\usepackage{graphicx} 

\usepackage{amsmath}
\usepackage{amsfonts}
\usepackage{amsthm}   
\usepackage{mathtools}
\usepackage{xcolor}

\usepackage{algorithm}
\usepackage{algpseudocode}

\usepackage{cite}		
\usepackage{hyperref}
\usepackage{cleveref}

\usepackage{graphicx}
\usepackage{tikz}
\usepackage{tkz-euclide}  
\usepackage{tikz-3dplot}
\usepackage{pgfplots}
\usetikzlibrary{pgfplots.groupplots}
\usetikzlibrary{arrows.meta}
\pgfplotsset{compat=1.17} 
\usepackage{circuitikz}
\usepackage{subcaption}

%% file: commands.tex
\theoremstyle{definition}
\newtheorem{theorem}{Theorem}[section]

\newtheorem{definition}[theorem]{Definition}

\crefname{equation}{}{}
\crefname{theorem}{Theorem}{Theorems}
\crefname{lemma}{Lemma}{Lemmas}
\crefname{corollary}{Corollary}{Corollaries}
\crefname{proposition}{Proposition}{Propositions}
\crefname{remark}{Remark}{Remarks}
\crefname{section}{Section}{Sections}
\crefname{figure}{Figure}{Figures}
\crefname{algorithm}{Algorithm}{Algorithms}
\crefname{table}{Table}{Tables}
\crefname{subsection}{Section}{Sections}

\DeclareMathOperator*{\argmin}{argmin}		

\DeclarePairedDelimiter{\abs}{\lvert}{\rvert}
\DeclarePairedDelimiter{\norm}{\lVert}{\rVert}

\DeclareMathOperator{\vspan}{span} 	
\DeclareMathOperator{\range}{Range}

\newcommand{\C}{\mathbb{C}}
\newcommand{\R}{\mathbb{R}}

\newcommand{\K}{\mathcal{K}} 

\newcommand{\subspace}{\mathcal{W}}
\newcommand{\Ker}[1]{\mathrm{Ker}(#1)}
\newcommand{\fl}{\mathbf{fl}}
\newcommand{\cond}{\mathrm{Cond}}	
\newcommand{\proj}[1]{\mathcal{P}_{#1}}	
\newcommand{\skproj}[1]{\mathcal{P}_{#1}^\Omega}	

\renewcommand{\vec}[1]{\boldsymbol{#1}}		
\newcommand{\mv}{\texttt{mv}}	

\newcommand{\xfom}{\widetilde{\vec x}^{\text{F}}}
\newcommand{\xgmres}{\widetilde{\vec x}^{\text{G}}}
\newcommand{\xsfom}{\vec x^{\text{F}}}
\newcommand{\xsgmres}{\vec x^{\text{G}}}
\newcommand{\yfom}{\widetilde{\vec y}^{\text{F}}}
\newcommand{\ygmres}{\widetilde{\vec y}^{\text{G}}}
\newcommand{\ysfom}{\vec y^{\text{F}}}
\newcommand{\ysgmres}{\vec y^{\text{G}}}

\graphicspath{{./}}

\newcommand{\comment}[1]{}

\newcommand{\e}{\mathbf{u}}
\newcommand{\f}{\mathbf{f}}

\definecolor{colW}{RGB}{190,80,80}    
\definecolor{colO}{RGB}{140,200,120}  
\definecolor{darkteal}{RGB}{0,96,115} 
\definecolor{colQ}{RGB}{58,107,181}   

%% file: introduction.tex
Krylov subspace methods are among the most powerful and widely used techniques for solving large-scale numerical linear algebra problems, such as linear systems of equations, eigenvalue problems, matrix equations, and matrix function evaluations. Given a matrix $A \in \R^{n \times n}$ and a vector $\vec b \in \R^{n}$, these methods iteratively construct the sequence of Krylov subspaces
\begin{equation*}
	\K_m(A, \vec b) = \vspan\{\vec b, A \vec b, \dots, A^{m-1} \vec b\},
\end{equation*}
which are used as search space to find an approximate solution to a given problem. An orthonormal basis of $\K_m(A, \vec b)$ is typically constructed with the Arnoldi process, which produces a decomposition 
$A Q_m = Q_{m+1} \underline{H}_m$, 
where the columns of $Q_m$ are the orthonormal basis vectors and  $\underline{H}_m$ is an upper Hessenberg matrix that contains the orthogonalization coefficients representing the projection of $A$ onto $\K_m(A, \vec b)$. Each iteration of the Arnoldi process involves a matrix-vector product with $A$ followed by orthogonalization of the new vector against all previous basis vectors through a Gram--Schmidt process, for a total computational cost of $\mathcal{O}(m \,\mv(A) + nm^2)$ for $m$ iterations, where $\mv(A)$ denotes the cost of a matrix-vector product with $A$. Consequently, for moderately large $m$, or when the computation is performed on a parallel computer, the orthogonalization step becomes the dominant cost and often limits the practical efficiency of Krylov subspace methods. 

Several strategies have been proposed in the literature to address this computational bottleneck. Incomplete or truncated orthogonalization schemes reduce the number of inner products at the expense of a loss of numerical stability. Restarting techniques limit the dimension of the Krylov subspace and periodically restart the iteration to control memory usage and computational cost, but they may incur convergence delays or stagnations \cite{Saad03}. 

Randomization has emerged as a powerful technique for solving large scale problems by enabling dimensionality reduction through random projections and subspace embeddings \cite{johnson1984extensions, sarlos2006improved, ailon2006fastjl}. It has been applied successfully to different linear algebra problems, including solving least squares problems \cite{drineas2006sampling, Rokhlin2008fastrandomized} and computing low-rank matrix approximations (see, e.g.  \cite{woodruff2014sketching, martinsson2020randomized, murray2023randomizednumericallinearalgebra} for details). More recently, randomized techniques have also been introduced in the context of Krylov subspace methods.  Randomized or sketched Krylov subspace methods replace exact orthogonalization in the Gram--Schmidt process with operations performed on vectors that belong to a low-dimensional sketched space, obtained by applying a random sketching matrix $\Omega \in \R^{d \times n}$ to the basis vectors, which acts as an oblivious subspace embedding. This approach has the potential to reduce the computational and communication costs of building the Krylov basis \cite{BalabanovGrigori2022, NakatsukasaTropp24}, and instead of an orthonormal basis $Q_m$, it produces a sketch-orthonormal basis $V_m$ such that its sketch $\Omega V_m$ contains orthonormal columns. It allows exploiting mixed-precision arithmetic and optimized computational kernels, while providing numerical guarantees with high probability. 

Variants of this randomized approach have been explored for the solution of linear systems \cite{BalabanovGrigori2022,BalabanovGrigori25block,TGB23,GuettelSimunec24,NakatsukasaTropp24}, eigenvalue problems \cite{NakatsukasaTropp24,BalabanovGrigori25block, DedamasGrigori25ira,DedamasGrigori25ks}, for the evaluation of matrix functions \cite{GuettelSchweitzer23,PSS25mf,CKN24} and the solution of matrix equations\cite{PSS25me}. This work provides a general introduction to the use of randomized orthogonalization within Krylov subspace methods and its application for the solution of different linear algebra problems. Section \ref{sec:sketchings} describes the oblivious subspace embedding property and different sketching matrices that satisfy this property.  Section \ref{sec:randomized-qr-introduction} reviews randomized orthogonalization techniques, including randomized Gram–Schmidt and randomized Householder QR. Randomized Krylov subspace methods are discussed in section \ref{sec:krylov}, while their usage to solve linear systems and eigenvalue problems is presented in Sections \ref{sec:linear-systems}, and \ref{sec:eigenvalues}, respectively. Randomized approaches for matrix functions and matrix equations are introduced in Sections \ref{sec:matrix-functions} and \ref{sec:matrix-equations}, respectively.

\subsection{Notation}
\label{subsec:notation}

We introduce here some general notation that we use throughout this manuscript. We denote matrices with uppercase letters, and vectors with bold lowercase letters. We denote by $I_n$ the identity matrix of dimension $n$, and omit the subscript when it can easily be inferred from the context. The columns of the identity matrix of dimension $n$ are denoted by $\vec e_1, \dots, \vec e_n$. We denote by $0_n$ the vector of zeros of length $n$, and by $0_{n \times m}$ a zero-matrix of size $n \times m$; occasionally, the subscripts may be omitted if there is no ambiguity on the dimensions.
We use calligraphic letters to denote a vector subspace $\mathcal{W} \subset \R^n$, and we denote by $A \mathcal{W}$ the image of $\mathcal{W}$ under the action of the matrix $A$.
We denote by $\norm{\vec x}$ the Euclidean norm of a vector $\vec x \in \R^n$, and by $\norm{A}_2$ and $\norm{A}_F$ the spectral and Frobenius norms of a matrix $A \in \R^{n\times n}$, respectively. The singular values of $A$ in nonincreasing order are denoted by $\sigma_1(A), \dots, \sigma_n(A)$.

%% file: sketching.tex
\section{Sketching and embeddings}
\label{sec:sketchings}

Randomized algorithms rely on sketching, a dimensionality reduction technique that allows to embed high dimensional subspaces into lower-dimensional ones while approximately preserving their geometry, such as inner products between vectors in the subspace \cite{sarlos2006improved}. These random linear maps (see  early references \cite{johnson1984extensions,dasgupta2003elementary,ailon2006fastjl,10.1007/3-540-45465-9_59}) preserve enough information to enable, with high probability, solving accurately a wide range of linear algebra problems. This is demonstrated, for example, in early work on overdetermined least squares problems \cite{drineas2006sampling,Rokhlin2008fastrandomized}.
We begin by introducing the definition of the $\epsilon$-embedding property \cite{sarlos2006improved, woodruff2014sketching}.
\begin{definition}
\label{defn:epsilon-embedding}
Let $\mathcal{W} \subset \R^n$ be an $m$-dimensional vector subspace and $\epsilon \in \left]0,1\right[$. We say that $\Omega \in \R^{\ell \times n}$, $m \leq \ell$, is a $\epsilon$-embedding of $\mathcal{W}$ if and only if
\begin{align} \label{eq:epsembedding}
        \forall \vec w \in \mathcal{W}, \quad (1-\epsilon) \| \vec w\|^2 \leq \| \Omega \vec w \|^2 \leq (1+\epsilon) \| \vec w\|^2. \qquad (\epsilon\text{-embedding property})
\end{align}
\end{definition}

\begin{figure}[t]
  \centering

  \begin{subfigure}{0.6\textwidth}
    \centering
    \begin{tikzpicture}[tdplot_main_coords, scale=2]

        \filldraw[fill=blue!8, draw=blue!60!black, opacity=0.8]
          (-1,-1,0) -- (1,-1,0) -- (1,1,0) -- (-1,1,0) -- cycle;
        \draw[->, thick, blue!60!black] (0,0,0) -- (0.8,0,0) node[anchor=east] {$\vec w$};
        \draw[->, thick, blue!60!black] (0,0,0) -- (0,0.8,0) node[anchor=south] {$\vec w'$};
        
        \draw (0,0,0) -- (0,0,1);
        \draw[dotted] (0,0,-1) -- (0,0,0);
        
        \node[blue!60!black, font=\Large] at (0.15,-0.45,0.2) {$\mathcal{W}$};
        \node[font=\Large] at (1.05,1.05,1.7) {$\mathbb{R}^n$};

        \draw[->, ultra thick] (0,2,0.5) -- (-0.35,3,0.5)
              node[midway, above, yshift=2pt, font=\Large] {$\Omega$};

    \end{tikzpicture}
  \end{subfigure}
  \hfill
  \begin{subfigure}{0.39\textwidth}
    \centering
    \begin{tikzpicture}[tdplot_main_coords, scale=2]

        \filldraw[fill=blue!30, draw=blue!60!black, opacity=0.8]
          (-1,-1,0) -- (1,-1,0) -- (1,1,0) -- (-1,1,0) -- cycle;
        \draw[->, blue!60!black, thick] (0,0,0) -- (0.5,-0.1,0) node[anchor=east] {$\Omega \vec w$};
        \draw[->, blue!60!black, thick] (0,0,0) -- (-0.2,0.9,0) node[anchor=south east] {$\Omega \vec w'$};
        
        \draw (0,0,0) -- (0,0,1);
        \draw[dotted] (0,0,-1) -- (0,0,0);
        
        \node[blue!60!black, font=\Large] at (0.15,-0.45,0.2) {$\Omega \mathcal{W}$};
        \node[font=\Large] at (1.05,1.05,1.7) {$\mathbb{R}^\ell$};

    \end{tikzpicture}
  \end{subfigure}
\caption{$\epsilon$-embedding of a vector subspace $\mathcal{W} \subset \R^n$, with minor distortion of norms and angles, and conservation of the dimension of $\mathcal{W}$.} \label{fig:embedding}

\end{figure}
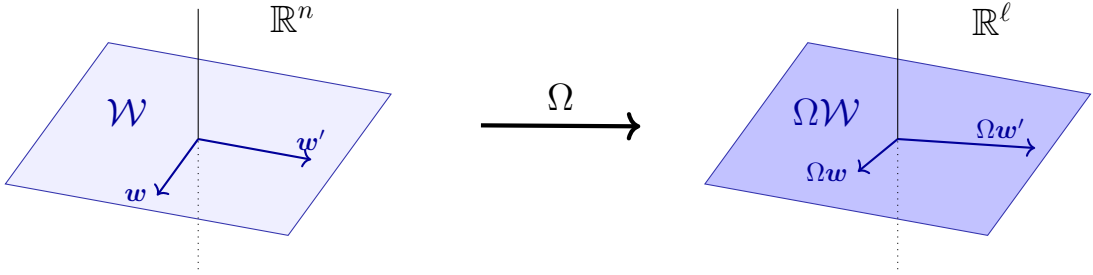

The vector $\Omega \vec w \in \Omega \mathcal{W} \subset \R^\ell$ is called the \textit{sketch} of $\vec w \in \mathcal{W} \subset \R^n$. The $\epsilon$-embedding property can be interpreted as the restriction of $\Omega$ to $\mathcal{W}$ being nearly isometric, that is, it maps $\mathcal{W}$ to a new vector subspace $\Omega \mathcal{W} \subset \R^\ell$ with very little distortion. While $\dim \mathcal{W} = \dim \Omega \mathcal{W}$, the latter lies in a vector space of much smaller dimension. We illustrate this property in~\Cref{fig:embedding}, where the norms of $\vec w,\vec w' \in \mathcal{W}$ and $\Omega \vec w, \Omega \vec w'$ are slightly different, the angles between $\vec w,\vec w'$ and $\Omega \vec w, \Omega \vec w'$ are slightly different, but the dimensions of $\mathcal{W}$, $\Omega \mathcal{W}$ are the same. Although $\Omega$ does not induce a proper inner product, as $\Omega^T \Omega$ is only positive semidefinite, it has been shown that it defines a proper norm when restricted to the embedded space $\mathcal{W}$ with $\norm{\vec w}_{\Omega^T \Omega}^2 = \vec w^T \Omega^T \Omega \vec w$ for $\vec w \in \mathcal{W}$; see, e.g., \cite{BalabanovNouyI}.

\tdplotsetmaincoords{60}{110}

When considering a matrix $W \in \R^{n \times m}$ such that $\mathcal{W} = \range(W)$, the $\epsilon$-embedding property allows to establish spectral relations between $W$ and its sketch $\Omega W$. For example, it is derived in \cite{Gilbert2012SketchedSVDRecovering} that for $j = 1,\dots,m$:
\begin{equation}
\label{eqn:singular-values-eps-embedding}
\sqrt{1-\epsilon} \leq \frac{\sigma_j(\Omega W)}{\sigma_j(W)} \leq \sqrt{1+\epsilon} \implies \mathrm{Cond}(W) \leq \sqrt{\frac{1+\epsilon}{1-\epsilon}} \cdot \mathrm{Cond}(\Omega W).
\end{equation}

Such $\epsilon$-embedding sketching matrices can be efficiently obtained by drawing from simple random distributions $\mathcal{D}$ over $\R^{\ell \times n}$ while satisfying  \cref{eq:epsembedding} with high probability, without knowledge of $\mathcal{W}$. 

\begin{definition} \cite[Definition 2.2]{woodruff2014sketching}
    Let $\mathcal{D}$ be a distribution over matrices of $\R^{\ell \times n}$ with $m \leq \ell$. We say that $\mathcal{D}$ is an $(\epsilon, \delta, m)$-oblivious subspace embedding (OSE) if and only if for any $\Omega \in \R^{\ell \times n}$ drawn from $\mathcal{D}$ and any given $m$-dimensional subspace $\mathcal{W} \subset \R^n$, $\Omega$ is an $\epsilon$-embedding of $\mathcal{W}$ with probability at least $1-\delta$.
\end{definition}
Depending on the distribution $\mathcal{D}$, it is possible to construct an OSE for which the order of magnitude of $\ell$ may be $O(\epsilon^{-2} m)$ or $O(\epsilon^{-2} m\log(m))$, the calculation $\Omega \vec x$ may be as cheap as $n \log(n)$ flops, and the storage cost of $\Omega$ may be as small as that of $\ell + n$ integers.  
Historically, such distributions $\mathcal{D}$ have been first described as yielding embeddings $\Omega$ of a finite set $E_d$ of $d$ vectors (and not of the subspace they generate). Indeed, for $\delta \in ]0,1[$, and for an integer $\ell$ greater than a modest multiple of $\epsilon^{-2}\log(d/\delta)$, the celebrated Johnson-Lindenstrauss lemma \cite{johnson1984extensions} shows that there exist distributions $\mathcal{D}$ over $\R^{\ell \times n}$ such that, for any given set $E_d$ of $d$ vectors, the following event occurs with probability at least $1-\delta$:
\[ \forall \; \vec x_i \neq \vec x_j \in E_d, \quad (1-\epsilon) \|\vec x_i - \vec x_j \|^2 \leq \|\Omega \vec x_i - \Omega \vec x_j \|^2 \leq (1+\epsilon) \|\vec x_i - \vec x_j \|^2.\]

Let us now outline some concrete OSEs. We first outline the Gaussian sketching distribution. To draw from this distribution, we simply draw each entry of the matrix $\Omega \in \R^{\ell \times n}$ independently from $\mathcal{N}(0,1)$, and scale the resulting matrix by $\ell^{-1/2}$. Provided that the sampling size $\ell$ is set to
\begin{equation}
    \ell = O \left( \frac{1}{\epsilon^2} (m + \log(1/\delta)) \right),
\end{equation}
this distribution is a $(\epsilon, \delta, m)$-OSE~\cite[Theorem 2.3]{woodruff2014sketching}. This requirement on the sampling size coincides with that of a Johnson-Lindenstrauss transform for an exponential number of arbitrary points~\cite[Theorem 2.1]{woodruff2014sketching}. In that sense, it is optimal~\cite{nelsonjloptimal}. In practice, setting the sampling size $\ell= O(m)$, we get an embedding matrix $\Omega$ with a parameter $\epsilon \approx 1/2$ with high probability \cite{martinsson2020randomized}. Despite its favorable theoretical properties, the main downside of a Gaussian sketching matrix is that it is dense and unstructured and thus costly to store and apply.

We next outline the $s$-hashing distribution. To draw $\Omega \in \R^{\ell \times n}$ from this distribution, we randomly choose $s$ entries in each column of $\Omega$, randomly set them to $\{- 1/s^{1/2},  \; + 1/s^{1/2}\}$, and set all other entries to zero, resulting in sparse columns of unit norm. This results in a sparse matrix $\Omega$ with exactly $ns$ nonzero entries. The more balanced the rows of $\mathcal{W}$, the lower we can set $s$ and $\ell$, as shown in~\cite{bourgainsparse}. 
From~\cite{chenakkod2025optimal}, the $s$-hashing ensemble is an ($\epsilon$,  $\delta$, $m$)-OSE provided that the sampling size $\ell$ and the parameter $s$ verify
\begin{equation}
    \ell = O \left(\frac{1}{\epsilon^2} (m + \log(m/\delta)) \right) \quad \text{and} \quad  s \geq O \left( \frac{1}{\epsilon}\log(m/\delta)^{5/2} + \log(m/\delta)^4 \right),
\end{equation}
which highlights a trade-off between the sampling size $\ell$ and the number $s$ of nonzero entries in each column. 
However, it has been experimentally observed that a constant parameter $s = 8$ and a sampling size $\ell = O(m \log(m))$ produce an embedding matrix $\Omega$ with parameter $\epsilon \approx 1/2$ with high probability for a wide variety of applications~\cite{Tropp2019streaminglowrank}. 

We finally outline the subsampled randomized Hadamard transform (SRHT), or SRHT distribution. Assuming that $n = 2^p$, to draw $\Omega \in \R^{\ell \times n}$ from the SRHT distribution we first draw a diagonal matrix $D \in \R^{n \times n}$ whose diagonal entries are signs $\pm 1$ drawn uniformly at random. We then apply the Walsh-Hadamard transform $H \in \R^{n \times n}$, defined by
\begin{align} \label{eq:hadamard} H_1 := \begin{bmatrix} 1 & 1 \\ 1 & -1 \end{bmatrix} \in \R^{2 \times 2}, \quad H_{j} := \begin{bmatrix} H_{j-1} & H_{j-1} \\ H_{j-1} & -H_{j-1} \end{bmatrix} \in \R^{2^j \times 2^j} \quad j\geq 2.\end{align}
\begin{align} \label{eq:hadamard2} H := \frac{1}{\sqrt{n}} \, H_p \in \R^{n \times n}.\end{align} 
When applied to a vector, this transformation uniformly distributes its mass across all its entries, with high probability \cite{Tropp2011improvedsrht}. We then draw $\ell$ rows of the identity matrix, uniformly at random and without replacement, to form $P \in \R^{\ell \times n}$ (applying $P$ to a vector is equivalent to sampling entries of this vector uniformly at random without replacement). A final scaling is required to compensate the sampling:
\[ \Omega = \sqrt{\frac{n}{\ell}} PHD \in \R^{\ell \times n}. \]
The matrix $H_p \in \R^{2^p \times 2^p}$ is a structured matrix, built recursively. For this reason, it can be applied without being formed by means of a fast recursive routine, such as the Fast-Walsh-Hadamard transform, which requires only $\mathcal{O}(n \log(n))$ flops. In the frequent case where $2^p < n < 2^{p+1}$, the input matrix can simply be padded with a block of zeros, so that it fits the application of $H_{p+1}$. 
The SRHT distribution is an ($\epsilon$,  $\delta$, $m$)-OSE for a sampling size $\ell$ such that \cite{BalabanovNouyI,BalabanovGrigori2022,woodruff2014sketching}
\begin{equation}
    \ell = O \left( \frac{1}{\epsilon^2} (\sqrt{m} + \sqrt{\log(n/\delta)})^2  \log(m/\delta) \right).
\end{equation}
Assuming that $m \gg \log(n)$, it is shown in~\cite{Tropp2011improvedsrht} that the sampling size $\ell$ can be set to $O(m \log(m))$, producing an embedding matrix $\Omega$ with parameter $\epsilon \approx 1/2$ with high probability for a wide variety of applications. The $\log(m)$ factor in the sampling size $\ell$ is necessary in the worst case: see, for instance, \cite[Remark~11.2]{HMT11}. 

%% file: randqr.tex
\section{Computation of a well-conditioned basis through randomization}
\label{sec:randomized-qr-introduction}

In this section, we outline the theoretical principles underlying the randomized orthogonalization framework. We then present three approaches for computing the randomized QR decomposition of a tall-and-skinny matrix $W$, namely the randomized Cholesky QR algorithm, the randomized Gram--Schmidt process, and the randomized Householder QR factorization.

\subsection{General discussion}
\label{subsec:sketch-orth-discussion}

Let $\mathcal{W} \subset \R^n$ be an $m$-dimensional subspace, and let $W \in \R^{n \times m}$ be a full-rank matrix such that $\range(W) = \mathcal{W}$. In many applications, the construction of an orthonormal basis of $\mathcal{W}$ is a key algorithmic component. For example, an orthonormal basis $Q$ of $\mathcal{W}$ is constructed when solving an overdetermined least squares problem with coefficient matrix $W$, or when $\mathcal{W}$ is a Krylov subspace employed in the solution of a linear system or eigenvalue problem. An orthonormal basis of $\mathcal{W}$ can be constructed via a Householder QR factorization or a Gram-Schmidt process. Both algorithms construct the factorization $W = QR$, where $Q \in \R^{n \times m}$ has orthonormal columns and $R \in \R^{m \times m}$ is upper triangular and have a computational cost of $\mathcal{O}(nm^2).$ 

Let us denote by $\mathcal{P}_{\mathcal{W}}$ the orthogonal projector to $\mathcal{W}$. We recall that this projector satisfies the following properties:
\begin{itemize}
	\item for any $\vec x \in \R^n$ we have $\mathcal{P}_{\mathcal{W}} \vec x = \argmin_{\vec w \in \mathcal{W}} \norm{\vec x - \vec w}$,
	\item we have $\mathcal{P}_{\mathcal{W}} = Q Q^T$, where $Q$ is an orthonormal basis of $\mathcal{W}$,
	\item we have $\proj{\subspace} = Z Z^+$ for an arbitrary basis $Z$ of $\mathcal{W}$, where $Z^+$ denotes the Moore-Penrose pseudoinverse of $Z$. 
\end{itemize}

\begin{figure}[tb]
	\begin{center}
	\begin{tikzpicture}[x={(-0.6cm,-0.6cm)}, y={(1.2cm,0cm)}, z={(0cm,1.2cm)}, font=\large]
		\coordinate (q1) at (1,0.1,0);
		\coordinate (q2) at (0,0.95, 0);
		\coordinate (q3) at (0,-0.1,1.1);
		
		\coordinate (Omq1) at (1,4,0);
		\coordinate (Omq2) at (0,5,0);
		\coordinate (Omq3) at (0,4,1);

		\coordinate (s1) at (0.1,4,0);
		\coordinate (s2) at (0,4.2,0);
		\coordinate (s3) at (0,4,0.2);
	
		\node at (0,1,1) [scale=1.2] {$\R^n$};
		\node at (0,5,1) [scale=1.2]{$\R^\ell$};
	
		\draw[-{Stealth[length=2mm, width=2mm]}, line width=1.2pt] (0,0,0) -- (q1) node[below]{$\vec q_1$};
		\draw[-{Stealth[length=2mm, width=2mm]}, line width=1.2pt] (0,0,0) -- (q2) node[below]{$\vec q_2$};
		\draw[-{Stealth[length=2mm, width=2mm]}, line width=1.2pt] (0,0,0) -- (q3) node[above]{$\vec q_3$};
		
		\draw[-{Stealth[length=2mm, width=2mm]}, line width=1.2pt] (0,4,0) -- (Omq1) node[below]{$\Omega \vec q_1$};
		\draw[-{Stealth[length=2mm, width=2mm]}, line width=1.2pt] (0,4,0) -- (Omq2) node[below]{$\Omega \vec q_2$};
		\draw[-{Stealth[length=2mm, width=2mm]}, line width=1.2pt] (0,4,0) -- (Omq3) node[above]{$\Omega \vec q_3$};
	
		\draw[-, color=red, thick] (s2) -- (0,4.2,0.2);
		\draw[-, color=red, thick] (0,4.2,0.2) -- (s3);
	\end{tikzpicture}
	\end{center}
	\caption{Sketch-orthogonal basis (left) and its orthogonal sketch (right).}
	\label{fig:sketchedorthogonalbasis}
	\end{figure}
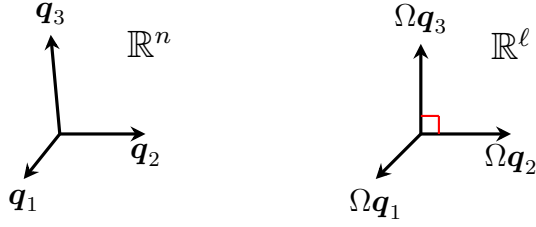

In this section, we introduce the concept of randomized QR factorization and present efficient algorithms for its computation. Assume that $\Omega \in \R^{\ell \times n}$ is an $\epsilon$-embedding for $\mathcal{W}$, and consider the decomposition: 
\begin{equation}
	\label{eq:randqr}
	W = QR, \qquad \Omega W =  \Omega Q \cdot R = S R, \quad (\Omega Q)^T \Omega Q = S^T S = I_m, \quad R \text{ upper triangular}.
\end{equation}
We refer to this decomposition as a \emph{randomized QR factorization} of $W$. Note that the columns of $Q = [\vec q_1 \cdots \vec q_m]$ are not orthogonal in general, but their sketches $\{\Omega \vec q_1, \dots, \Omega \vec q_m\}$ form an orthonormal basis of $\range(\Omega W)$, as illustrated in \cref{fig:sketchedorthogonalbasis}. Although $Q$ is not orthonormal, it is guaranteed to be extremely well-conditioned due to the $\epsilon$-embedding property of $\Omega$. Indeed, \cref{eqn:singular-values-eps-embedding} implies that we have
\begin{equation*}
	\cond(Q) \le \sqrt{\frac{1+\epsilon}{1-\epsilon}} \cdot \cond(\Omega Q) = \sqrt{\frac{1+\epsilon}{1-\epsilon}}.
\end{equation*}  
This property is fundamental to successfully applying the decomposition \cref{eq:randqr} to the solution of a variety of linear algebra problems, as we discuss in the following sections.
 
Before we present the algorithms to compute a randomized QR factorization, we introduce the \emph{sketch-orthogonal projector} $\skproj{\subspace}$ to $\subspace$. Given an embedding $\Omega$ of $\mathcal{W} = \range(W)$, given a randomized QR factorization~\eqref{eq:randqr}, we define:
\[ \forall \vec x \in \R^n, \quad \skproj{\subspace} \vec x = Q (\Omega Q)^T \Omega x. \] 
Using the sketch-orthogonality relation $(\Omega Q)^T \Omega Q = I$, we verify that $\skproj{\subspace}$ is indeed a projector. More precisely, it is an \textit{oblique} projector. We remark that it satisfies the identity
\begin{equation}
	\label{eqn:sketch-orth-projector-identity}
	\Omega \cdot \skproj{\subspace} = \proj{\Omega \subspace} \cdot \Omega,
\end{equation} 
where $\proj{\Omega \subspace} = \Omega Q (\Omega Q)^T$ is the orthogonal projector onto $\Omega \subspace \subset \R^\ell$. 
The sketch-orthogonal projector satisfies the following properties:
\begin{itemize}
	\item for any $\vec x \in \R^n$ we have $\skproj{\subspace} \vec x = \argmin_{\vec w \in \subspace} \norm{\Omega (\vec x - \vec w)}$,
	\item for any $\vec x \in \R^n$ we have $\vec x - \skproj{\subspace} \vec x \perp^\Omega \subspace$, where we use $\perp^\Omega$ to denote the sketch-orthogonality condition $\Omega (\vec x - \skproj{\subspace}\vec x) \perp \Omega \subspace$,
	\item given an arbitrary basis $Z$ of $\mathcal{W}$, the sketch-orthogonal projector can be equivalently written as $\skproj{\subspace} = Z (\Omega Z)^+ \Omega$.
\end{itemize}

\begin{figure}[t]
    \centering
\begin{tikzpicture}[tdplot_main_coords, scale=3]

    \filldraw[fill=blue!20, draw=blue!50!black, opacity=0.5]
      (-1.5,-0.5,0) -- (0.5,-0.5,0) -- (0.5,1.5,0) -- (-1.5,1.5,0) -- cycle;

    \node[blue!50!black, font=\Large] at (-0.8,-0.25,0.05) {$\mathcal W$};

    \coordinate (B) at (-0.6,0.9,1);
    \draw[->, thick] (0,0,0) -- (B) node[anchor=west] {$\vec b$};

    \coordinate (Porth) at (-0.6,0.9,0);

    \begin{scope}
        \draw[fill=blue!20, draw=blue!40, opacity=1]
            (Porth) circle (0.3);   
    \end{scope}

    \draw[dashed] (B) -- (Porth);
    \filldraw[black] (Porth) circle (0.6pt);
    \node[anchor=south west] at (Porth) {$\proj{\mathcal W} \vec b$};

    \coordinate (Psk) at (-0.45,0.79,0);
    \filldraw[blue!70!black] (Psk) circle (0.6pt);
    \node[anchor=south east, blue!70!black] at (Psk)
        {$\skproj{\mathcal W} \vec b$};

    \draw[dotted, thick, blue!70!black] (B) -- (Psk);

    \coordinate (R) at ($(Porth) + (0,0.3,0)$); 
    
    \draw (Porth) -- (R);
    
    \draw[<->] ($(Porth) + (0.05,0,-0.05)$) --
               ($(R) + (0.05,0,-0.05)$)
               node[anchor=north] {$\displaystyle O(\varepsilon) \cdot \min_{\vec w \in \mathcal{W}} \norm{\vec b - \vec w}$};
    
    \draw ($(Porth)$) --
          ($(Porth) + (0,0,0.08)$) --
          ($(Porth) + (0,0.08,0.08)$) --
          ($(Porth) + (0,0.08,0)$);



\end{tikzpicture}
    \caption{Quasi optimality of the sketched projection $\skproj{\mathcal{W}} \vec b$. }
	\label{fig:sketchedproj}
\end{figure}
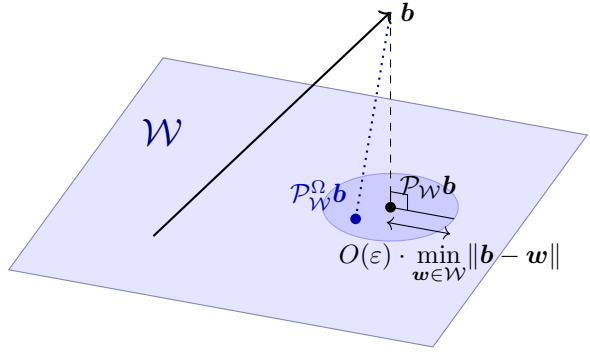

The sketch-orthogonal projector can be used as an approximation of the standard orthogonal projector \cite{BalabanovGrigori2022,BalabanovGrigori25block}.  Let $\vec b \in \R^n$, and assume that $\Omega$ is an $\epsilon$-embedding for $\mathcal{W} + \range(\vec b)$. Recalling the $\varepsilon$-embedding property and \cref{eqn:sketch-orth-projector-identity}, we have 
\begin{equation*}
	\norm{\vec b - \skproj{\subspace} \vec b} \leq \frac{1}{\sqrt{1-\epsilon}} \norm{\Omega \vec b - \Omega \skproj{\subspace} \vec b} = \frac{1}{\sqrt{1-\epsilon}} \norm{\Omega \vec b - \proj{\Omega \subspace} \Omega \vec b} = \frac{1}{\sqrt{1-\epsilon}} \min_{\vec z \in \Omega \subspace} \norm{\Omega \vec b - \vec z},
\end{equation*}
where for the last equality, we used the optimality property of the orthogonal projector $\proj{\Omega \subspace}$. 
For any fixed $\vec w \in \subspace$, we have $\vec z = \Omega \vec w \in \Omega \subspace$, so again using the $\varepsilon$-embedding property, we get
\begin{equation*}
	\norm{\vec b - \skproj{\subspace} \vec b} \le \frac{1}{\sqrt{1-\epsilon}} \norm{\Omega(\vec b - \vec w)} \le \sqrt{\frac{1+\epsilon}{1-\epsilon}} \norm{\vec b - \vec w}.
\end{equation*}   
By taking the minimum over $\vec w \in \mathcal{W}$, we finally obtain
\begin{equation}
	\label{eqn:sketch-projector-approx-least-squares}
	\norm{\vec b - \skproj{\subspace} \vec b} \le \sqrt{\frac{1+\epsilon}{1-\epsilon}} \cdot \min_{\vec w \in \subspace}\norm{\vec b - \vec w}.
\end{equation} 
This shows that the sketch-orthogonal projection of $\vec b$ onto $\subspace$ is a \emph{quasi-optimal minimizer} of the distance between $\vec b$ and $\subspace$, so the sketch-orthogonal projector $\skproj{\subspace}$ acts as an approximation of the orthogonal projector $\proj{\subspace}$. This property is illustrated in \cref{fig:sketchedproj}.

This property has a straightforward implication for the solution of a least squares problem. Indeed, let us denote by $\vec x^\star$ the solution of the sketched least squares problem
\begin{equation*}
	\vec x^\star = \argmin_{\vec x \in \R^m} \norm{\Omega W \vec x - \Omega \vec b}.
\end{equation*}
Then we have $W \vec x^\star = \skproj{\mathcal{W}} \vec b$, and from \cref{eqn:sketch-projector-approx-least-squares} it follows that $\vec x^\star$ satisfies 
\begin{equation*}
	\norm{W \vec x^\star - \vec b} \le \sqrt{\frac{1+\epsilon}{1-\epsilon}} \cdot \min_{\vec x \in \R^m} \norm{W \vec x - \vec b},
\end{equation*}
i.e., it is an approximate solution of the corresponding non-sketched least squares problem (see, e.g., \cite{sarlos2006improved} and \cite[Section~2.2]{NakatsukasaTropp24}). Furthermore, the randomized QR factorization in~\eqref{eq:randqr} yields a closed form formula for the computation of $\vec x^\star$:
\begin{equation*}
    W \vec x^\star = \skproj{\mathcal{W}} \vec b \quad \implies \quad QR \vec x^\star = Q (\Omega Q)^T \Omega \vec b \quad \implies \quad \vec x^\star = R^{-1} (\Omega Q)^T \Omega \vec b.
\end{equation*}

We conclude this introductory section on the randomized QR factorization by outlining a framework for computing it, which underpins the randomized QR processes that we present in the following sections. Given a full-rank matrix $W \in \R^{n \times m}$ and an $\epsilon$-embedding $\Omega \in \R^{\ell \times n}$ for $\mathcal{W} = \range(W)$, a randomized QR factorization of $W$ can be computed using the following simple procedure: 
\begin{enumerate}
    \item Compute the sketch $\Omega W \in \R^{\ell \times m}$.
	\item Compute a QR factorization $\Omega W = S R$.
	\item Set $Q = W R^{-1}$.
\end{enumerate}
Then it follows that $\Omega Q = S$ and $Q$, $S$ and $R$ satisfy \cref{eq:randqr}. The idea for this randomized orthogonalization procedure originates from \cite{Rokhlin2008fastrandomized}, where the $R$ factor is used as a preconditioner for the solution of a least squares problem, and is presented in \cite[Algorithm~2.3]{BalabanovGrigori25block}. This framework is sometimes simply called \textit{randomized QR}~\cite{Rokhlin2008fastrandomized,Higgins2025multisketching}, or sometimes \textit{randomized Cholesky QR}~\cite{BalabanovGrigori25block,Balabanov2022rcholesky}, or \textit{randomized preconditioning} or \textit{sketch-and-precondition}~\cite{Garrison2024rcholesky,Meier2024sketch}, and it is often used as a preconditioner for subsequent deterministic algorithms. For instance, in~\cite{Rokhlin2008fastrandomized}, where the authors propose to solve a least squares problem, the obtained factor $R$ and the solution $\vec x^\star$ to the sketched least squares are used, respectively, as a preconditioner and a starting point for conjugate gradient iterations. 

In the foundational work~\cite{Rokhlin2008fastrandomized}, and in many sketch-and-precondition papers~\cite{Melnichenko2025rcholesky,Garrison2024rcholesky,Meier2024sketch}, the authors obtain the triangular factor $R$ from $\Omega W$ by pivoted (strong) rank-revealing factorization, rather than a simple QR factorization. Moreover, the randomized QR algorithm can be followed by an efficient algorithm for the computation of a QR factorization, such as standard deterministic CholeskyQR, to obtain an orthogonal factor~$Q$~\cite{BalabanovGrigori25block}. For example, in~\cite{Melnichenko2025rcholesky,Garrison2024rcholesky,Meier2024sketch}, the matrix $W$ is preconditioned with the truncated, pivoted factor $R$ and this is followed by a standard deterministic CholeskyQR factorization. 
All of these algorithms are very efficient, since they require a constant number of synchronizations. In addition, randomized QR + CholeskyQR hybrids perform most of their flops through BLAS3 kernels.

Instead of applying $R^{-1}$ to $W$ in order to obtain $Q$ explicitly, one can also keep the basis $W$ and compute $Q \vec x$ as $W(R^{-1} \vec x)$, i.e., applying $R^{-1}$ to the input vectors instead. This approach, often called \emph{whitening of the basis}, is widely used in randomized Krylov subspace methods \cite{NakatsukasaTropp24,GuettelSchweitzer23,PSS25mf,PSS25me}. We refer to \cref{subsubsec:whitening} for further details.

To close this section, we emphasize that, as in standard orthogonalization processes, there are multiple methodologies for the sketch orthogonalization of $W$. Although they are all equivalent in exact arithmetic, they accumulate rounding errors in different ways when performed in floating point arithmetic, and they suffer from these errors in various ways.

%% file: rgs.tex
\subsection{Randomized Gram--Schmidt}

In this section we present the randomized Gram--Schmidt process \cite{BalabanovGrigori2022, BalabanovGrigori25block} for computing the randomized QR decomposition \cref{eq:randqr} of a tall-and-skinny matrix $W \in \R^{n \times m}$. 
This randomized process is inspired by the deterministic Gram--Schmidt process, which we briefly recall here. 

Let us denote by $\vec w_1, \dots \vec w_m$ the columns of $W$, and by $\mathcal{W}_j = \range([\vec w_1, \dots, \vec w_j])$. The Gram--Schmidt process constructs an orthonormal basis $Q = [\vec q_1, \dots, \vec q_m]$ of $\range(W)$ by iteratively subtracting from $\vec w_{j+1}$ its projection onto $\mathcal{W}_j$. More precisely, the algorithm sets $\vec q_1 = \vec w_1 / \norm{\vec w_1}$, and then for $j = 1, \dots, m-1$ we set 
\begin{equation}
	\label{eqn:standard-gram-schmidt}
	\widetilde{\vec q}_{j+1} = (I - \proj{\mathcal{W}_j}) \vec w_{j+1}, \qquad \vec q_{j+1} = \widetilde{\vec q}_{j+1} / \norm{\widetilde{\vec q}_{j+1}}.
\end{equation}  
The practical implementation of the Gram--Schmidt process depends on the specific implementation of the projector $\proj{\mathcal{W}_j}$. Letting $Q_j = [\vec q_1, \dots, \vec q_j]$, we have $\proj{\mathcal{W}_j} = Q_j Q_j^T$ and thus we can implement \cref{eqn:standard-gram-schmidt} as
\begin{equation*}
	\widetilde{\vec q}_{j+1} = (I - Q_j Q_j^T) \vec w_{j+1}, \qquad \vec q_{j+1} = \widetilde{\vec q}_{j+1} / \norm{\widetilde{\vec q}_{j+1}},
\end{equation*}
which corresponds to the \emph{classical Gram--Schmidt} (CGS) algorithm. The main advantage of the CGS implementation is that it performs the inner products $Q_j \vec w_{j+1}$ by exploiting matrix-vector BLAS2 routines. Using the orthogonality of the columns of $Q_j$, we have $I - Q_j Q_j^T = \prod_{k = 1}^{j} (I - \vec q_k \vec q_k^T)$, so we can also write \cref{eqn:standard-gram-schmidt} as
\begin{equation*}
	\widetilde{\vec q}_{j+1} = \prod_{k = 1}^{j} (I - \vec q_k \vec q_k^T) \vec w_{j+1}, \qquad \vec q_{j+1} = \widetilde{\vec q}_{j+1} / \norm{\widetilde{\vec q}_{j+1}},
\end{equation*}
which corresponds to the \emph{modified Gram--Schmidt} (MGS) algorithm. This algorithm computes the inner products with the columns of $Q_j$ sequentially, so it is slower than CGS on modern computational architectures, but it has better numerical stability. In general, the stability of CGS and MGS can be improved by applying the projector $I - \proj{\mathcal{W}_j}$  twice, leading to the CGS2 and MGS2 algorithms. The numerical stability of different implementations of the Gram--Schmidt process and its relation with the condition number of $W$ has been extensively studied in the literature, see, e.g., \cite{CLRT22}.

The randomized Gram--Schmidt process essentially replaces the orthogonal projector $\proj{\mathcal{W}_j}$ in \cref{eqn:standard-gram-schmidt} with the oblique projector $\skproj{\mathcal{W}_j}$ to construct a basis $Q$ that is now sketch-orthogonal. The first column of $Q$ is set as $\vec q_1 = \vec w_1 / \norm{\vec w_1}$, and for $j = 1, \dots, m-1$ we compute
\begin{equation}
	\label{eqn:randomized-gram-schmidt}
	\widetilde{\vec q}_{j+1} = (I - \skproj{\mathcal{W}_j}) \vec w_{j+1}, \qquad \vec q_{j+1} = \widetilde{\vec q}_{j+1} / \norm{\Omega \widetilde{\vec q}_{j+1}}.
\end{equation}
It immediately follows from the properties of the sketch-orthogonal projector $\skproj{\mathcal{W}_j}$ that $\vec q_{j+1} \perp^\Omega \mathcal{W}_j$, which in turn implies that $\Omega Q = [\Omega \vec q_1, \dots, \Omega \vec q_m]$ is an orthonormal basis of $\range(\Omega W)$. Recall that since $Q_j = [\vec q_1, \dots, \vec q_j]$ is a sketch-orthogonal basis of $\mathcal{W}_j$, we have $\skproj{\mathcal{W}_j} = Q_j (\Omega Q_j)^T \Omega = Q_j (\Omega Q_j)^+ \Omega$ and we can more explicitly write
\begin{equation}
	\label{eqn:randomized-gram-schmidt-least-squares}
	\widetilde{\vec q}_{j+1} = \vec w_{j+1} - Q_j \vec h_j, \qquad \vec h_j = (\Omega Q_j)^+ \Omega \vec w_{j+1} = \argmin_{h \in \R^j} \norm{\Omega Q_j \vec h - \Omega \vec w_{j+1}}.
\end{equation}
From \cref{eqn:randomized-gram-schmidt-least-squares} we see that the main difference between the standard and randomized Gram--Schmidt processes is that the latter replaces inner products of vectors of length $n$ with inner products of much shorter sketched vectors of length $\ell$ or the solution of an $\ell \times j$ least squares problem, drastically reducing the cost of this operation. Since the product $Q_j \vec h_j$ still needs to be performed, the overall computational cost is roughly half that of the deterministic Gram--Schmidt process.

The solution of the least squares problem for the computation of $\vec h_j$ in \cref{eqn:randomized-gram-schmidt-least-squares} is crucial for the implementation of the randomized Gram--Schmidt process. By exploiting the fact that $\Omega \vec Q_j$ is sketch-orthogonal, we can simply compute $\vec h_j = (\Omega Q_j)^T \Omega \vec w_{j+1}$ and obtain the \emph{randomized classical Gram--Schmidt} algorithm. However, since $\Omega Q_j$ is a small $\ell \times j$ matrix, we can afford to solve the least squares problem with a more expensive method to achieve better numerical stability, without assuming that $\Omega Q_j$ has orthonormal columns, which in general is not true in finite-precision arithmetic. We refer to \cite[Section~2.3]{BalabanovGrigori2022} for further details. We also mention that this process can be combined with deterministic reorthogonalization to obtain a basis with orthonormal columns; see for instance \cite[Section~3.1]{JangGrigori25}, where an orthogonal projector is obtained as a combination of the randomized Gram--Schmidt projector and either CGS or MGS.

The randomized Gram--Schmidt process (RGS) is presented in \cref{algo:rgs}. As detailed in~\cite{BalabanovGrigori2022}, its flop cost is dominated by $nm^2 + 2mt$ flops, where $t$ is the flop cost of sketching one vector. With SRHT, we thus get $nm^2 + 2nm \log(n)$ flops, namely half the flops of Gram--Schmidt processes. Depending on the availability of the vectors $\vec w_1, \hdots, \vec w_m$ during factorization, this algorithm requires between $1$ and $2$ synchronizations per iteration, similar to the cost of communication of CGS. As in CGS, most of the flops between synchronizations in RGS can be carried out by BLAS2 routines.

\begin{algorithm}[t]
\caption{Randomized Gram-Schmidt process}
\label{algo:rgs}
\begin{algorithmic}[1]
\Require $W \in \R^{n \times m}$ full-rank, $\Omega \in \R^{\ell \times n}$ that is an $\epsilon$-embedding for $\range(W)$
\Ensure $Q \in \R^{n \times m}$, $S \in \R^{\ell \times m}$ and $R \in \R^{m \times m}$ such that $W=QR$, $S = \Omega Q$,  $S^T S = I_m$, $R$ upper triangular
\Function{Randomized-Gram-Schmidt}{$W, \Omega$}
    \State $\vec z \gets \Omega \vec w_1$
    \State $R_1 \gets [\|\vec z\|]$, $Q_1 \gets [\vec w_1 / \|\vec z\|]$, $S_1 \gets [\vec z / \|\vec z\|]$
    \For{$j = 2:m$}
        \State $\vec z \gets \Omega \vec w_j$ \label{algo:rgs:sketch1}
        \State $\vec r \gets S_{j-1}^+ \vec z$ \label{algo:rgs:lls} \# use a stable method to solve the least squares problem
        \State $\vec w \gets \vec w_j - Q_{j-1} \vec r$ \label{algo:rgs:refresh}
        \State $\vec z \gets \Omega \vec w$ \label{algo:rgs:sketch2}
        \State $R_j \gets \left[ \begin{array}{ccc|c} & & & \\ & R_{j-1} &  & \vec r \\ & & & \\\hline  & 0_{1\times (j-1)} & &  \|\vec z\| \end{array} \right]$, $\quad Q_j \gets [Q_{j-1} \; | \; \vec w / \|\vec z\| ] $, $\quad S_j \gets [S_{j-1} \; | \; \vec z / \|\vec z\|]$ \label{algo:rgs:scalew}
    \EndFor
    \State \Return $Q = Q_m$, $S = S_m$, $R = R_m$
\EndFunction
\end{algorithmic}
\end{algorithm}

We emphasize that the sketch in line~\ref{algo:rgs:sketch2} is crucial for numerical stability, as shown in the finite-precision analysis in~\cite{BalabanovGrigori2022}. An algorithm that replaces this sketch with a formula inferring $\Omega \vec w$ from $\vec z - S_{j-1} \vec r \in \R^\ell$ would not qualify as an implementation of RGS, but rather as an implementation of the \textit{sketch-and-precondition} framework. We also emphasize that line~\ref{algo:rgs:lls} is specified by which orthogonalization method is chosen by the user to orthogonalize $\Omega W$. It is crucial to select one that is stable enough to handle the successful orthogonalization of $\Omega W$. 

\begin{figure}[t]
    \centering
    \begin{subfigure}{0.48\textwidth}
        \centering
        \includegraphics[width=\linewidth]{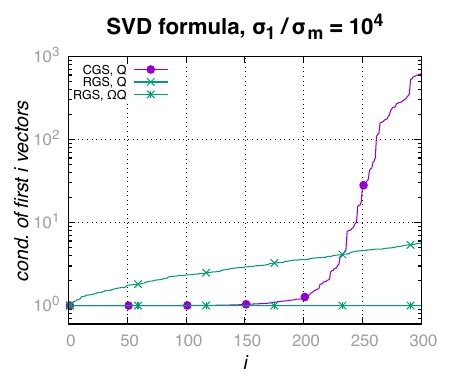}
        \caption{Condition number of basis $Q$ and \textit{a posteriori} sketch $\Omega Q$ built by RGS, compared with CGS.}
        \label{fig:condseasy}
    \end{subfigure}\hfill
    \begin{subfigure}{0.48\textwidth}
        \centering
        \includegraphics[width=\linewidth]{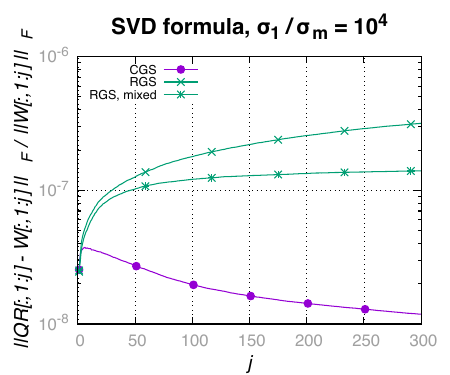}
        \caption{Relative Frobenius factorization errors of RGS in both single and mixed precision, compared with CGS.}
        \label{fig:errseasy}
    \end{subfigure}
    \caption{Comparison of CGS and RGS on a medium difficulty example.}
    \label{fig:rgscgseasy}
\end{figure}

We illustrate in~\Cref{fig:rgscgseasy} the numerical efficiency of RGS compared to CGS on a medium difficulty example. The input matrix $W$ is initialized in double precision through an SVD formula, with singular values decreasing exponentially from $10^2$ to $10^{-2}$. CGS is tested in single precision. RGS is tested in single precision, and in mixed precision (with low-dimensional operations done in double precision and remaining operations in single precision), with the results cast into single precision. We see in~\Cref{fig:condseasy} that the basis obtained with CGS loses orthogonality slowly in the first iterations, and then much more quickly, near the point where $W^T W$ becomes numerically singular. At the same time, the sketch of the basis generated by RGS is numerically orthogonal, which in turn explains the small condition number of the resulting basis. As illustrated in~\Cref{fig:errseasy}, all factorization errors remain very small, with CGS performing marginally better. Employing mixed precision for RGS further enhances the factorization accuracy.

%% file: householder.tex
\subsection{Randomized Householder QR}

We now introduce a randomized version of the celebrated Householder QR factorization. We only outline here the elements that are directly used in this factorization. More general properties can be found in~\cite{rhqr}. We give first a brief summary of the standard Householder QR.

\begin{figure}[t]
    \centering
    \begin{tikzpicture}
    \fill[purple] (0,0) rectangle ++(1,2);
    \draw[-] (0,0) -- (1,0) -- (1,2) -- (0,2) -- cycle;

    \draw[->] (1.2,1) -- (2.3,1);
    \node at (1.8,1.4) {$P_1$};
    
    \fill[blue] (2.5,1.8) rectangle ++(1,0.2);
    \fill[red] (2.7,0) rectangle ++(0.8,1.8);
    \draw[-] (2.5,0) -- (3.5,0) -- (3.5,2) -- (2.5,2) -- cycle;

    \draw[->] (3.7,1) -- (4.8,1);
    \node at (4.2,1.4) {$P_2$};

    \fill[blue] (5,1.8) rectangle ++(1,0.2);
    \fill[blue] (5.2,1.6) rectangle ++(0.8,0.2);
    \fill[orange] (5.4,0) rectangle ++(0.6,1.6);
    \draw[-] (5,0) -- (6,0) -- (6,2) -- (5,2) -- cycle;

    \draw[->] (6.2,1) -- (6.8,1);
    \node at (7.2,1) {$\ldots$};
    \draw[->] (7.7,1) -- (8.3,1);
    
    \fill[blue] (8.8,1.8) rectangle ++(1,0.2);
    \fill[blue] (9,1.6) rectangle ++(0.8,0.2);
    \fill[blue] (9.2,1.4) rectangle ++(0.6,0.2);
    \fill[blue] (9.4,1.2) rectangle ++(0.4,0.2);
    \fill[blue] (9.6,1) rectangle ++(0.2,0.2);
    \draw[-] (8.8,0) -- (9.8,0) -- (9.8,2) -- (8.8,2) -- cycle;

    \node at (9.5,1.7) {$R$};
    \node at (0.5,1.0) {$W$};
\end{tikzpicture}
\caption{Triangularizing $W$ through Householder reflections.}
\label{fig:hqrortho}
\end{figure}
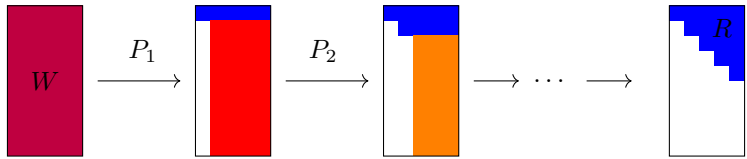

The Householder QR is an orthogonalization process alternative to the Gram--Schmidt process. The central operator of the process is the \textit{Householder reflector}:
\[ P = I_n - \beta \vec u \vec u^T, \quad \vec u \in \R^n \setminus \{0\}, \quad \beta = 2/\|\vec u\|^2.\] 
It is an orthogonal reflector, i.e., it verifies $P^T P = P^2 = P P^T = I_n$. The Householder process is derived from the ability to easily generate a Householder reflector $P$ that annihilates all the entries in a given vector $\vec w$ below a given index. Indeed, for some $\vec w_j \in \R^n$, denoting by $\vec w'$ the vector formed by $j-1$ zeros followed by the last $n-j+1$ entries in $\vec w_j$, we may define
\begin{align} \label{eq:hvector} \rho_j := \| \vec w' \|, \quad \sigma_j := \mathrm{sign}(\vec e_j^T \vec w'), \quad \vec u_j := \vec w' + \sigma_j \|\vec w'\| \vec e_j,  \quad \beta_j := 2/\|\vec u_j\|^2 \end{align}
and verify that $P_j = I_n - \beta_j \vec u_j \vec u_j^T$ annihilates all entries of $\vec w_j$ strictly below the $j$-th index, while not modifying the first $j-1$ entries of any vector. We can thus triangularize an arbitrary matrix $W \in \R^{n \times m}$ with Householder reflectors: we generate $P_1$ that annihilates the first column $\vec w_1$ below the first index, and apply it to the whole matrix; then we generate $P_2$ that annihilates the updated second column $P_1 \vec w_2$ below the second index and does not modify the first row of $P_1 W$, and apply it to the whole matrix, and continue similarly on the following columns, as illustrated in~\Cref{fig:hqrortho}. We obtain:
\[ P_m P_{m-1} \cdots P_1 W = \begin{bmatrix}
    R \\ 0_{(n-m)\times m}
\end{bmatrix} \implies W = P_1 \cdots P_m \begin{bmatrix}
    R \\ 0_{(n-m)\times m}
\end{bmatrix}.\]
As shown in~\cite{Schreiber1989storage}, the reflectors can be aggregated as follows:
\[ W = (I_n - U T U^T)\begin{bmatrix}
    R \\ 0_{(n-m)\times m}
\end{bmatrix}, \quad U = [\vec u_1 \; \cdots \; \vec u_m] \in \R^{n \times m}, \quad T \in \R^{m \times m} \text{ upper triangular. }\]

\begin{figure}[t]
\centering
\begin{minipage}[b]{0.45\textwidth}
    \hspace{2.1cm}\begin{tikzpicture}[scale=1.6]
        \coordinate (A) at (0.5,0.866);
        \coordinate (B) at (0.575, 0.433);
        \coordinate (C) at (0.65,0);
        \coordinate (D) at (-0.08,0.433);
        \coordinate (E) at (0,0);
    
        \draw[dashed, blue] (-0.5,-0.35) -- (1.4722,1.05) node[right, above] {\footnotesize$\{ \vec x \, : \, (\Psi \vec u)^T \Psi \vec x = 0\}$};
        \draw[dashed] (-0.08,0.433) -- (0.5,0.866) -- (0.65,0);
        \draw[-{Stealth[length=3mm, width=2mm]}, thick,blue] (0,0) -- (-0.15,0.866) node[above] {$\vec u$};
        \draw[-{Stealth[length=3mm, width=2mm]}, thick] (0,0) -- (0.65,0) node[below=2mm] {$P \vec w$};
        \tkzMarkSegment[color=red,pos=.58,mark=|](A,B)
        \tkzMarkSegment[color=red,pos=.5,mark=|](B,C)
        \tkzMarkSegment[color=red,pos=.45,mark=|](D,E)

        \draw[-{Stealth[length=3mm, width=2mm]}, thick, thick] (0,0) -- (0.5,0.866) node[right] {$\vec w$};
        
    \end{tikzpicture}
\end{minipage}
\hfill
\begin{minipage}[b]{0.45\textwidth}
    \hspace{0.5cm}\begin{tikzpicture}[scale=1.6]

    \coordinate (A) at (-0.25,0.433);
    \coordinate (B) at (0.5,0.866);
    \coordinate (C) at (0.75,0.433);
    \coordinate (D) at (1,0);
    \coordinate (E) at (0,0);
        \draw[dashed, blue] (-0.5,-0.29) -- (1.4722,0.85) node[below, right] {$(\Psi \vec u)^\perp$};
        \draw[-{Stealth[length=3mm, width=2mm]}, thick,blue] (0,0) -- (-0.5,0.866) node[above] {$\Psi \vec u$};
        \draw[-{Stealth[length=3mm, width=2mm]}, thick] (0,0) -- (1,0) node[below=2mm] {$\Psi P \vec w$};
        \draw[dashed] (-0.25,0.433) -- (0.5,0.866) -- (1,0);

        \draw[-{Stealth[length=3mm, width=2mm]}, thick, thick] (0,0) -- (0.5,0.866) node[right] {$\Psi \vec w$};
        \draw[blue, line width=0.8pt] (-0.05,0.0866) -- (0.0366, 0.1366) -- (0.0866, 0.05);

        \tkzMarkSegment[color=red,pos=.5,mark=|](E,A)
        \tkzMarkSegment[color=red,pos=.5,mark=|](B,C)
        \tkzMarkSegment[color=red,pos=.5,mark=|](C,D)
        
    \end{tikzpicture}
\end{minipage}
\caption{Randomized Householder reflector $P \in \R^{n \times n}$ (left) and induced Householder reflector $P' \in \R^{(\ell+m)\times(\ell+m)}$ on the sketched space (right).} \label{fig:reflectors}
\end{figure}
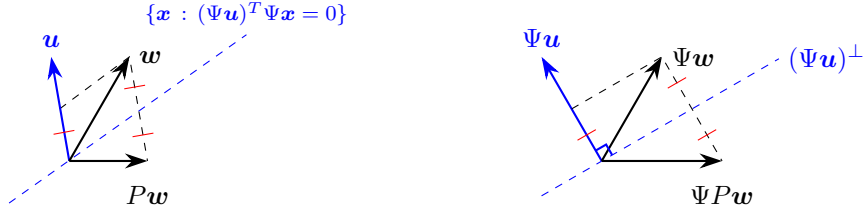

Let us now outline how this process can be randomized in order to obtain a randomized QR factorization equivalent to \cref{eq:randqr}. We use the following wrapper $\Psi$ of any sketching  matrix $\Omega$:
\begin{align}\label{eq:psi} 
\Psi = \left[ \begin{array}{c|ccccc} I_m & & & & & \\ \hline & & & & & \\ & & & \Omega & & \\ & & & & &\end{array} \right] \in \R^{(\ell+m) \times n}, \quad \Omega \in \R^{\ell \times (n-m)}.
\end{align}
Computing the sketch of a vector $\vec x \in \R^n$ with $\Psi$ thus consists of sketching the last $n-m$ coordinates of $\vec x$ with $\Omega \in \R^{\ell \times (n-m)}$, and concatenating the result to the first $m$ coordinates of $\vec x$. For simplicity, we denote $\ell' = \ell+m$ throughout this section. We then define the randomized Householder reflector associated with $\Psi$ as any matrix of the form
\[ P = I_n - \beta \, \vec u \, (\Psi \vec u)^T \Psi \in \R^{n \times n}, \quad \vec u \in \R^n \setminus \Ker{\Psi}, \quad \beta = 2/\|\Psi \vec u\|^2.  \]
We can verify that $P \in \R^{n \times n}$ defined in this way verifies $P^2 = I_n \neq P^T P$, i.e., it is an oblique reflection, with respect to the hyperplane $\{\vec x \in \R^n \, : \,  (\Psi \vec u)^T \Psi \vec x = 0\}$, and with $\vec u$ being mapped to $- \vec u$. We also verify that $ \Psi P = P' \Psi$ for some standard Householder reflector $P' \in \R^{\ell'\times \ell'}$. We illustrate the two related reflectors in~\Cref{fig:reflectors}. 

For any vector $\vec w_j \in \R^n$, $j \leq m$, using the same notation for $\vec w'$, we can straightforwardly \textit{randomize} the formulas in~\eqref{eq:hvector}: 
\begin{align} 
\label{eq:rhvector} 
\rho_j := \| \Psi \vec w' \|, \quad \sigma_j := \mathrm{sign}(\vec e_j^T \vec w'), \quad \vec u_j := \vec w' + \sigma_j \|\Psi \vec w'\| \vec e_j,  \quad \beta_j := 2/\|\Psi \vec u_j\|^2. 
\end{align}
The associated randomized Householder reflector $P_j = I_n - \beta_j \vec u_j (\Psi \vec u_j)^T \Psi$ also annihilates the last $n-j$ coordinates of $\vec w$, while not modifying the first $j-1$ entries of any vector. We can thus proceed as in the standard Householder process and generate randomized Householder reflectors $P_1 \hdots P_m \in \R^{n \times n}$ verifying
\[ W = P_1 \cdots P_m \begin{bmatrix}
    R \\ 0_{(n-m)\times m}
\end{bmatrix} = \left(I_n - U T (\Psi U)^T \Psi\right) \begin{bmatrix}
    R \\ 0_{(n-m)\times m}
\end{bmatrix} \in \R^{n\times n}, \] 
\[ \Psi W = \Psi \bigl( \underbrace{I_n - U T (\Psi U)^T \Psi \bigr) \begin{bmatrix}
    R \\ 0_{(n-m)\times m}
\end{bmatrix}}_{\text{RHQR fact. of W}} = \underbrace{ \bigl( I_{\ell'} - \Psi U \, T (\Psi U)^T \bigr) \begin{bmatrix}
    R \\ 0_{\ell \times m}
\end{bmatrix}}_{\text{HQR fact. of $\Psi W$}}  \in \R^{\ell' \times m}.\]

All these elements yield~\Cref{algo:rhqr_leftlooking} ($\mathrm{RHQR}$). As detailed in~\cite{rhqr}, this algorithm has the same computational and communication cost as RGS, while also leveraging mainly BLAS2 routines between sketches. Compared with the original Householder process, the randomized paradigm allows to aggregate the reflectors without synchronizations and for a negligible arithmetic cost. RHQR with SRHT is thus twice as cheap as non-aggregated Householder QR, and thrice as cheap as aggregated Householder QR.

\begin{algorithm}[t]
\caption{Randomized Householder QR (left-looking)}
\label{algo:rhqr_leftlooking}
\begin{algorithmic}[1]
\Require Matrix $W = \left[\vec w_1 \; | \cdots \; | \; \vec w_m \right]\in \R^{n \times m}$, $\Omega \in \R^{\ell \times (n-m)}$, $m < \ell \ll n-m$
\Ensure $U \in \R^{n \times m},S \in \R^{(\ell+m) \times m}, \, T,R \in \R^{m \times m}$ such that $S = \Psi U$ and $W = (I-UT (\Psi U)^T \Psi) \cdot \left[R; \; 0_{(n-m)\times m} \right]$
\Function{RHQR}{$W, \Omega$}
    \State $\vec z \gets \Psi \vec w_1$
    \State Define $\rho_1, \sigma_1, \vec u_1, \; \vec s_1 = \Psi \vec u_1, \beta_1$ as in equations~\eqref{eq:rhvector}
    \State $U_1 \gets \left[ \vec u_1 \right]$, $\; S_1 \gets \left[ \vec s_1\right]$, $\; T_1 \gets \left[\beta_1\right]$, $\; R_1 \gets \left[-\sigma_1 \rho_1\right]$
    \For{$j = 2:m$}
        \State $\vec z \gets \Psi \vec w_j$ \label{algo:leftrhqr:firstsketch}
        \State $\vec w \gets \vec w_j - U_{j-1} T_{j-1}^T S_{j-1}^T \vec z$  \label{algo:leftrhqr:refreshwj} 
        \State $\vec z \gets \vec z_j - S_{j-1} T_{j-1}^T S_{j-1}^T \vec z$  \label{algo:leftrhqr:refreshzj} 
        \State $\vec z \gets \Psi \vec w'$ \label{algo:leftrhqr:secondsketch}
        \State Define $\rho_j, \sigma_j, \vec u_j, \; \vec s_j = \Psi \vec u_j, \beta_j$ as in equations~\eqref{eq:rhvector}
        \State $U_j \gets \left[ U_{j-1} \; | \; \vec u_j \right]$, $\; S_j \gets \left[ S_{j-1} \; | \; \vec s_j \right]$, $\; R_j \gets \begin{bmatrix} R_{j-1} & (\vec z)_{1:j-1} \\ 0_{1 \times (j-1)} & -\sigma_j \rho_j \end{bmatrix}$, \label{algo:leftrhqr:updatecompact1}
        
        $\; \; \; T_j \gets \begin{bmatrix} T_{j-1} & - \beta_j T_{j-1} S_{j-1}^T \vec s_j \\ 0_{1 \times (j-1)} & \beta_j \end{bmatrix}$
    \EndFor
    \State \Return $R_m, \; U_m, \; S_m, \; T_m$
\EndFunction
\end{algorithmic}
\end{algorithm}

\begin{figure}[t]
    \centering
    \begin{subfigure}{0.48\textwidth}
        \centering
        \includegraphics[width=\linewidth]{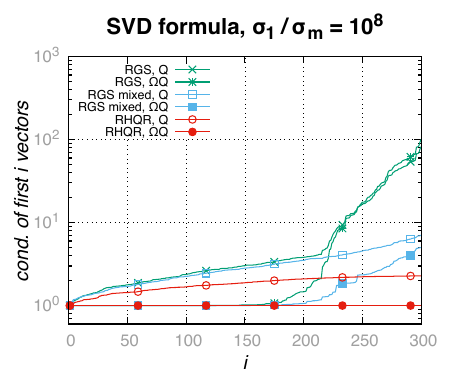}
        \caption{Condition number of basis $Q$ and \textit{a posteriori} sketch $\Omega Q$ built by RHQR, compared with RGS.}
        \label{fig:condsrhqr}
    \end{subfigure}\hfill
    \begin{subfigure}{0.48\textwidth}
        \centering
        \includegraphics[width=\linewidth]{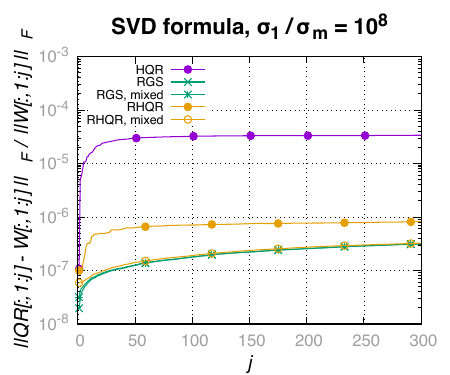}
        \caption{Relative Frobenius factorization errors of RHQR in both single and mixed precisions, compared with RGS and Householder QR.}
        \label{fig:errsrhqr}
    \end{subfigure}
    \caption{Comparison of RHQR and RGS on a difficult example.}
    \label{fig:rhqr}
\end{figure}

We illustrate the numerical stability of RHQR in~\Cref{fig:rhqr} for a difficult example. The input matrix $W$ is initialized in double precision with an SVD formula, with its singular values decreasing exponentially from $10^4$ to $10^{-4}$. It is then cast in single precision. Householder QR is tested in single precision. RGS and RHQR are tested in single and in mixed precisions (with low-dimensional operations done in double precision). We see in~\Cref{fig:condsrhqr} that the sketch of the basis computed by RGS loses orthogonality. When using mixed precision, this phenomenon occurs later. The loss of orthogonality in RGS, and the growth of the condition number of the basis, are very well mitigated by the use of mixed precision. Meanwhile, RHQR maintains the orthogonality of the sketch of the basis, which explains the small condition number of the basis itself. The use of mixed precision in RHQR does not substantially improve the condition number of the basis, as numerical sketched orthogonality is already achieved with single precision. The basis obtained with Householder QR also achieves numerical orthogonality. We see in~\Cref{fig:errsrhqr} that all factorization errors of randomized algorithms are small, with a noticeable advantage when compared to Householder QR. In single precision, the factorization error of RGS is slightly better than that of RHQR. The use of mixed precision allows RHQR to attain the same factorization error as RGS in both precision settings.

%% file: randorth-extra.tex
\subsection{Block sketch-orthogonalization}

When computing multiple matrix--vector products, substantial speedups can be achieved by replacing successive BLAS-2 operations $\vec c_1 = A \vec b_1$, $\vec c_2 = A \vec b_2, \hdots, \vec c_b = A \vec b_b$ with a single BLAS-3 matrix--matrix multiplication $C = AB$, where $B = [\vec b_1, \dots, \vec b_b]$. Indeed, this allows to reduce data movement between different levels of the memory hierarchy. BLAS-3 kernels can be exploited by using block algorithms that partition the input matrix $W \in \R^{n \times bm}$ into blocks $W_1, \dots, W_m \in \R^{n \times b}$.  Such a  block strategy is used in LAPACK's \texttt{xgeqrf} routines for computing the QR decomposition of dense matrices.  Rather than performing a single iterative loop with $bm$ orthogonalization steps for $W \in \mathbb{R}^{n \times bm}$, the computation is organized into two nested loops: an outer loop over the $m$ blocks and an inner loop over the $b$ vectors within each block. At step $jb+1$ (first vector of the $(j+1)$-th block), having built the matrix $\mathbf{Q}_j = [\vec q_1 \; \vec q_2 \; \cdots \; \vec q_{bj}] \in \R^{n \times bj}$ and its sketch $\Omega \mathbf{Q}_j \in \R^{\ell \times bj}$, assuming that the whole block $W_{j+1} \in \R^{n \times m}$ is available, we orthogonalize all vectors $W_{j+1} \vec e_1, \hdots , W_{j+1} \vec e_m$ with matrix-matrix operations only (BLAS3):
\begin{enumerate}
    \item Sketch $Z_{j+1} \gets \Omega W_{j+1}$.
    \item Orthogonalize against current basis $W_{j+1}' \gets W_{j+1} - \mathbf{Q}_j (\Omega \mathbf{Q}_j)^T Z_{j+1}$.
    \item Orthogonalize $W_{j+1}'$ with a single loop, BLAS2 algorithm.
\end{enumerate}
This results in the randomized block Gram-Schmidt process described in~\cite{BalabanovGrigori25block}. The RHQR process, thanks to its cost-free aggregation of reflectors, can be easily expressed as a block algorithm~\cite{rhqr}. We emphasize that these block algorithms are mathematically equivalent to single-loop BLAS2 algorithms and performs the same number of flops. The only difference is that more flops are performed through BLAS3 kernels.

The sketch $Z_{j+1}$ is necessary in all randomized algorithms (the input matrix $W_{j+1}$ must be sketched at least once). Since $Z_{j+1}$ is a matrix of small dimensions that approximately preserves the condition number of $W_{j+1}$, this condition number can thus be efficiently estimated. If it is small, less stable but faster algorithms can be used for the orthogonalization of the block, leading to even more significant speedups.

\subsection{Bi-orthogonalization}

In \cite{GPS25}, a randomized two-sided Gram--Schmidt algorithm is introduced to compute sketch-biorthogonal bases associated with two subspaces $\mathcal{X}$ and $\mathcal{Y}$ of the same dimension. This algorithm computes two bases $P$ and $Q$ such that $\range(Q) = \mathcal{X}$ and $\range(P) = \mathcal{Y}$, satisfying the sketch-biorthogonality condition $(\Omega P)^T \Omega Q = I$. This approach is computationally cheaper and more numerically stable than the two-sided Gram--Schmidt deterministic process, and it often constructs bases that are better conditioned than those obtained by deterministic algorithms which impose the biorthogonality condition $P^T Q = I$. We refer to \cite{GPS25} for further details.

\subsection{Computation in mixed precision and on parallel computers} 
\begin{figure}[t]
\centering
\begin{tikzpicture}[>=Latex, font=\sffamily, scale = 1]

\foreach \i/\y in {1/3,2/2,3/1,4/0} {
  \node[anchor=east] at (-2.0,\y) {Process \i};
}

\foreach \y in {0.5,1.5,2.5} {
  \draw[dotted, thick] (-3.2,\y) -- (8.6,\y);
}

\foreach \i/\y in {1/3,2/2,3/1,4/0} {
  \draw[fill=colO, rounded corners=2pt]
      (-1.3,\y-0.15) rectangle (-0.4,\y+0.45);
  \node[font=\itshape\large]
      at (-0.85,\y+0.15) {$\Omega_{\i}$};
}

\foreach \i/\y in {1/3,2/2,3/1,4/0} {
  \draw[fill=colW, rounded corners=2pt]
      (-0.3,\y-0.45) rectangle (0.3,\y+0.45);
  \node[text=black, font=\itshape\large]
      at (0,\y) {$W_{\i}$};
}

\node[circle, fill=darkteal, minimum size=8pt] (c) at (1.5,1.5) {};

\foreach \y in {3,2,1,0} {
  \draw[thick, darkteal]
      (0.5,\y) .. controls (0.8,\y) and (1.5,1.5) .. (c);
  \draw[thick, darkteal]
      (c) .. controls (1.5,1.5) and (2.2,\y) .. (2.5,\y);
}

\foreach \i/\y in {1/3,2/2,3/1,4/0} {
  \node[text=black, font=\itshape\large]
      at (3,\y) {$\Omega W$};
}

\foreach \i/\y in {1/3,2/2,3/1,4/0} {
  \draw[->, thick] (3.7,\y) -- (4.6,\y);
  \node[font=\large]
      at (5.1,\y) {$R,$};
}

\foreach \i/\y in {1/3,2/2,3/1,4/0} {
  \draw[fill=colQ, rounded corners=2pt]
      (5.8,\y-0.45) rectangle (6.4,\y+0.45);
  \node[text=black, font=\itshape\large]
      at (6.1,\y) {$Q_{\i}$};
  \node[text=black, font=\itshape\large]
      at (7.4,\y) {$= W_{\i} R^{-1}$};
}

\node[text=darkteal] at (1.5,3.7) {Sketch};
\node                at (4.2,3.7) {Stable QR};

\end{tikzpicture}
    \caption{RCholeskyQR on 4 processes.}
     \label{fig:parallelrcholesky}
\end{figure}
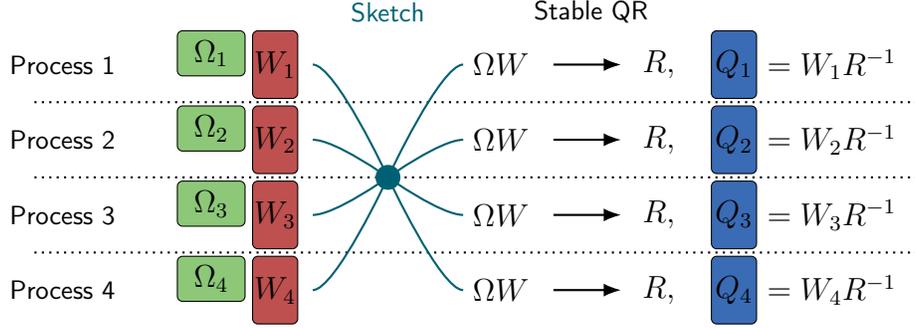

Randomized algorithms benefit not only from optimized kernels, but also from mixed precision and reduced communication on parallel architectures. On a parallel computer, the matrix $W \in \mathbb{R}^{n \times m}$ is typically distributed over processes by using a block row distribution, as in~\cite{Demmel2012tsqr}.  Sketching $W$ can be performed efficiently in parallel. Consider for example $\Omega$ to be a dense Gaussian or an s-hashing sketching. As displayed in Figure~\ref{fig:parallelsketching}, the sketching matrix is partitioned into blocks of columns and can be generated locally on each process with no communication. Once each process computes a local sketch $\Omega_i W_i$, an Allreduce communication among processes is required to sum the local sketches and compute $\Omega W = \sum_{i=1}^p \Omega_i W_i$, where $p$ is the number of processes.  The RCholeskyQR process can be performed very efficiently using a single synchronization (required by the sketch), as outlined in~\Cref{fig:parallelrcholesky}. Once $\Omega W$ is computed, each process computes its QR factorization, and then the computed $R$ factor is used locally to compute a block of the orthogonal factor as $Q_i = W_i R^{-1}$. 

For each vector, the Gram-Scmidt process involves computing the projection coefficients onto the current basis, and then updating the vector by removing these projected components. In a parallel setting, the goal is to compute the projection coefficients and have them available on every process so that the second stage of the algorithm can be performed independently, without further communication. For randomized methods, this second stage has the same asymptotic cost as deterministic algorithms, namely, $nm^2$ floating-point operations. Their main advantage typically lies in how efficiently they build and replicate  the small matrix $\Omega Q$ from which these projection coefficients are computed.

\begin{figure}[t]
\centering
\begin{tikzpicture}[>=Latex, font=\sffamily, scale=1]

\foreach \i/\y in {1/3,2/2,3/1,4/0} {
  \node[anchor=east] at (-2.0,\y) {Process \i};
}

\foreach \y in {0.5,1.5,2.5} {
  \draw[dotted, thick] (-3.2,\y) -- (8.6,\y);
}

\foreach \i/\y in {1/3,2/2,3/1,4/0} {
  \draw[fill=colO, rounded corners=2pt]
      (-1.3,\y-0.15) rectangle (-0.4,\y+0.45);
  \node[font=\itshape\large]
      at (-0.85,\y+0.15) {$\Omega_{\i}$};
}

\foreach \i/\y in {1/3,2/2,3/1,4/0} {
  \draw[fill=colW, rounded corners=2pt]
      (-0.3,\y-0.45) rectangle (0.3,\y+0.45);
  \node[text=black, font=\itshape\large]
      at (0,\y) {$W_{\i}$};
}

\foreach \i/\y in {1/3,2/2,3/1,4/0} {
  \draw[->, thick] (0.8,\y) -- (1.7,\y);
  \node[anchor=west, font=\large]
      at (2,\y) {$\Omega_{\i} W_{\i}$};
}

\node[font=\large, text=darkteal] at (6.1,3.6) {Allreduce (+)};

\node[circle, fill=darkteal, minimum size=8pt] (c) at (5.5,1.5) {};

\foreach \y in {3,2,1,0} {
  \draw[thick, darkteal]
      (3.5,\y) .. controls (4.2,\y) and (4.8,1.5) .. (c);
  \draw[thick, darkteal]
      (c) .. controls (6.2,1.5) and (6.8,\y) .. (7.5,\y);
  \node[anchor=west, font=\large] at (7.6,\y) {$\Omega W$};
}

\end{tikzpicture}
\caption{Sketching a matrix $W$ partitioned into 4 blocks of rows over 4 processes. The sketching matrix $\Omega$ (Gaussian or s-hashing) is partitioned into 4 blocks of columns $\Omega_1, \Omega_2, \Omega_3, \Omega_4$.}
    \label{fig:parallelsketching}
\end{figure}
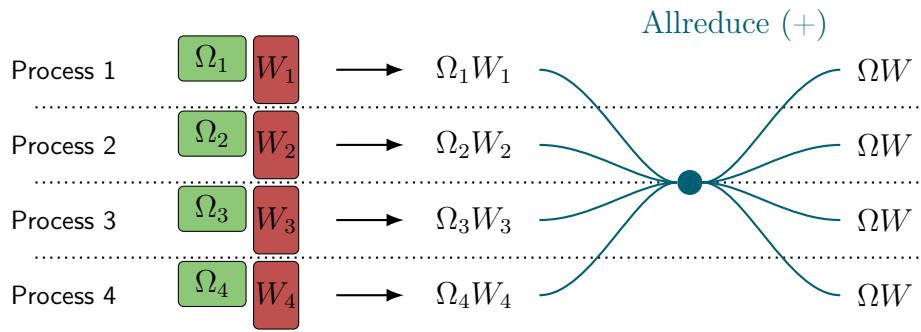

\begin{figure}[t]
\begin{center}
\begin{tikzpicture}
    \node at (-1,0.875) {Process 4};
    \node at (-1,2.625) {Process 3};
    \node at (-1,4.375) {Process 2};
    \node at (-1,6.125) {Process 1};
    
    \fill[lime] (0,0) rectangle ++(0.4,7);
    \draw[-] (0,0) -- (0.4,0) -- (0.4,7) -- (0,7) -- cycle;
    
    \node at (0.2,0.875) {$Q_4$};
    \node at (0.2,2.625) {$Q_3$};
    \node at (0.2,4.375) {$Q_2$};
    \node at (0.2,6.125) {$Q_1$};

    \fill[lime] (0.8,0.9) rectangle ++(0.33,0.6) node[below=0.6cm,color=black] {$\Omega Q$};
    \draw[-] (0.8,0.9) -- (1.13,0.9) -- (1.13,1.5) -- (0.8,1.5) -- cycle;
    
    \fill[lime] (0.8,2.65) rectangle ++(0.33,0.6) node[below=0.6cm,color=black] {$\Omega Q$};
    \draw[-] (0.8,2.65) -- (1.13,2.65) -- (1.13,3.25) -- (0.8,3.25) -- cycle;
    
    \fill[lime] (0.8,4.4) rectangle ++(0.33,0.6) node[below=0.6cm,color=black] {$\Omega Q$};
    \draw[-] (0.8,4.4) -- (1.13,4.4) -- (1.13,5.0) -- (0.8,5.0) -- cycle;
    
    \fill[lime] (0.8,6.15) rectangle ++(0.33,0.6) node[below=0.6cm,color=black] {$\Omega Q$};
    \draw[-] (0.8,6.15) -- (1.13,6.15) -- (1.13,6.75) -- (0.8,6.75) -- cycle;

    \fill[lime] (1.8,1.16) rectangle ++(0.33,0.33) node[below=0.4cm,color=black] {$R$};
    \draw[-] (1.8,1.16) -- (2.13,1.16) -- (2.13,1.49) -- (1.8,1.49) -- cycle;
    
    \fill[lime] (1.8,2.91) rectangle ++(0.33,0.33) node[below=0.4cm,color=black] {$R$};
    \draw[-] (1.8,2.91) -- (2.13,2.91) -- (2.13,3.24) -- (1.8,3.24) -- cycle;
    
    \fill[lime] (1.8,4.66) rectangle ++(0.33,0.33) node[below=0.4cm,color=black] {$R$};
    \draw[-] (1.8,4.66) -- (2.13,4.66) -- (2.13,5) -- (1.8,5) -- cycle;
    
    \fill[lime] (1.8,6.41) rectangle ++(0.33,0.33) node[below=0.4cm,color=black] {$R$};
    \draw[-] (1.8,6.41) -- (2.13,6.41) -- (2.13,6.74) -- (1.8,6.74) -- cycle;

    \fill[pink] (3,0) rectangle ++(0.1,7);
    \draw[-] (3,0) -- (3.1,0) -- (3.1,7) -- (3,7) -- cycle;
    \node at (3.1,4.375) {etc.};
    \node at (3.1,6.125) {$(\vec w_j)_{1:n/p}$};
    
    \fill[pink] (4.2,0.9) rectangle ++(0.1,0.6) node[below=0.6cm,color=black] {$\Omega \vec w_j$};
    \fill[pink] (4.2,2.65) rectangle ++(0.1,0.6) node[below=0.6cm,color=black] {$\Omega \vec w_j$};
    \fill[pink] (4.2,4.4) rectangle ++(0.1,0.6) node[below=0.6cm,color=black] {$\Omega \vec w_j$};
    \fill[pink] (4.2,6.15) rectangle ++(0.1,0.6) node[below=0.6cm,color=black] {$\Omega \vec w_j$};
    \draw[-] (4.2,0.9) -- (4.3,0.9) -- (4.3,1.5) -- (4.2,1.5) -- cycle;
    \draw[-] (4.2,2.65) --(4.3,2.65) --(4.3,3.25) --(4.2,3.25) -- cycle;
    \draw[-] (4.2,4.4) -- (4.3,4.4) -- (4.3,5.0) -- (4.2,5.0) -- cycle;
    \draw[-] (4.2,6.15) --(4.3,6.15) --(4.3,6.75) --(4.2,6.75) -- cycle;

    \node at (5.7,4) [color=blue] {\textbf{1 SYNC.}};
    \draw[->, blue, line width=1.5pt] (4.8,3.3) to[out=35,in=145] (6.5,3.3);

    \fill[lime] (7,0) rectangle ++(0.4,7);
    \fill[red] (7.4,0) rectangle ++(0.1,7);
    \draw[-] (7,0) -- (7.5,0) -- (7.5,7) -- (7,7) -- cycle;
    %
    \node at (7.2,0.875) {$Q_4$};
    \node at (7.2,2.625) {$Q_3$};
    \node at (7.2,4.375) {$Q_2$};
    \node at (7.2,6.125) {$Q_1$};

%
    \fill[lime] (7.8,0.9) rectangle ++(0.33,0.6) node[below=0.6cm,color=black] {$\Omega Q$};
    \fill[lime] (7.8,2.65) rectangle ++(0.33,0.6) node[below=0.6cm,color=black] {$\Omega Q$};
    \fill[lime] (7.8,4.4) rectangle ++(0.33,0.6) node[below=0.6cm,color=black] {$\Omega Q$};
    \fill[lime] (7.8,6.15) rectangle ++(0.33,0.6) node[below=0.6cm,color=black] {$\Omega Q$};
    \fill[red] (8.13,0.9) rectangle ++(0.1,0.6) ;
    \fill[red] (8.13,2.65) rectangle ++(0.1,0.6);
    \fill[red] (8.13,4.4) rectangle ++(0.1,0.6) ;
    \fill[red] (8.13,6.15) rectangle ++(0.1,0.6);
    \draw[-] (7.8,0.9) -- (8.23,0.9) -- (8.23,1.5) -- (7.8,1.5) -- cycle;
    \draw[-] (7.8,2.65) --(8.23,2.65) --(8.23,3.25) --(7.8,3.25) -- cycle;
    \draw[-] (7.8,4.4) -- (8.23,4.4) -- (8.23,5.0) -- (7.8,5.0) -- cycle;
    \draw[-] (7.8,6.15) --(8.23,6.15) --(8.23,6.75) --(7.8,6.75) -- cycle;
%
%
%
    \fill[lime] (8.8,1.16) rectangle ++(0.33,0.33) node[below=0.5cm,color=black] {$R$};
    \fill[lime] (8.8,2.91) rectangle ++(0.33,0.33) node[below=0.5cm,color=black] {$R$};
    \fill[lime] (8.8,4.66) rectangle ++(0.33,0.33) node[below=0.5cm,color=black] {$R$};
    \fill[lime] (8.8,6.41) rectangle ++(0.33,0.33) node[below=0.5cm,color=black] {$R$};
    \fill[red]  (8.8,1.06) rectangle ++(0.43,0.1);
    \fill[red]  (9.13,1.06) rectangle ++(0.1,0.43);
    \fill[red]  (8.8,2.81) rectangle ++(0.43,0.1);
    \fill[red]  (9.13,2.81) rectangle ++(0.1,0.43);
    \fill[red]  (8.8,4.56) rectangle ++(0.43,0.1);
    \fill[red]  (9.13,4.56) rectangle ++(0.1,0.43);
    \fill[red]  (8.8,6.31) rectangle ++(0.43,0.1);
    \fill[red]  (9.13,6.31) rectangle ++(0.1,0.43);
    
    \draw[-] (8.8,1.06) -- (9.23,1.06) -- (9.23,1.49) -- (8.8,1.49) -- cycle;
    \draw[-] (8.8,2.81) -- (9.23,2.81) -- (9.23,3.24) -- (8.8,3.24) -- cycle;
    \draw[-] (8.8,4.56) -- (9.23,4.56) -- (9.23,5   ) -- (8.8,5) -- cycle;
    \draw[-] (8.8,6.31) -- (9.23,6.31) -- (9.23,6.74) -- (8.8,6.74) -- cycle;


    \draw[dashed] (-1,1.75) -- (11.5,1.75); 
    \draw[dashed] (-1,3.5) -- (11.5,3.5); 
    \draw[dashed] (-1,5.25) -- (11.5,5.25); 
    
\end{tikzpicture}
\end{center}
\caption{Scattering of data in one iteration of RGS}
\label{fig:rgsparallel}
\end{figure}
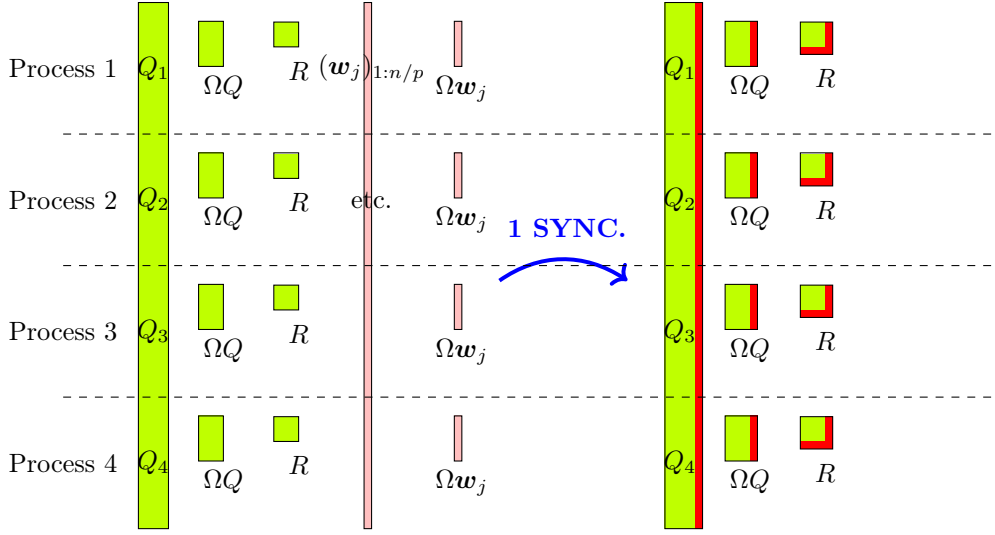

We remark that both RGS in~\Cref{algo:rgs} and RHQR in~\Cref{algo:rhqr_leftlooking} require between one and two synchronizations per iteration, depending on how the sketches can be grouped. This is $m$ or $2m$ synchronizations overall, which is the same computational cost as Classical Gram-Schmidt. As mentioned before, the computational cost of the two algorithms is dominated by the substitution process ($nm^2$ flops), and the total sketching cost, namely $2mt$ flops, where $t$ is the cost of sketching one vector. Thus, with SRHT, the total cost of both algorithms is $nm^2 + 2nm \log(m)$ flops asymptotically, which is half the flops of standard Gram-Schmidt processes. We illustrate the distribution of data in one iteration of RGS when performed on a parallel computer in~\Cref{fig:rgsparallel}. The matrix $Q \in \R^{n \times m}$ is distributed on each node by blocks of contiguous rows, and so is the next incoming input basis vector $\vec w_j \in \R^n$. The sketches $\Omega Q$ and $\Omega \vec w$ are available on each processor, and so is $R \in \R^{m \times m}$. Each processor can compute locally $(\Omega Q_{1:j-1})^T \Omega \vec  w_j$ and thus update the first $j-1$ entries of $R \vec e_j$. It also allows them to perform their part of the vector refresh $\vec w_j \gets \vec w_j - Q_{1:j-1} (\Omega Q_{1:j-1})^T \Omega \vec w_j$. The result is then sketched, which allows every processor to compute the $j$-th entry of $R \vec e_j$, scale their part of the refreshed $\vec w_j$, and thus obtain the new basis vector $\vec q_j \in \R^n$. All the matrices are updated with the vectors computed in this iteration.

\begin{figure}[t]
\begin{center}
\begin{tikzpicture}
    \node at (-1,0.875) {Process 4};
    \node at (-1,2.625) {Process 3};
    \node at (-1,4.375) {Process 2};
    \node at (-1,6.125) {Process 1};
    
    \fill[lime] (0,0) rectangle ++(0.4,7);
    \draw[-] (0,0) -- (0.4,0) -- (0.4,7) -- (0,7) -- cycle;
    
    \node at (0.2,0.875) {$U_4$};
    \node at (0.2,2.625) {$U_3$};
    \node at (0.2,4.375) {$U_2$};
    \node at (0.2,6.125) {$U_1$};

    \fill[lime] (0.8,0.9) rectangle ++(0.33,0.6) node[below=0.6cm,color=black] {$\Psi U$};
    \draw[-] (0.8,0.9) -- (1.13,0.9) -- (1.13,1.5) -- (0.8,1.5) -- cycle;
    
    \fill[lime] (0.8,2.65) rectangle ++(0.33,0.6) node[below=0.6cm,color=black] {$\Psi U$};
    \draw[-] (0.8,2.65) -- (1.13,2.65) -- (1.13,3.25) -- (0.8,3.25) -- cycle;
    
    \fill[lime] (0.8,4.4) rectangle ++(0.33,0.6) node[below=0.6cm,color=black] {$\Psi U$};
    \draw[-] (0.8,4.4) -- (1.13,4.4) -- (1.13,5.0) -- (0.8,5.0) -- cycle;
    
    \fill[lime] (0.8,6.15) rectangle ++(0.33,0.6) node[below=0.6cm,color=black] {$\Psi U$};
    \draw[-] (0.8,6.15) -- (1.13,6.15) -- (1.13,6.75) -- (0.8,6.75) -- cycle;

    \fill[lime] (1.8,1.16) rectangle ++(0.33,0.33) node[below=0.4cm,color=black] {$T$};
    \draw[-] (1.8,1.16) -- (2.13,1.16) -- (2.13,1.49) -- (1.8,1.49) -- cycle;
    
    \fill[lime] (1.8,2.91) rectangle ++(0.33,0.33) node[below=0.4cm,color=black] {$T$};
    \draw[-] (1.8,2.91) -- (2.13,2.91) -- (2.13,3.24) -- (1.8,3.24) -- cycle;
    
    \fill[lime] (1.8,4.66) rectangle ++(0.33,0.33) node[below=0.4cm,color=black] {$T$};
    \draw[-] (1.8,4.66) -- (2.13,4.66) -- (2.13,5) -- (1.8,5) -- cycle;
    
    \fill[lime] (1.8,6.41) rectangle ++(0.33,0.33) node[below=0.4cm,color=black] {$T$};
    \draw[-] (1.8,6.41) -- (2.13,6.41) -- (2.13,6.74) -- (1.8,6.74) -- cycle;

    \fill[pink] (3,0) rectangle ++(0.1,7);
    \draw[-] (3,0) -- (3.1,0) -- (3.1,7) -- (3,7) -- cycle;
    \node at (3.1,4.375) {etc.};
    \node at (3.1,6.125) {$(\vec w_j)_{1:n/p}$};
    
    \fill[pink] (4.2,0.9) rectangle ++(0.1,0.6) node[below=0.6cm,color=black] {$\Psi \vec w_j$};
    \fill[pink] (4.2,2.65) rectangle ++(0.1,0.6) node[below=0.6cm,color=black] {$\Psi \vec w_j$};
    \fill[pink] (4.2,4.4) rectangle ++(0.1,0.6) node[below=0.6cm,color=black] {$\Psi \vec w_j$};
    \fill[pink] (4.2,6.15) rectangle ++(0.1,0.6) node[below=0.6cm,color=black] {$\Psi \vec w_j$};
    \draw[-] (4.2,0.9) -- (4.3,0.9) -- (4.3,1.5) -- (4.2,1.5) -- cycle;
    \draw[-] (4.2,2.65) --(4.3,2.65) --(4.3,3.25) --(4.2,3.25) -- cycle;
    \draw[-] (4.2,4.4) -- (4.3,4.4) -- (4.3,5.0) -- (4.2,5.0) -- cycle;
    \draw[-] (4.2,6.15) --(4.3,6.15) --(4.3,6.75) --(4.2,6.75) -- cycle;

    \node at (5.7,4) [color=blue] {\textbf{1 SYNC.}};
    \draw[->, blue, line width=1.5pt] (4.8,3.3) to[out=35,in=145] (6.5,3.3);

    \fill[lime] (7,0) rectangle ++(0.4,7);
    \fill[red] (7.4,0) rectangle ++(0.1,7);
    \draw[-] (7,0) -- (7.5,0) -- (7.5,7) -- (7,7) -- cycle;
    %
    \node at (7.2,0.875) {$U_4$};
    \node at (7.2,2.625) {$U_3$};
    \node at (7.2,4.375) {$U_2$};
    \node at (7.2,6.125) {$U_1$};

%
    \fill[lime] (7.8,0.9) rectangle ++(0.33,0.6) node[below=0.6cm,color=black] {$\Psi U$};
    \fill[lime] (7.8,2.65) rectangle ++(0.33,0.6) node[below=0.6cm,color=black] {$\Psi U$};
    \fill[lime] (7.8,4.4) rectangle ++(0.33,0.6) node[below=0.6cm,color=black] {$\Psi U$};
    \fill[lime] (7.8,6.15) rectangle ++(0.33,0.6) node[below=0.6cm,color=black] {$\Psi U$};
    \fill[red] (8.13,0.9) rectangle ++(0.1,0.6) ;
    \fill[red] (8.13,2.65) rectangle ++(0.1,0.6);
    \fill[red] (8.13,4.4) rectangle ++(0.1,0.6) ;
    \fill[red] (8.13,6.15) rectangle ++(0.1,0.6);
    \draw[-] (7.8,0.9) -- (8.23,0.9) -- (8.23,1.5) -- (7.8,1.5) -- cycle;
    \draw[-] (7.8,2.65) --(8.23,2.65) --(8.23,3.25) --(7.8,3.25) -- cycle;
    \draw[-] (7.8,4.4) -- (8.23,4.4) -- (8.23,5.0) -- (7.8,5.0) -- cycle;
    \draw[-] (7.8,6.15) --(8.23,6.15) --(8.23,6.75) --(7.8,6.75) -- cycle;
%
%
%
    \fill[lime] (8.8,1.16) rectangle ++(0.33,0.33) node[below=0.5cm,color=black] {$T$};
    \fill[lime] (8.8,2.91) rectangle ++(0.33,0.33) node[below=0.5cm,color=black] {$T$};
    \fill[lime] (8.8,4.66) rectangle ++(0.33,0.33) node[below=0.5cm,color=black] {$T$};
    \fill[lime] (8.8,6.41) rectangle ++(0.33,0.33) node[below=0.5cm,color=black] {$T$};
    \fill[red]  (8.8,1.06) rectangle ++(0.43,0.1);
    \fill[red]  (9.13,1.06) rectangle ++(0.1,0.43);
    \fill[red]  (8.8,2.81) rectangle ++(0.43,0.1);
    \fill[red]  (9.13,2.81) rectangle ++(0.1,0.43);
    \fill[red]  (8.8,4.56) rectangle ++(0.43,0.1);
    \fill[red]  (9.13,4.56) rectangle ++(0.1,0.43);
    \fill[red]  (8.8,6.31) rectangle ++(0.43,0.1);
    \fill[red]  (9.13,6.31) rectangle ++(0.1,0.43);
    
    \draw[-] (8.8,1.06) -- (9.23,1.06) -- (9.23,1.49) -- (8.8,1.49) -- cycle;
    \draw[-] (8.8,2.81) -- (9.23,2.81) -- (9.23,3.24) -- (8.8,3.24) -- cycle;
    \draw[-] (8.8,4.56) -- (9.23,4.56) -- (9.23,5   ) -- (8.8,5) -- cycle;
    \draw[-] (8.8,6.31) -- (9.23,6.31) -- (9.23,6.74) -- (8.8,6.74) -- cycle;

    \fill[blue] (10,6.6) rectangle ++(0.1,0.4);
    \draw[-] (10,0) -- (10.1,0) -- (10.1,7) -- (10,7) -- cycle;
    \node at (10.3,6.5) {$r_j$};

    \draw[dashed] (-1,1.75) -- (11.5,1.75); 
    \draw[dashed] (-1,3.5) -- (11.5,3.5); 
    \draw[dashed] (-1,5.25) -- (11.5,5.25); 
    
\end{tikzpicture}
\end{center}
\caption{Scattering of data in one iteration of RHQR}
\label{fig:rhqrparallel}
\end{figure}
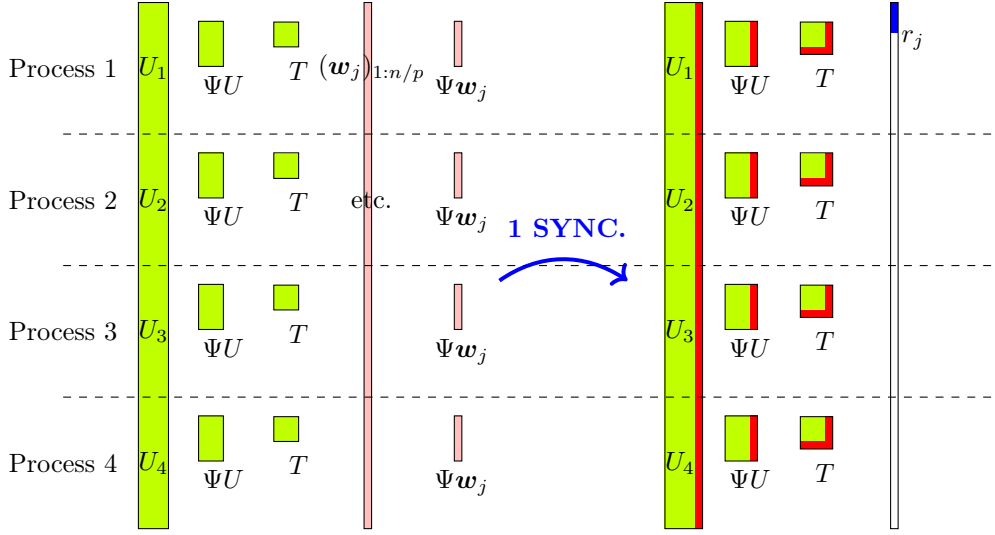

We illustrate the distribution of data in one iteration of RHQR when performed on a parallel computer in~\Cref{fig:rhqrparallel}. The tall-and-skinny matrix $U \in \R^{n \times m}$ is distributed by blocks of contiguous rows, as the input (it can be written in its place), and so is the next basis vector $\vec w \in \R^n$. The sketches $\Psi U \in \R^{(\ell+m) \times m}$ and $\Psi \vec w \in \R^{\ell + m}$ is available on each processors, and so is $T \in \R^{m \times m}$. All processors can thus compute locally $T^T (\Psi U)^T \Psi \vec w$, which allows each processors to perform its part of the vector refresh $\vec w \gets \vec w - U T^T (\Psi U)^T \Psi \vec w$. The resulting, refreshed vector is then sketched, which allows every processors to compute its share of the associated randomized Householder vector. All the matrices are then updated with the vectors computed in this iteration.

It is very natural to use a mixed-precision setting in the case of randomized orthogonalization. A common approach is to store the high-dimensional matrix in a \textit{coarse} floating point format (typically 32 bits, or even 16 bits when the CPU supports it), while casting and storing low-dimensional matrices (mainly the sketches and the triangular factors) in a \textit{fine} floating point format (typically 64, or even 128 bits when the CPU supports it). In the case of RCholeskyQR in~\Cref{fig:parallelrcholesky} for instance, the matrix $W$ is stored in coarse precision. Then $\Omega$ (often not stored explicitly) is applied, using either coarse or fine arithmetic operations. The result $\Omega W$ is stored in fine precision, and the $R$ factor is computed from $\Omega W$ in fine precision. The computed $R$ factor is finally cast to coarse precision, and the substitution for $Q$ is performed in coarse precision.

A critical aspect of finite-precision analysis of randomized algorithms is the \textit{forward accuracy} of the sketching step, namely the bounding of the magnitude of $\| \Omega \vec x - \fl[\Omega \vec x] \| / \|\Omega \vec x\|$, where $\fl[\Omega \vec x]$ denotes the finite-precision result of the routine applying $\Omega$ to $\vec x$. If this forward error is sufficiently small, the stability of the process can potentially reduce to the accuracy of the  orthogonalization of the sketch~\cite{BalabanovGrigori2022,rhqr,Higgins2025multisketching}. Authors in~\cite{BalabanovGrigori2022} describe the accuracy of the sketching process, performed in fine precision $\e$, for any sketching matrix $\Omega \in \R^{\ell \times n}$, basing their analysis on the following backward error~\cite{Higham08} 
\[ \fl [\Omega \vec x] = (\Omega + \Delta \Omega) \cdot \vec x, \quad |\Delta \Omega| \leq \mathcal{O}(n^{1/2}) |\Omega | \e \quad \text{(w.h.p)}\]
By restricting $\Omega$ to some specific OSE distributions, the accuracy of the sketching step can be even more favorable. Authors in~\cite{rhqr} base their analysis on the observation that the application of the SRHT, in fine precision $\e$, using the Fast Walsh-Hadamard Transform, is \textit{backward stable}, with fine constants even in the worst-case:
\[ \forall \vec x \in \R^n, \quad \fl [ \Omega \vec x ] = \Omega (\vec x + \Delta \vec x), \quad \|\Delta \vec x\| \leq \mathcal{O}(\log(n)) \| \vec x \| \e \quad \text{(SRHT)}. \]

The authors in~\cite{BalabanovGrigori2022} show that, if $\mathcal{O}(m^2) \mathrm{Cond}(W) \f < 1$, if $\mathcal{O}(m^{1/2} n^{3/2} + m^{3/2} \ell) \e < \f$, the output of RGS in~\Cref{algo:rgs} verifies
\begin{equation}
    \mathrm{Cond}(\fl[Q]) \leq (1+\mathcal{O}(\epsilon)) \left(1 + \mathcal{O}(m^2) \mathrm{Cond}(W) \e \right).
\end{equation} 
Finally, using SRHT sketching and assuming that $\mathcal{O}(\log(n) + \ell + m) \e \leq \f$, the authors in~\cite{rhqr} show that the output of the RHQR process in~\Cref{algo:rhqr_leftlooking} verify
\[ \mathrm{Cond}(\fl [Q]) \leq (1 + \mathcal{O}(\epsilon)) (1 + \mathcal{O}\ell m^{3/2} \f), \quad \| (\fl[QR] - W) \vec e_j \| \leq (1 + \mathcal{O}(\epsilon)) \mathcal{O}(\ell m^{3/2}) \| W \vec e_j\| \f.\]
It is often observed that the factorization error of randomized algorithm is smaller in practice than that of standard orthogonalization. Indeed, with fast and stable sketching processes, the coefficients of the matrix $R$ driving the substitution are potentially obtained with much less flops and better numerical stability than those of the determinsitic $R$. This phenomenon has very concrete consequences, for instance, in the case of randomized GMRES, as we outline in the next sections.

%% file: krylov.tex
\section{Krylov subspace methods}
\label{sec:krylov} 

Krylov subspace methods are a valuable tool for the solution of various problems in numerical linear algebra. We review here the the fundamental definitions and properties associated with Krylov subspaces, and we introduce the corresponding randomized versions.

Given a matrix $A \in \R^{n \times n}$ and a vector $\vec b \in \R^n$, the associated Krylov subspace of dimension $m$ is defined as $\K_m(A, \vec b)= \vspan\{\vec b, A \vec b, \dots, A^{m-1} \vec b\}$. Defining $\vec q_1 = \vec b / \norm{\vec b}_2$, the Arnoldi process \cite[Algorithm~6.2]{Saad03} can be used to generate an orthonormal basis $Q_m = [\vec q_1 \,\cdots\, \vec q_m]$ for $\K_m(A, \vec b)$. This basis satisfies the Arnoldi relation
\begin{equation}
	\label{eqn:arnoldi-relation}
	A Q_m = Q_{m+1} \underline{H}_m = Q_m H_m + h_{m+1,m} \vec q_{m+1} \vec e_m^T, \qquad \text{with} \quad \underline{H}_m = \begin{bmatrix}
		H_m \\ h_{m+1, m} \vec e_m^T
	\end{bmatrix},
\end{equation}   
where $\vec q_{m+1} \perp Q_m$ and $H_m \in \R^{m \times m}$ is an upper Hessenberg matrix containing the orthogonalization coefficients. Since $Q_{m+1}$ has orthonormal columns, it follows from \cref{eqn:arnoldi-relation} that $Q_m^T A Q_m = H_m$.

\subsection{Randomized Arnoldi and Krylov factorizations}
\label{subsec:rand-Arnoldi}

Let $W_m$ be any basis of $\K_m(A, \vec b)$ where each new basis vector $\vec w_{j+1}$ is generated iteratively as a linear combination of $A \vec w_{j}$ and the columns of the current basis $W_j$. Such a basis satisfies the Arnoldi-like relation
\begin{equation}
	\label{eqn:arnoldi-like-relation}
	A W_m = W_{m+1} \underline{L}_m = W_m L_m + \vec w_{m+1}  \vec \ell_m^T, \qquad \vec \ell_m := \ell_{m+1,m} \vec e_m,
\end{equation}
where $\underline{L}_m \in \R^{(m+1) \times m}$ is upper Hessenberg but $W_m$ does not necessarily have orthonormal columns. In this case, we have the identity
\begin{equation}
	\label{eqn:arnoldi-projection-identity--general-basis}
	W_m^+ A W_m = L_m + W_m^+ \vec w_{m+1} \vec \ell_m^T.
\end{equation}
The randomized Arnoldi process \cite{BalabanovGrigori2022} given in \cref{algo:rando-arnoldi} constructs a basis of $\K_m(A, \vec b)$ employing the randomized Gram--Schmidt algorithm to generate a sketch-orthonormal basis $V_m$ and an associated Hessenberg matrix $\underline{G}_m$, which also satisfy the randomized Arnoldi relation
\begin{equation}
	\label{eqn:arnoldi-like-relation--sketched-basis}
	A V_m = V_{m+1} \underline{G}_m = V_m G_m + \vec v_{m+1}  \vec g_m^T, \qquad \vec g_m := g_{m+1,m} \vec e_m.
\end{equation}
In this case, by multiplying \cref{eqn:arnoldi-like-relation--sketched-basis} from the left by $(\Omega V_m)^T\Omega$ and using the fact that $V_m$ satisfies $(\Omega V_{m+1})^T \Omega V_{m+1} = I$, we obtain the following
\begin{equation*}
	(\Omega V_m)^T \Omega A V_m = G_m.
\end{equation*} 
This identity allows us to use the Hessenberg matrix $G_m$ generated by the randomized Arnoldi process to efficiently construct approximate solutions for different problems, ranging from the solution of linear systems and matrix equations to the evaluation of matrix functions and the computation of eigenvalues.  
\begin{algorithm}[t]
\caption{Randomized Arnoldi process \cite{BalabanovGrigori2022}}
\label{algo:rando-arnoldi}
\begin{algorithmic}[1]
\Require $A \in \R^{n \times n}$, a starting vector $\vec b \in \R^{n}$, a Krylov dimension $m$, a matrix $\Omega \in \R^{\ell \times n}$ such that $\Ker{\Omega} \cap \K_m(A, \vec b) = \emptyset$.
\Ensure A randomized Arnoldi factorization $A V_m = V_{m+1} \underline{G}_m = V_m G_m + \vec v_{m+1}  \vec g_m^T$ as in \cref{eqn:arnoldi-like-relation--sketched-basis} with $S_{m+1} = \Omega V_{m+1}$ explicitly computed. 
\Function{Randomized Arnoldi}{$A, \vec b, m, \Omega$}
    \State $\vec z_1 \gets \Omega \vec b$
    \State $V_1 \gets [\vec b / \| \vec z\|]$, $S_1 \gets [\vec z / \|\vec z\|]$
    \For{$j = 1:m$}
        \State $\vec v_{j+1} \gets  A \vec v_{j}$ \label{algo:rarno:applyA}
        \State $\vec z \gets \Omega \vec v_{j+1}$ \label{algo:rarno:sketch1}
        \State $\vec g_{j} \gets S_{j}^+ \vec z \in \R^{j}$ \label{algo:rarno:lls} \# use a method stable enough to handle $S_{j}^+$
        \State $\vec v_{j+1} \gets \vec v_{j+1} - V_{j} \vec g_{j}$ \label{algo:rarno:refresh}
        \State $\vec z \gets \Omega \vec v_{j+1}$ \label{algo:rarno:sketch2}
		\State $g_{j+1,j} \gets  \| \vec z\|$
        \State $\underline{G}_{j} \gets \left[ \begin{array}{ccc|c} & & & \\ & \underline{G}_{j-1} &  & \vec g_{j} \\ & & & \\\hline  & 0_{1\times (j-1)} & &  g_{j+1,j} \end{array} \right]$ 
		\State $V_{j+1} \gets [V_{j} \; | \; \vec v_{j+1} / g_{j+1,j} ] $, $\quad S_{j+1} \gets [S_{j} \; | \; \vec z / g_{j+1,j}]$ \label{algo:rarno:scalew}
    \EndFor
    \State \Return $V_{m+1}, S_{m+1}, \underline{G}_m$
\EndFunction
\end{algorithmic}
\end{algorithm}

A major distinction of the randomized Arnoldi algorithm compared to deterministic Arnoldi is that the method does not reduce to a short recurrence when the matrix $A$ is symmetric. In the deterministic case $H_m = Q_m^T A Q_m$ from \cref{eqn:arnoldi-relation} is symmetric and upper Hessenberg, hence it is necessarily tridiagonal. Since its entries are inner products between successive Krylov basis vectors, it yields the three‑term Lanczos recurrence and allows orthogonalization only against the two most recent basis vectors at each iteration; see \cite{Lanczos1950iterationmethodsolution}. By contrast, in the randomized Arnoldi process the Hessenberg matrix $G_m = (\Omega V_m)^T \Omega A V_m$ from \cref{eqn:arnoldi-like-relation--sketched-basis} arising from the sketched oblique projection need not be symmetric. Consequently, no short recurrence is obtained in general, and one must exercise care when using $G_m$ to approximate eigenvalues of $A$ as later developed in \cref{subsec:sketch-RR}, since spurious or complex eigenvalues may appear.

\subsection{Whitening}
\label{subsubsec:whitening}

In several recent works \cite{NakatsukasaTropp24,GuettelSchweitzer23,PSS25mf,PSS25me}, it has been proposed to construct the sketch-orthonormal basis $V_m$ implicitly. First, a non-orthogonal basis $W_m$ is generated with a cheap procedure, such as a $k$-truncated Arnoldi process (see, e.g., \cite[Chapter 6.4.2]{Saad03}), in which the orthogonalization of a new vector is only performed against the last $k$ basis vectors, and then a QR factorization of the sketched basis $\Omega W_m = S_m R_m$ is computed. A sketch-orthonormal basis of $\K_m(A, \vec b)$ can then be obtained as $V_m = W_m R_m^{-1}$. This process is often called \emph{whitening} in the literature and coincides with an implicit application of the randomized Cholesky QR framework described at the end of \cref{subsec:sketch-orth-discussion}. 

The main advantage of this approach is that it allows for a cheaper, implicit computation of the sketch-orthonormal basis~$V_m$. 
Indeed, note that explicitly forming $V_m = W_m R_m^{-1}$ requires $\mathcal{O}(nm^2)$ operations, which coincides with the computational cost of a direct sketch-orthogonalization via randomized Gram--Schmidt. However, the main advantage of the whitening strategy is that, in many cases, there is no need to form the basis $V_m$ explicitly. For example, if the desired solution to a certain problem has the form $V_m \vec y_m$ for some $\vec y_m \in \R^m$, it can be computed as $W_m (R_m^{-1} \vec y_m)$. This operation only costs $\mathcal{O}(mn + m^2)$, since $R_m^{-1}$ is applied to a vector of length $m$, thus reducing the asymptotic computational complexity of the method, provided that the basis $W_m$ is computed using a cheap procedure.
When $W_m$ is constructed using the $k$-truncated Arnoldi process, the orthogonalization cost in the computation of $W_m$ drops to $\mathcal{O}(kmn)$, and implicit whitening of the basis requires $\mathcal{O}(m^3)$ operations, making this approach asymptotically cheaper than directly performing the randomized Gram--Schmidt process. In exact arithmetic, the (implicitly computed) whitened basis $V_m$ coincides with the basis obtained by randomized Gram--Schmidt, ensuring that the two approaches produce the same approximations for any task that employs a sketch-orthonormal basis of $\K_m(A, \vec b)$. An ulterior advantage of the $k$-truncated Arnoldi process is that only the $k+1$ vectors need to be kept in memory during the generation of the basis $W_m$, making the whitening approach viable in a low-memory setting. If the full basis $W_m$ is required to form the approximate solution, it can be generated again with a two-pass approach, without the need to store it in full: see, for instance, \cite[Section~4.2]{GuettelSchweitzer23} and \cite{Borici00,FrommerSimoncini08}.

The main limitation of whitening is that the non-orthogonal basis $W_m$ typically becomes severely ill-conditioned even for moderate $m$, and as a consequence both the QR factorization $\Omega W_m = S_m R_m$ and the application of $R_m^{-1}$ may suffer from numerical instability, eventually yielding approximations that in finite precision may diverge significantly from those obtained with a sketch-orthonormal basis $V_m$ constructed explicitly through a randomized Gram--Schmidt process. Although the numerical behavior of whitening within Krylov subspace methods is still not completely understood theoretically, encouraging numerical results have been observed in several applications, often even better than what one would expect from the growth of $\kappa(W_m)$ \cite{CKN24,PSS25mf}. Various selective orthogonalization strategies have been explored in \cite{GuettelSimunec24} as alternatives to $k$-truncated Arnoldi, aiming to mitigate the ill-conditioning of the non-orthogonal basis $W_m$ before applying whitening. 

When using whitening to sketch-orthogonalize the basis $W_m$, the following relations may be useful. Similarly to \cref{eqn:arnoldi-projection-identity--general-basis}, we have
\begin{equation*}
	(\Omega W_m)^+ \Omega A W_m = L_m + \vec z_m \vec \ell_m^T, \qquad \vec z_m := (\Omega W_{m})^+ \Omega \vec w_{m+1}
\end{equation*}
which gives us the following alternative expression for the Arnoldi relation \cref{eqn:arnoldi-like-relation},
\begin{equation}
	\label{eqn:arnoldi-like-relation--sketchorth-lastvec}
	A W_m = W_m (L_m + \vec z_m \vec \ell_m^T) + \widehat{\vec v}_{m+1}  \vec \ell_m^T, \qquad \widehat{\vec v}_{m+1} := \vec w_{m+1} - W_m \vec z_m,
\end{equation}
where $\vec{\widehat v}_{m+1} \perp^\Omega W_m$. In other words, we can add a rank-one perturbation to $L_m$ in order to obtain an Arnoldi relation in which the last basis vector $\widehat{\vec v}_{m+1}$ is sketch-orthogonal to the basis $W_m$. 

The vector $\vec z_m$ can be alternatively written in terms of a QR factorization of the sketched basis $\Omega W_{m+1}$. Indeed, assuming that we have the QR factorization
\begin{equation*}
	\Omega W_{m+1} = \begin{bmatrix}
		S_m & \vec s_{m+1}
	\end{bmatrix}
	\begin{bmatrix}
		R_m & \vec r_{m} \\
		& \rho_{m+1}
	\end{bmatrix},
\end{equation*} 
we have
\begin{equation}
	\label{eqn:zm-expression-from-sketch-qr}
	\vec z_m = (\Omega W_m)^+ \Omega \vec w_{m+1} = R_m^{-1} S_m^T (S_m \vec r_{m} + \rho_{m+1} \vec s_{m+1}) = R_m^{-1} \vec r_m.
\end{equation}
A variant of the identity \cref{eqn:arnoldi-like-relation--sketchorth-lastvec} with the expression for $\vec z_m$ in \cref{eqn:zm-expression-from-sketch-qr} is given in \cite[Proposition~3.1]{PSS25mf}, where it is referred to with the name \emph{sketched Arnoldi relation}. 

Multiplying \cref{eqn:arnoldi-like-relation--sketchorth-lastvec} by $R_m^{-1}$ on the right, we obtain
\begin{equation*}
	A W_m R_m^{-1} = W_m R_m^{-1} (R_m L_m R_m^{-1} + \vec z_m \vec \ell_m^T R_m^{-1}) + \widehat{\vec v}_{m+1}  \vec \ell_m^T R_m^{-1},
\end{equation*}
and using $V_m = W_m R_m^{-1}$ and $\vec \ell_m^T R_m^{-1} = \rho_m^{-1} \vec \ell_m^T$, where $\rho_m = [R_m]_{m,m}$, we get   
\begin{equation}
	\label{eqn:whitened-sketched-arnoldi-relation}
	A V_m = V_m (\widehat{L}_m + \rho_m^{-1} \vec z_m \vec \ell_m^T) + \rho_m^{-1} \widehat{\vec v}_{m+1}  \vec \ell_m^T, \qquad \widehat{L}_m := R_m^{-1} L_m R_m.
\end{equation}
Multiplying \cref{eqn:whitened-sketched-arnoldi-relation} by $(\Omega V_m)^T \Omega$ from the left and using $(\Omega V_m)^T \Omega V_m = I$ and $(\Omega V_m)^T \Omega \widehat{\vec v}_{m+1}= 0$, we obtain
\begin{equation}
	\label{eqn:whitened-hessenberg-matrix}
	G_m = (\Omega V_m)^T \Omega A V_m = \widehat{L}_m + \rho_m^{-1} \vec z_m \vec \ell_m^T.
\end{equation}
The identity \cref{eqn:whitened-sketched-arnoldi-relation} is called \emph{whitened-sketched Arnoldi relation} in \cite{PSS25mf}, and \cref{eqn:whitened-hessenberg-matrix} provides an explicit expression for the Hessenberg matrix associated with the whitened basis $V_m$, which can be evaluated cheaply by only using the $\underline{L}_m$ and the upper triangular factor from the QR factorzization of $\Omega W_{m+1}$.

%% file: linearsystems.tex
\section{Solution of linear systems}
\label{sec:linear-systems}
Let us consider the linear system $A \vec x = \vec b$, with $A \in \R^{n \times n}$ and $\vec b \in \R^n$. Subspace projection methods \cite{Saad03} are a very effective tool for solving large-scale linear systems. These methods seek an approximate solution $\vec x_m$ which satisfies the two following conditions:
\begin{itemize}
	\item $\vec x_m$ is contained in the Krylov subspace $\K_m(A, \vec b)$,
	\item the residual $\vec b - A\vec x_m$ satisfies the \emph{Petrov-Galerkin condition} $\vec b - A \vec x_m \perp \mathcal{L}_m$, where $\mathcal{L}_m \subset \R^n$ is an $m$-dimensional subspace.
\end{itemize} 
When $\mathcal{L}_m = \K_m(A, \vec b)$, the resulting methods is usually known as an orthogonal projection method, while with a more general $\mathcal{L}_m$ we obtain an oblique projection method.

\subsection{Krylov methods for linear systems}
We start by briefly reviewing the well-known GMRES and FOM methods for the solution of $A \vec x = \vec b$. 

GMRES \cite{SaadSchultz86} takes $\mathcal{L}_m = A \K_m(A, \vec b)$, so it seeks an approximate solution $\xgmres_m \in \K_m(A, \vec b)$ whose residual satisfies 
\begin{equation*}
	\vec b - A \xgmres_m \perp A \K_m(A, \vec b).
\end{equation*} 
Recalling the Arnoldi relation \cref{eqn:arnoldi-relation} and writing $\xgmres_m = Q_m \ygmres_m$ with $\ygmres_m \in \R^m$, we can equivalently rewrite this condition as
\begin{align*}
	0 = (A Q_m)^T(\vec b - A Q_m \ygmres_m) &= (Q_{m+1} \underline{H}_m)^T (\vec b - Q_{m+1}\underline{H}_m \ygmres_m) \\
	&= \underline{H}_{m}^T \widetilde{\beta} \vec e_1 - \underline{H}_m^T \underline{H}_m \ygmres_m,  
\end{align*}
where $\widetilde{\beta} \in \R$ is defined from the identity $\vec b = \widetilde{\beta} Q_m \vec e_1$.
Observe that these are the normal equations associated with a least squares problem, so we can write the GMRES approximate solution as 
\begin{equation}
	\label{eqn:linear-system-gmres-sol}
	\xgmres_m = Q_m \ygmres_m \qquad \text{where} \qquad \ygmres_m = \argmin_{\vec y \in \R^m} \norm{\underline{H}_m \ygmres_m - \widetilde{\beta} \vec e_1}.
\end{equation}
From the condition $(A Q_m)^T (\vec b - A Q_m \ygmres_m) = 0$, we also obtain that  
\begin{equation*}
	\xgmres_m = \argmin_{\vec x \in \K_m(A, \vec b)} \norm{A \vec x - \vec b},
\end{equation*}
i.e., the GMRES solution $\xgmres_m$ minimizes the residual over the Krylov subspace $\K_m(A,\vec b)$.  

On the other hand, the Full Orthogonalization Method (FOM) \cite{Arnoldi1951principleminimizediterations, Saad03} employs $\mathcal{L}_m = \K_m(A, \vec b)$, thus requiring that the residual is orthogonal to the Krylov subspace $\K_m(A, \vec b)$, a condition that is also known as \emph{Galerkin condition}. Let us denote by $\xfom_m$ the approximate solution after $m$ iterations of FOM. Recalling the Arnoldi relation \cref{eqn:arnoldi-relation} and writing $\xfom_m = Q_m \yfom_m$ with $\yfom_m \in \R^m$, we can rewrite the Galerkin condition as
\begin{align*}
	0 = Q_m^T (\vec b - A Q_m \yfom_m) &= Q_m^T (\vec b - (Q_m H_m + h_{m+1,m} \vec q_{m+1} \vec e_m^T) \yfom_m) \\
	&= \widetilde{\beta} \vec e_1 - H_m \yfom_m,
\end{align*}
where we used the orthogonality of the columns of $Q_{m+1}$ in the last equality. This yields the expression for the FOM approximate solution
\begin{equation}
	\label{eqn:linear-system-fom-sol}
	\xfom_m = Q_m \yfom_m, \qquad \text{where} \qquad H_m \yfom_m = \widetilde{\beta} \vec e_1.
\end{equation}
When $A$ is symmetric positive definite, FOM simplifies to the conjugate gradient (CG)~\cite{HestenesStiefel52}. In this case, the approximate solution $\xfom_m$ can be iteratively updated using short recurrences and it additionally satisfies the optimality property
\begin{equation*}
	\xfom_m = \argmin_{\widehat{\vec x} \in \K_m(A, \vec b)} \norm{\vec x - \widehat{\vec x}}_A, \qquad \text{where} \qquad \norm{\vec z}_A = (\vec z^T A \vec z)^{1/2}.
\end{equation*}

\subsection{Randomized Krylov methods for linear systems}

In this section, we present algorithms that employ a sketch-orthonormal basis of the Krylov subspace $\K_m(A, \vec b)$ to solve the linear system $A \vec x = \vec b$, following the presentation in \cite{NakatsukasaTropp24} and \cite{TGB23}. These algorithms replace the Petrov-Galerkin imposed by standard Krylov methods with a similar condition on the sketched residual. Specifically, they seek an approximate solution $\vec x_m$ that satisfies:
\begin{itemize}
	\item $\vec x_m$ belongs to the Krylov subspace $\K_m(A, \vec b)$,
	\item the residual $\vec b - A\vec x_m$ satisfies the \emph{sketched Petrov-Galerkin condition} $\vec b - A \vec x_m \perp^\Omega \mathcal{L}_m$, where $\mathcal{L}_m \subset \R^n$ is an $m$-dimensional subspace and $\Omega$ is an $\epsilon$-embedding for $\K_m(A, \vec b)$.       
\end{itemize} 
In analogy with the deterministic Krylov methods, choosing $\mathcal{L}_m = \K_m(A, \vec b)$ or $\mathcal{L}_m = A \K_m(A, \vec b)$ leads, respectively, to the \emph{randomized FOM} and \emph{randomized GMRES} algorithms (also known as sketched FOM and sketched GMRES in the literature). In the following, we assume that we have a sketch-orthonormal basis $V_m$ of $\K_m(A, \vec b)$ and an associated Hessenberg matrix $\underline{G}_m$ which satisfy the Arnoldi-like relation \cref{eqn:arnoldi-like-relation--sketched-basis}. We define $\beta \in \R$ according to the identity $\vec b = \beta V_m \vec e_1$. We mention that the sketched Galerkin and sketched Petrov-Galerkin conditions have been used in the context of model order reduction in \cite{BalabanovNouyI, BalabanovNouyII}.

\subsubsection{Randomized GMRES}

Randomized GMRES (rGMRES) seeks an approximate solution $\xsgmres_m \in \K_m(A, \vec b)$ such that its associated residual satisfies the \emph{sketched Petrov-Galerkin condition}
\begin{equation*}
	\Omega (\vec b - A \xsgmres_m) \perp \Omega A \K_m(A, \vec b).
\end{equation*}
Using \cref{eqn:arnoldi-like-relation--sketched-basis} and writing $\xsgmres_m = V_m \ysgmres_m$,  we have 
\begin{align*}
	0 = (\Omega A V_m)^T \Omega (\vec b - A V_m \ysgmres_m) &= (\Omega V_{m+1} \underline{G}_m)^T (\Omega \vec b - \Omega V_{m+1} \underline{G}_m \ysgmres_m) \\
	&= \underline{G}_m^T \beta \vec e_1 - \underline{G}_m^T \underline{G}_m \ysgmres_m.
\end{align*}
Similarly to the above derivation of GMRES, these are the normal equation associated with the least squares problem
\begin{equation}
	\label{eqn:linear-system-sketched-gmres-sol}
	\ysgmres_m = \argmin_{\vec y \in \R^m} \norm{\underline{G}_m \ysgmres_m - \beta \vec e_1}.
\end{equation}
From the condition $(\Omega A V_m)^T (\Omega \vec b - \Omega A V_m \ysgmres_m) = 0$, we also see that $\xsgmres_m$ solves the least squares problem
\begin{equation*}
	\xsgmres_m = \argmin_{\vec x \in \K_m(A, \vec b)} \norm{\Omega(\vec b -  A \vec x)}.
\end{equation*}
In other words, the rGMRES solution $\xsgmres_m$ minimizes the sketched residual $\Omega(\vec b - A \vec x)$ among all $\vec x$ in the Krylov subspace $\K_m(A, \vec b)$. If $\Omega$ is an $\epsilon$-embedding of $\K_{m+1}(A, \vec b)$, we get
\[ \norm{\vec b - A \xsgmres_m} \leq \sqrt{\frac{1+\epsilon}{1-\epsilon}} \argmin_{ \widehat{\vec x} \in \K_m(A, \vec b)} \|\vec b - A \widehat{\vec x}\|,\]
so rGMRES achieves a quasi-optimal residual, up to the multiplicative factor $(1+\epsilon)/(1-\epsilon)$.

We emphasize that, even though the sequence of sketched residuals of rGMRES is decreasing, it is not true in general that the sequence of residuals of rGMRES is also decreasing, especially in consecutive iterations where the sketched residual stagnates.

As in GMRES, we remark that the residual norm can be evaluated with cheap formulas in the Krylov basis:
\[ \| \vec b - A \xsgmres_m \| \leq \frac{1}{\sqrt{1- \epsilon}} \| \Omega (\vec b - A \xsgmres_m) \| = \frac{1}{\sqrt{1-\epsilon}}\|\Omega V_{m+1} (\beta \vec e_1 - \underline{G}_m \ysgmres_m)\| = \frac{1}{\sqrt{1-\epsilon}}\| \beta \vec e_1 - \underline{G}_m \ysgmres_m\|.\]
Of course, in finite-precision arithmetic, the two equalities that we used are only as good as the factorization error of the orthogonalization process used, and the true condition number of the computed basis. As to the first aspect, the factorization error of randomized QR processes is often better than that of deterministic processes (fewer flops). As to the second aspect, randomization often allow users to choose stabler methods. For these reasons, we often observe experimentally the true residual going slightly lower with rGMRES than with GMRES.

\subsubsection{Randomized FOM}

Randomized FOM (rFOM) searches for an approximate solution $\xsfom_m \in \K_m(A, \vec b)$ such that its residual satisfies the \emph{sketched Galerkin condition}
\begin{equation*}
	\Omega (\vec b - A \xsfom_m) \perp \Omega \K_m(A, \vec b).
\end{equation*}
Using \cref{eqn:arnoldi-like-relation--sketched-basis} and writing $\xsfom_m = V_m \ysfom_m$, we can rewrite this condition as 
\begin{align*}
	0 = (\Omega V_m)^T \Omega (\vec b - A V_m \ysfom_m) &= (\Omega V_m)^T (\Omega \vec b - (\Omega V_m G_m + g_{m+1,m} \Omega \vec q_{m+1} \vec e_m^T) \ysfom_m) \\
	&= \beta \vec e_1 - G_m \ysfom_m,
\end{align*}
where for the last equality we exploited the orthogonality of the columns of $\Omega Q_{m+1}$. This yields the expression for the randomized FOM approximate solution
\begin{equation}
	\label{eqn:linear-system-sketched-fom-sol}
	\xsfom_m = V_m \ysfom_m, \qquad \text{where} \qquad G_m \ysfom_m = \beta \vec e_1.
\end{equation}
Note that when $A$ is symmetric positive definite, the randomized FOM approximation \cref{eqn:linear-system-sketched-fom-sol} cannot in general be obtained with short-term recurrences, in contrast with the deterministic case, where FOM reduces to the CG method. Indeed, as mentioned already in \cref{subsec:rand-Arnoldi} for the Lanczos process, the core reason behind the short recurrence in CG is the fact that the Hessenberg matrix $H_m$ in the Arnoldi relation \cref{eqn:arnoldi-relation} is tridiagonal when $A$ is symmetric. 
On the other hand, when we generate a sketch-orthogonal basis $V_m$ through the randomized Arnoldi process we have the relation \cref{eqn:arnoldi-like-relation--sketched-basis}, which implies $G_m = (\Omega V_m)^T \Omega A V_m = V_m^T \Omega^T \Omega A V_m$, hence this matrix is not symmetric in general, unless $\Omega^T \Omega$ commutes with $A$, which is usually not the case. We refer to \cite[Section~4]{TGB23} for a more in-depth discussion.

It is shown in~\cite{Cullum1996relations} that the sequence of FOM approximants $\xfom_1, \xfom_2 \hdots $, for an arbitrary operator $A \in \R^{n \times n}$, yields a sequence of quasi-optimal residual norms $\|\vec b - A \xfom_1\|_2,  \|\vec b - A \xfom_2 \|_2, \dots$, with occasional spikes when the minimum residual sequence (produced by GMRES) stagnates. This property is founded on~\cite[Theorem 1]{Cullum1996relations}, which only uses the Hessenberg matrix.  As long as the randomized Arnoldi process does not break, this property extends to the sequence of residuals associated with randomized FOM approximations: 
\[ \forall \; 2 \leq k \leq m, \quad \| \Omega \vec r_{k-1}^{\text{G}} \| < \| \Omega \vec r_k^{\text{G}} \|,\]
\[ \| \vec r_k^{\text{F}} \| \leq \sqrt{\frac{1+\epsilon}{1-\epsilon}} \cdot \frac{1}{\sqrt{1 - (\|\Omega \vec r_k^{\text{G}}\| / \| \Omega \vec r_{k-1}^{\text{G}}\|)^2}} =: \sqrt{\frac{1+\epsilon}{1-\epsilon}} \cdot \alpha_k \| \Omega \vec r_k^{\text{G}}\|, \quad 2 \leq k \leq m.\]
This bound is valid for both the symmetric and non-symmetric cases $A\in\R^{n\times n}.$

In the case where $A \in \R^{n \times n}$ is symmetric and positive definite, FOM (which is equivalent to CG) produces a sequence that minimizes the $A$-norm of the error $\|\vec x - \xfom_k\|_A$. On the other hand, the behavior of $\| \vec x - \xsfom_k\|_A$ is a more complex topic. One of the main difficulties surrounding this question is that the standard sketching framework does not yield a simple concept of \textit{sketched energy norm}, as $(\vec x,\vec y) \mapsto (\Omega A \vec x)^T \Omega \vec y$ is not even symmetric in general. 

For inputs of medium difficulty $A \in \R^{n \times n}$ where $A$ is symmetric or non-symmetric, the cost of rRFOM might be higher than that of CG or BiCG. However, it is well known that these algorithms, based on short recurrences, lose orthogonality and hence on very-ill conditioned symmetric inputs FOM is often preferred to CG. In these cases, where we may require the Householder-Arnoldi iteration, randomized FOM based on RGS2-Arnoldi or RHQR-Arnoldi perform well, achieving the same stability as MGS2-Arnoldi and Householder-Arnoldi, for half the flops (or even a third) and much fewer synchronizations of a parallel computer. 

For both rGMRES and rFOM, the solution at the $m$-th iteration is given in the form $\vec x_m = V_m \vec y_m$, where $\vec y_m \in \R^m$. This implies that, if the sketch-orthogonal basis $V_m$ is obtained implicitly via a whitening procedure, we can compute the approximate solution as $\vec x_m = W_m (R_m^{-1} \vec y_m)$, where $W_m$ is a basis constructed, for instance, with the $k$-truncated Arnoldi process and $R_m$ is the upper triangular factor of the QR factorization of $\Omega W_m$. The Hessenberg matrix associated with $V_m$, required for the solution of \cref{eqn:linear-system-sketched-gmres-sol} and \cref{eqn:linear-system-sketched-fom-sol}, can be obtained from the whitened-sketched Arnoldi relation \cref{eqn:whitened-sketched-arnoldi-relation}.

%% file: eigenvalues.tex
\section{Solution of eigenvalue problems}
\label{sec:eigenvalues}

Finding the eigenvalues and eigenvectors of a matrix $A \in \mathbb{R}^{n \times n}$ is a fundamental task in numerical linear algebra, with applications ranging from structural engineering, where eigenvectors represent natural modes of vibration whose frequencies are given by the eigenvalues, to quantum chemistry, where the lowest eigenvector corresponds to the ground state of a molecule, and to machine learning, where principal component analysis reduces dimensionality by focusing on dominant eigendirections. \cite{Pain2005PhysicsVibrationsWaves,Jolliffe2016Principalcomponentanalysis}

The standard eigenvalue problem consists in finding pairs $(\lambda_i, \vec x_i)$ indexed by $i \in \mathcal{I}$ of scalars $\lambda_i \in \mathbb{C}$ and unit‑norm eigenvectors $\vec x_i \in \mathbb{C}^n$ such that
\begin{equation}
    A \vec x_i = \lambda_i \vec x_i.
\end{equation}

Since eigenvalues are the roots of the characteristic polynomial of $A$, the Abel–Ruffini theorem implies that there is no general algebraic formula in radicals for polynomials of degree greater than four. Consequently, for matrices of dimension $n \ge 5$ the eigenvalue problem is solved in practice by iterative numerical methods that compute approximate eigenpairs \cite{saad2011numerical, Kressner2005NumericalMethodsGeneral}.

The set $\mathcal{I}$ of eigenpairs that can be computed in practice depends on the size of the matrix and the structure of the problem. For moderate $n$, the state‑of‑the‑art method to compute all eigenpairs $\mathcal{I}=\{1,\dots,n\}$ is the shifted QR algorithm. This algorithm produces a sequence of matrices $A_k$ similar to $A$ that converges (up to ordering) to the real Schur form of $A$; in particular, one obtains an upper quasi‑triangular matrix with eigenvalues on the diagonal (and $2\times 2$ blocks for complex conjugate pairs):
\begin{equation}
    \lim_{k \to \infty} A_k = \begin{bmatrix}
        \lambda_1 & & \star  \\
        & \ddots & \\
        0 &  & \lambda_n 
    \end{bmatrix}.
\end{equation}
In step $k$, approximate eigenvalues are read from the diagonal (or block eigenvalues) of $A_k$, and eigenvectors can be recovered from accumulated similarity transformations; convergence is typically monitored by the magnitudes of the off‑diagonal entries. See \cite{Miminis1991ImplicitShiftingQR,Watkins1982UnderstandindQRAlgorithm,Kressner2005NumericalMethodsGeneral} for details on the shifted QR algorithm and \cite{Golub2013Matrixcomputations} for the real Schur form. Each step of the QR algorithm consists of computing the QR factorization of $A_k$, resulting in an arithmetic cost of $O(n^3)$ for the method, which makes it impractical for large‑scale eigenvalue problems.

\subsection{Rayleigh-Ritz}
For large eigenvalue problems, a small subset $\mathcal{I}\subset\{1,\dots,n\}$ of $m$ eigenpairs is generally sought, with $m$ small relative to $n$. The Rayleigh-Ritz method projects $A$ onto an $m$‑dimensional subspace and solves the reduced eigenvalue problem to obtain approximations to the desired eigenpairs; its convergence and accuracy depend on the choice of the projection subspace and on the spectral properties of $A$. More precisely, given an $m$-dimensional subspace $\mathcal{K}_m$, the Rayleigh-Ritz method seeks an approximate eigenvector $\widetilde{\vec x}$ and an approximate eigenvalue $\widetilde{\lambda}$ by imposing the following two constraints:
\begin{enumerate}
    \item The approximate eigenvector (Ritz vector) $\widetilde{\vec x}$ belongs to $\mathcal{K}_m$.
    \item The residual vector $A \widetilde{\vec x} - \widetilde{\lambda} \widetilde{\vec x}$ is orthogonal to $\mathcal{K}_m$.
\end{enumerate}
The orthogonality of the residual to $\mathcal{K}_m$ fixes the $m$ degrees of freedom that arise when seeking $\widetilde{\vec{x}}$ in $\mathcal{K}_m$ and is known as the Galerkin condition. Let $Q_m \in \mathbb{R}^{n \times m}$ be an orthogonal basis for $\mathcal{K}_m$ and write $\widetilde{\vec{x}} = Q_m \vec y \in \mathcal{K}_m$. The Galerkin condition can be written as
\begin{align}
    Q_m^T(A Q_m \vec y - \widetilde{\lambda} Q_m \vec y) = 0, \\
    Q_m^T A Q_m \vec y = \widetilde{\lambda} \vec y. \label{eq:RR_exactEVP}
\end{align}
Hence, the pair $(\widetilde{\lambda},\vec y)$ is an exact eigenpair of the small operator $H_m = Q_m^T A Q_m \in \mathbb{R}^{m \times m}$, which represents the orthogonal projection of $A$ onto $\mathcal{K}_m$. Multiplying \cref{eq:RR_exactEVP} by $Q_m$ and using $\vec y = Q_m^T Q_m \vec y$ gives
\begin{equation}
    \proj{\mathcal{K}_m} A \proj{\mathcal{K}_m} \widetilde{\vec{x}} = \widetilde{\lambda} \widetilde{\vec{x}},
\end{equation}
where $\proj{\mathcal{K}_m}$ is the orthogonal projector onto $\mathcal{K}_m$ represented by $Q_m Q_m^T$. This reinforces the characterization of the Rayleigh-Ritz method as an orthogonal projection approach that solves a small eigenvalue problem for $H_m$ to obtain approximate eigenpairs of $A$.

The quality of a Ritz pair, namely the Ritz eigenvalue $\widetilde{\lambda}$ and the Ritz vector $\widetilde{\vec{x}}$, as an approximation of an eigenpair of $A$ depends strongly on the subspace $\mathcal{K}_m$, in particular on the distance between an exact eigenvector $\vec x$ of $A$ and $\mathcal{K}_m$. In practice, Krylov subspaces are a natural and effective choice for $\mathcal{K}_m$. The Arnoldi procedure introduced in \cref{sec:krylov} constructs, from a starting vector $\vec b$, an orthogonal basis $Q_m$ for the Krylov subspace $\K_m(A ,\vec b)$ and simultaneously computes the small projected Hessenberg matrix $H_m$, as given by \cref{eqn:arnoldi-relation}.

Since $H_m$ is an $m\times m$ matrix, it is possible to compute its eigenpairs $(\widetilde{\lambda_i}, \vec y_i)$ (for instance, with a shifted QR algorithm) and form Ritz vectors $\widetilde{\vec{x}}_i = Q_m \vec y_i$, for $i = 1, \dots m$. Multiplying \cref{eqn:arnoldi-relation} by $\vec{y}_i$ yields
\begin{equation}
    A \widetilde{\vec{x}}_i - \widetilde{\lambda}_i \widetilde{\vec{x}}_i =  h_{m+1,m} \vec q_{m+1} \vec{e}_m^T \vec y_i.
\end{equation} 
Consequently, the residual norm satisfies $\norm{A \widetilde{\vec{x}}_i - \widetilde{\lambda}_i \widetilde{\vec{x}}_i} = \norm{h_{m+1,m} \vec q_{m+1} \vec{e}_m^T \vec y_i}$. Using $\norm{\vec q_{m+1}} = 1$, one obtains a simple and inexpensive expression to monitor the quality of the approximate $i$‑th eigenpair at the end of the Arnoldi method:
\begin{equation}
    \norm{A \widetilde{\vec{x}}_i - \widetilde{\lambda}_i \widetilde{\vec{x}}_i} =  \abs{h_{m+1,m}} \cdot \abs{\vec{e}_m^T \vec y_i} \quad \text{for $i=1,\dots,m$.}
\end{equation}

The Arnoldi method \cite{Arnoldi1951principleminimizediterations}, originally developed as an extension of the Lanczos method for non-symmetric matrices \cite{Lanczos1950iterationmethodsolution}, underpins many restarted schemes with filtering, such as the implicitly restarted Arnoldi algorithm and the Krylov-Schur method. The convergence results and the bounds on the distance between the eigenvectors of $A$ and the Krylov subspace $\K_m(A,\vec b)$ are given in particular in \cite{saad2011numerical,Bellalij2010FurtherAnalysisArnoldi,Bellalij2016distanceeigenvectorKrylov}.

\subsection{Randomized Rayleigh-Ritz}
\label{subsec:sketch-RR}

In the Arnoldi algorithm, each step consists in applying $A$ to the last vector of the Krylov basis and orthogonalizing the result against the existing basis. Provided that $A$ is structured so that the cost of a matrix-vector product scales linearly with $n$, as is the case for a sparse $A$, the dominant arithmetic cost of the method is the orthogonalization step. Reducing this cost with randomized orthogonalization procedures has been the subject of recent work, including \cite{BalabanovGrigori25block,NakatsukasaTropp24,GuettelSimunec24,DedamasGrigori25ira}. The approach in \cite{BalabanovGrigori25block} uses a randomized (block) Gram–Schmidt process \cite{BalabanovGrigori2022} within the Arnoldi procedure to obtain a sketch‑orthonormal basis for the Krylov subspace, as described in \cref{subsec:rand-Arnoldi} with \cref{algo:rando-arnoldi}.

Let $V_m$ and $G_m$ denote the sketch-orthogonal basis and the associated Hessenberg matrix generated by the randomized Arnoldi process \cref{algo:rando-arnoldi}.
The structure of the resulting randomized Arnoldi decomposition \cref{eqn:arnoldi-like-relation--sketched-basis} mirrors the deterministic Arnoldi factorization \cref{eqn:arnoldi-relation}, and therefore a randomized Rayleigh–Ritz procedure follows naturally. The theory relies on the oblique projector $\skproj{\K_m}$ defined in \cite{BalabanovGrigori25block,DedamasGrigori25ira} by
\begin{equation}
    \skproj{\K_m} \vec x = \argmin_{\vec v \in \K_m} \norm{\Omega(\vec x - \vec v)}
\end{equation}
for $\vec x \in \R^n$, which we introduced in \cref{subsec:sketch-orth-discussion}. If $(\widetilde{\lambda},\vec y)$ is an eigenpair of $G_m$, multiplying \cref{eqn:arnoldi-like-relation--sketched-basis} by $\vec y$ gives
\begin{equation}
    \label{eqn:sketch-residual-Arnoldi}
    A (V_m \vec y) - \widetilde{\lambda} (V_m \vec y) = g_{m+1,m} \vec v_{m+1} \vec e_m^T \vec y.
\end{equation}
Defining the Ritz vector $\widetilde{\vec{x}} = V_m \vec y \in \K_m$, the residual $A \widetilde{\vec{x}} - \widetilde{\lambda} \widetilde{\vec{x}}$ is sketch‑orthogonal to $\K_m(A, \vec b)$ because it is proportional to $\vec v_{m+1}$:
\begin{equation}
    \label{eqn:sketch-Galerkin}
    \Omega (A \widetilde{\vec{x}} - \widetilde{\lambda} \widetilde{\vec{x}}) \perp \Omega \K_m.
\end{equation}
Equation \cref{eqn:sketch-Galerkin} is the sketched Galerkin condition. Since $\widetilde{\vec{x}} \in \K_m(A, \vec b)$, this framework is naturally called a sketched or randomized Rayleigh–Ritz approach. It is shown in \cite{BalabanovGrigori25block} that a pair satisfying \cref{eqn:sketch-Galerkin} is an exact eigenpair of $\skproj{\K_m} A \skproj{\K_m}$:
\begin{equation}
    \skproj{\K_m} A \skproj{\K_m} \widetilde{\vec{x}} = \widetilde{\lambda} \widetilde{\vec{x}}.
\end{equation}
Since $G_m = (\Omega V_m)^T \Omega A V_m$, \cite{DedamasGrigori25ira} further shows that $G_m$ is the representation of $\skproj{\K_m} A \skproj{\K_m}$ in the basis $V_m$, i.e. $\skproj{\K_m} A \skproj{\K_m} \vec w = V_m G_m \vec z$ when $\vec w = V_m \vec z$. In particular, the characteristic polynomial $p_{G_m}$ of $G_m$ verifies that
\begin{equation}
    p_{G_m} = \argmin_{p \in \mathcal{PM}_m} \norm{\Omega p(A) \vec b}
\end{equation}
where $ \mathcal{PM}_m$ is the set of monic polynomials of degree $m$. If $p_{H_m}$ is the characteristic polynomial of the Hessenberg matrix $H_m$ that originates from an orthogonal projection, then 
\begin{equation}
    \norm{p_{G_m}(A) \vec b} \leq \sqrt{\frac{1+\epsilon}{1-\epsilon}}  \norm{p_{H_m}(A) \vec b} = \sqrt{\frac{1+\epsilon}{1-\epsilon}} \min_{p \in \mathcal{PM}_m} \norm{p(A) \vec b},
\end{equation}
which follows from the $\epsilon$-embedding property \cref{eq:epsembedding}.
As in the deterministic setting, the quality of a Ritz pair $(\widetilde{\lambda},\widetilde{\vec{x}})$ depends on the spectral properties of $A$ and on the Krylov subspace $\K_m(A,\vec b)$. Bounds on the residual $\norm{(A - \widetilde{\lambda} I) \widetilde{\vec{x}}}$ and on the distance of an exact eigenvector $\vec{x}$ from $\K_m$ are derived in \cite{BalabanovGrigori25block,DedamasGrigori25ira}; these results are outside the scope of this review but are useful to characterize the convergence of randomized Arnoldi. The main conclusion is that the randomized Arnoldi method is an oblique projection technique on a Krylov subspace, delivering accuracy comparable to deterministic Arnoldi while reducing cost via sketch‑orthogonalization.

We conclude this section by describing practical convergence monitoring. Given an eigenpair $(\widetilde{\lambda}_i,\vec y_i)$ of $G_m$, \cref{eqn:sketch-residual-Arnoldi} implies
\begin{equation}
    \norm{A \widetilde{\vec{x}}_i - \widetilde{\lambda}_i \widetilde{\vec{x}}_i} = \abs{g_{m+1,m}} \cdot \norm{\vec v_{m+1}} \cdot \abs{\vec e_m^T \vec y_i}.
\end{equation}
Unlike the deterministic case, $\norm{\vec v_{m+1}}$ need not be equal to one. Instead, \cite{NakatsukasaTropp24} provides the bound, for $i=1,\dots,m$,
\begin{equation}
    \label{eqn:sketch-Arnoldi-residual-error}
    \sqrt{\frac{1 - \epsilon}{1 + \epsilon}}  \norm{ \Omega (A \widetilde{\vec{x}}_i - \widetilde{\lambda}_i \widetilde{\vec{x}}_i)} \leq \frac{\norm{A \widetilde{\vec{x}}_i - \widetilde{\lambda}_i \widetilde{\vec{x}}_i}}{\norm{\widetilde{\vec{x}}_i}} \leq \sqrt{\frac{1 + \epsilon}{1 - \epsilon}}  \norm{ \Omega (A \widetilde{\vec{x}}_i - \widetilde{\lambda}_i \widetilde{\vec{x}}_i)}.
\end{equation}
Since $\norm{ \Omega (A \widetilde{\vec{x}}_i - \widetilde{\lambda}_i \widetilde{\vec{x}}_i)} = \abs{g_{m+1,m}} \cdot \abs{\vec e_m^T \vec{y}_i}$, this quantity is readily available during the randomized Arnoldi iteration and provides a good approximation of the relative residual $\norm{A \widetilde{\vec{x}}_i - \widetilde{\lambda}_i \widetilde{\vec{x}}_i} / \norm{\widetilde{\vec{x}}_i}$ in view of \cref{eqn:sketch-Arnoldi-residual-error}.

\subsection{Restarting strategies for Krylov subspace methods}

A typical issue arising from Krylov subspace methods such as the Arnoldi procedure is the quadratic growth, with the number of steps $m$, of the arithmetic cost of orthogonalization together with the accompanying memory required to store a growing Krylov basis $V_m$ of vectors in $\R^n$. To mitigate these problems, several restarting strategies have been proposed to produce a new Arnoldi factorization from a new starting vector $\vec{b}^+$ that retains much of the information from a previous length‑$m$ factorization. Popular methods include the implicitly restarted Arnoldi (IRA) and the Krylov–Schur (KS) approaches \cite{Sorensen1992ImplicitApplicationPolynomial,Lehoucq1995DeflationTechniquesImplicitly,Stewart2002KrylovSchurAlgorithm}. These approaches have been extended to the randomized Arnoldi setting in \cite{DedamasGrigori25ira,DedamasGrigori25ks}.

Assume we have a randomized Arnoldi decomposition \cref{eqn:arnoldi-like-relation--sketched-basis} with Krylov basis $\K_m(A,\vec b)$ and we wish to restart it. The randomized implicitly restarted Arnoldi method (rIRA) \cite{DedamasGrigori25ira} is based on polynomial filtering; see \cite[Chapter 7]{saad2011numerical}. If the starting vector $\vec b \in \R^n$ admits the expansion
\begin{equation}
    \vec b = \sum_{i=1}^{n} \alpha_i \vec x_i,
\end{equation}
with $\vec x_i$ the eigenvectors of $A$, then for any polynomial $p$ one has
\begin{equation}
    p(A) \vec b = \sum_{i=1}^{n} \alpha_i p(\lambda_i) \vec x_i.
\end{equation}
If the goal is to approximate the $k$ dominant eigenvectors $\vec x_1,\dots,\vec x_k$ of $A$ (with $k \le m$), one may use $\vec b^+ = p(A)\vec b$ as a new starting vector, where $p$ is chosen to be large on $\lambda_1,\dots,\lambda_k$ and small on the remaining eigenvalues. Then $\K_m(A,\vec b^+)$ is likely to contain high‑quality Ritz vectors for the desired eigenpairs. The key advantage of IRA is that $p(A)\vec b$ can be applied implicitly by performing shifted QR steps on the small Hessenberg matrix $G_m$ obtained from the previous Arnoldi factorization; this property carries over to rIRA. In practice the polynomial $p$ has degree $q$, equal to the number of shifted QR steps performed, and its roots are the shifts. There exist different strategies to choose these shifts (for example Chebyshev polynomials \cite[Chapter 4]{Sorensen1992ImplicitApplicationPolynomial}), but a common and successful approach is to pick the $q = m-k$ unwanted Ritz values of $G_m$ as the shifts. When aiming for the $k$ largest eigenvalues of $A$, one typically designates the $q$ Ritz values of smallest modulus as unwanted. This promotes the Ritz vectors associated with the largest Ritz values in the subsequent Arnoldi iteration and yields practical convergence to the desired eigenvalues. Convergence is proven in \cite{DedamasGrigori25ira} in a more restrictive setting of fixed shifts over the iteration. It is also shown how sketch‑orthonormalization preserves IRA’s beneficial properties while reducing computational cost.

\begin{figure}[t]
\centering
\resizebox{1.0\textwidth}{!}{%
\begin{circuitikz}
\tikzstyle{every node}=[font=\LARGE]
\draw  (-2.5,12) rectangle (0,7);
\node [font=\LARGE] at (-4,9.5) {$A$};
\node [font=\LARGE] at (0.75,9.5) {$=$};
\draw  (1.75,12) rectangle (4.25,7);
\draw [short] (5,12) -- (7.5,12);
\draw [short] (7.5,12) -- (7.5,9.5);
\draw [short] (5,12) -- (7.5,9.5);
\draw [short] (5,12) -- (5,11.5);
\draw [short] (5,11.5) -- (7,9.5);
\draw [short] (7,9.5) -- (7.5,9.5);
\node [font=\LARGE] at (9.75,9.5) {$ + g_{m+1,m} \vec v_{m+1} \vec e_{m}^T$};
\draw [->, >=Stealth] (12.5,9.5) -- (22.75,9.5);
\draw [short] (25,10.75) -- (27.5,10.75);
\draw [short] (27.5,10.75) -- (27.5,8.25);
\draw [short] (25,10.75) -- (27.5,8.25);
\draw [short] (25,10.75) -- (25,10.25);
\draw [short] (25,10.25) -- (27,8.25);
\draw [short] (27,8.25) -- (27.5,8.25);
\draw [->, >=Stealth, dashed] (30.25,8.25) .. controls (28.75,6.25) and (27.75,6.25) .. (26.25,8.25);
\draw [->, >=Stealth, dashed] (26.25,11.5) .. controls (27.75,13.25) and (28.5,13.25) .. (30,11.5);
\node [font=\LARGE] at (30,9.5) {$ - \; \mu_i I = Q^i R^i$};
\node [font=\LARGE] at (31,7.0) {$i = 1,\dots,q$};
\draw [->, >=Stealth] (27.5,6) -- (27.5,1.75);
\node [font=\LARGE] at (-1.5,13) {$V_m$};
\node [font=\LARGE] at (2.75,13) {$V_m$};
\node [font=\LARGE] at (6.5,13) {$G_m$};
\draw  (16.25,-1) rectangle (18.75,-6);
\node [font=\LARGE] at (14.75,-3.5) {A};
\node [font=\LARGE] at (22,-3.5) {=};
\draw  (22.5,-1) rectangle (25,-6);
\draw [short] (31.25,-1) -- (33.75,-1);
\draw [short] (33.75,-1) -- (33.75,-3.5);
\draw [short] (31.25,-1) -- (33.75,-3.5);
\draw [short] (31.25,-1) -- (31.25,-1.5);
\draw [short] (31.25,-1.5) -- (33.25,-3.5);
\draw [short] (33.25,-3.5) -- (33.75,-3.5);
\node [font=\LARGE] at (39.50,-3.5) {$ +g_{m+1,m} \vec v_{m+1} \vec e_{m}^T Q$};
\node [font=\LARGE] at (17.25,0) {$V_m$};
\node [font=\LARGE] at (23.5,0) {$V_m$};
\node [font=\LARGE] at (32.5,0) {$G_m$};
\draw  (25.5,-1) rectangle (27.75,-3.5);
\draw  (28.5,-1) rectangle (30.75,-3.5);
\draw  (34.5,-1) rectangle (36.75,-3.5);
\draw  (19.25,-1) rectangle (21.5,-3.5);
\node [font=\LARGE] at (20.5,0) {$Q$};
\node [font=\LARGE] at (26.5,0) {$Q$};
\node [font=\LARGE] at (29.5,0) {$Q^T$};
\node [font=\LARGE] at (35.5,0) {$Q$};
\node [font=\LARGE] at (17.5,10.5) {Run the shifted QR algorithm on $G_m$};
\node [font=\LARGE] at (23.5,5.5) {Multiply the factorization};
\draw [->, >=Stealth] (14.5,-1.25) -- (10.25,-1.25);
\draw [dashed] (16.25,-6.5) .. controls (17.75,-7.25) and (19.75,-7.25) .. (21.25,-6.5);
\draw [dashed] (22.75,-6.5) .. controls (24.25,-7.25) and (26.25,-7.25) .. (27.75,-6.5);
\draw [dashed] (30.25,-4) .. controls (31.75,-4.75) and (33.75,-4.75) .. (35.25,-4);
\node [font=\LARGE] at (18.5,-7.5) {$\widetilde{V}_m$};
\node [font=\LARGE] at (25,-7.5) {$\widetilde{V}_m$};
\node [font=\LARGE] at (32.5,-5) {$\widetilde{G}_m$};
\draw  (-1.25,-1) rectangle (0,-6);
\node [font=\LARGE] at (-2.75,-3.5) {$A$};
\node [font=\LARGE] at (2,-3.5) {$=$};
\draw  (3,-1) rectangle (4.25,-6);
\draw [short] (6.25,-1) -- (8.75,-1);
\draw [short] (8.75,-1) -- (8.75,-3.5);
\draw [short] (6.25,-1) -- (8.75,-3.5);
\draw [short] (6.25,-1) -- (6.25,-1.5);
\draw [short] (6.25,-1.5) -- (8.25,-3.5);
\draw [short] (8.25,-3.5) -- (8.75,-3.5);
\node [font=\LARGE] at (10.5,-3.5) {$ + \widetilde{\vec{r}}_k \vec e_k^T$};
\node [font=\LARGE] at (-0.50,0.5) {$\widetilde{V}_m(:,1:k)$};
\node [font=\LARGE] at (3.50,0.5) {$\widetilde{V}_m(:,1:k)$};
\node [font=\LARGE] at (7.0,0) {$\widetilde{G}_m(1:k,1:k)$};
\draw [dashed] (0,0) -- (0,-6.75);
\draw [dashed] (4.25,-0.25) -- (4.25,-6.75);
\draw [ dashed] (6,-0.5) rectangle  (7.5,-2.5);
\draw [dashed] (0,-1) -- (1.25,-1);
\draw [dashed] (1.25,-1) -- (1.25,-6);
\draw [dashed] (1.25,-6) -- (0,-6);
\draw [dashed] (4.25,-1) -- (5.25,-1);
\draw [dashed] (5.25,-1) -- (5.25,-6);
\draw [dashed] (5.25,-6) -- (4.25,-6);
\draw [->, >=Stealth] (2,1.5) -- (2,6);
\node [font=\LARGE] at (24.25,4.5) {by the accumulated};
\node [font=\LARGE] at (24,3.25) {$Q = Q^1 \dots Q^q$ factor};
\node [font=\LARGE] at (12.75,0.25) {Truncate to $k$ first columns};
\node [font=\LARGE] at (3.75,5) {Extend the};
\node [font=\LARGE] at (5,4) {randomized Arnoldi};
\node [font=\LARGE] at (4,3) {factorization};
\end{circuitikz}
}%
\caption{Visual representation of a cycle of the rIRA method.}
\label{fig:rIRA_schema}
\end{figure}
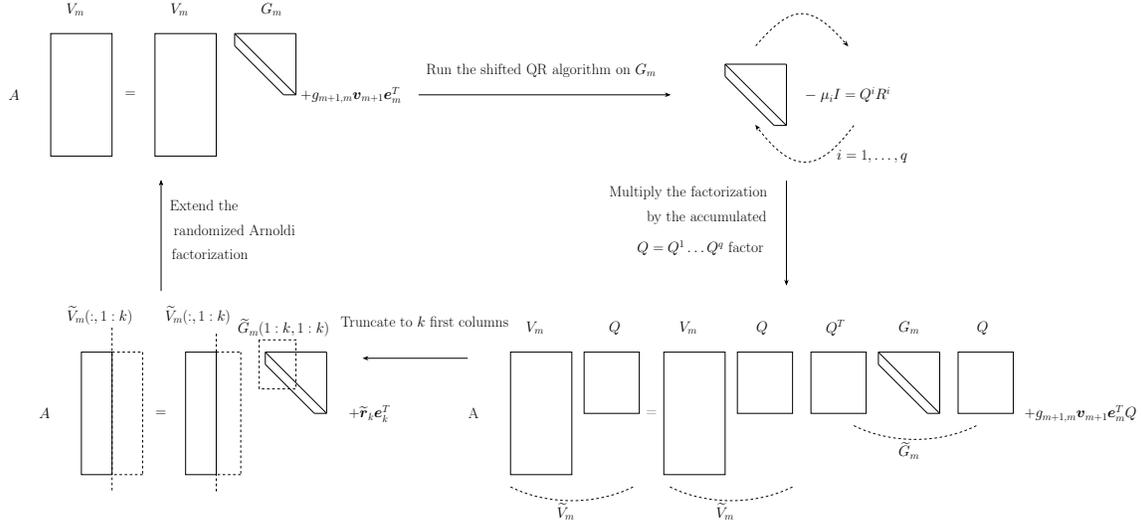

We summarize the main steps of the randomized Implicitly Restarted Arnoldi algorithm below. The method is illustrated in \cref{fig:rIRA_schema} and further details appear in \cite{DedamasGrigori25ira}.
\begin{enumerate}
    \item Compute the eigenvalues $\mu_1,\dots,\mu_m$ of $G_m$ and select $k$ wanted eigenvalues among them (e.g. those of largest or smallest modulus).
    \item Apply $q = m - k$ steps of the shifted QR algorithm on $G_m$ using the $q$ unwanted eigenvalues as shifts.
    \item Recover a new length‑$k$ randomized Arnoldi factorization by multiplying the previous factorization by the accumulated orthogonal transformations from the shifted QR steps and truncating to length $k$.
    \item Extend this factorization to length $m$ using the randomized Arnoldi method and repeat from step 1.
\end{enumerate}

The IRA method has been implemented in the ARPACK library \cite{Lehoucq1998ARPACKusersguide} and has seen wide adoption. However, \cite{Lehoucq1995DeflationTechniquesImplicitly,Stewart2002KrylovSchurAlgorithm} observed that the shifted QR algorithm can suffer from loss of forward stability, which motivated the development of the Krylov–Schur method \cite{Stewart2002KrylovSchurAlgorithm}. 

The Krylov–Schur method is based on a generalization of the Arnoldi decomposition, referred to as the Krylov decomposition:
\begin{equation}
    \label{eqn:Krylov-decomposition}
    A W_m = W_m B_m + \vec w_{m+1} \vec z_m^T.
\end{equation}
The distinction from an Arnoldi decomposition \cref{eqn:arnoldi-relation} is that the columns of $[W_m \; \vec w_{m+1}] \in \R^{n\times(m+1)}$ are only required to be linearly independent rather than orthonormal, $B_m$ must only be full‑rank and $\vec z_m$ may be arbitrary. This relaxation removes structural constraints present in Arnoldi factorizations; indeed, both Arnoldi factorization \cref{eqn:arnoldi-relation} and randomized Arnoldi factorization \cref{eqn:arnoldi-like-relation--sketched-basis} are special cases of Krylov decompositions. Working with Krylov decompositions permits replacing the shifted QR steps used in IRA by numerically stable Schur form re-orderings; see \cite{Kressner2006Blockalgorithmsreordering} for details. This idea is the basis for Stewart’s Krylov–Schur algorithm.

In \cite{DedamasGrigori25ks} it is shown that a randomized Arnoldi decomposition can be obtained from a Krylov decomposition and that a sketch‑orthonormal Krylov basis integrates naturally into the Krylov–Schur framework. Sketch‑orthonormalizing $W_m$ and translating $\vec w_{m+1}$ so that it is sketch‑orthogonal to the basis yields a randomized Arnoldi factorization. This observation leads to a randomized Krylov–Schur (rKS) algorithm that combines the stability of Schur reordering with the efficiency and scalability of sketch‑orthonormalization.

The Krylov–Schur method also incorporates a simple deflation procedure \cite{Stewart2002KrylovSchurAlgorithm,Kressner2005NumericalMethodsGeneral}, which was extended in rKS in \cite{DedamasGrigori25ks}. When eigenpairs have converged, they can be removed from the active subspace by sketch‑orthogonalizing the remaining Krylov vectors against the converged vectors, producing a partial sketch‑orthonormal Schur factorization for $A$. That is, $A Q_m = Q_m T_m$ where $Q_m$ is sketch‑orthonormal and $T_m$ is block upper‑triangular. If among the $k$ sought eigenpairs, there are $q$ converged eigenpairs whose residual errors are smaller than a value $\eta$, it is shown in \cite{DedamasGrigori25ks} that the sketch-orthonormal deflation procedure is equivalent to continuing the rKS method by seeking $k -q$ eigenpairs of a slightly perturbed matrix $A+E$ with 
\begin{equation}
    \norm{E}_{F,2} \leq \sqrt{q} \sqrt{\frac{1+\epsilon}{1-\epsilon}} \eta.
\end{equation}

We conclude this section with numerical experiments for the rKS method on its efficiency and accuracy, that are similar to those in \cite{DedamasGrigori25ks}. A number of $k = 40$ eigenvalues  of tri-diagonal synthetic matrices are sought, with two types of spectra, harmonic and geometric. This means that their diagonal entries are $1 + 1/i^2$ and $0.99^i$ for $i=1,\dots,n$, respectively. Their off-diagonal entries are Gaussian noise:
\begin{equation}
    \label{eq:synthMathSubdiags}
    \quad A_{i+1,i} = \frac{g_{i+}}{100}, \; A_{i-1,i} = \frac{g_{i-}}{100},
\end{equation}
where $g_{i \pm}$ are drawn from $\mathcal{N}(0,1)$. The eigenvalues are sought with a precision of $10^{-10}$ on the sketch residuals, that is $\norm{\Omega(A \widetilde{\vec{x}}_i  - \widetilde{\lambda} \widetilde{\vec{x}}_i) } \leq 10^{-10}$ from \cref{eqn:sketch-Arnoldi-residual-error}. The Ritz vectors are sought in Krylov subspaces of dimension $m = 2k = 80$. For the  harmonic spectrum, the $40$ eigenvalues of smallest modulus (SM) are sought, whereas for the geometric spectrum it is the ones of largest modulus (LM). Experiments are conducted with the Julia programming language in its version 1.10 \cite{Bezanson2017Juliafreshapproach}. 

In \cref{fig:rKS-exectime} the input dimension $n$ of $A \in \R^{n \times n}$ increases from $n_{min} = 10^5$ to $n_{max} = 5 \times 10^6$ and the execution times of IRA, KS and rKS are compared. It can be observed that the rKS method runs faster than IRA and KS, with a speed-up of 2-3 times faster thanks to sketch-orthonormalization. In \cref{fig:rKS-accuracy}, the quality of the obtained Ritz eigenpairs is compared. The left panel shows the evolution of the maximum over $i = 1,\dots,40$ of $\norm{\Omega(A \widetilde{\vec{x}}_i - \widetilde{\lambda} \widetilde{\vec{x}}_i) }$ at the end of each restart for KS and rKS, demonstrating that both methods converge in roughly the same number of iterations with similar convergence behavior. The right panel displays the real parts of the obtained eigenvalues, using the IRA method (via Julia's \texttt{eigs} function from ARPACK) as a reference. All three methods find the same approximate eigenvalues, demonstrating that the faster randomized approach rKS delivers accurate solutions for these problems.
\begin{figure}[t]
    \centering
    \input{figures/mosaics_synth_timings.pgf}
    \caption{Execution time of KS and rKS for an increasing input dimension $n$. A total of $40$ eigenvalues are sought, those of smallest modulus for the harmonic spectrum (left) and those of largest modulus for the geometric spectrum (right). The time label is in logarithmic scale, meaning that rKS being constantly below KS and IRA represents here a speed-up of around 3 times faster.}
    \label{fig:rKS-exectime}
\end{figure}
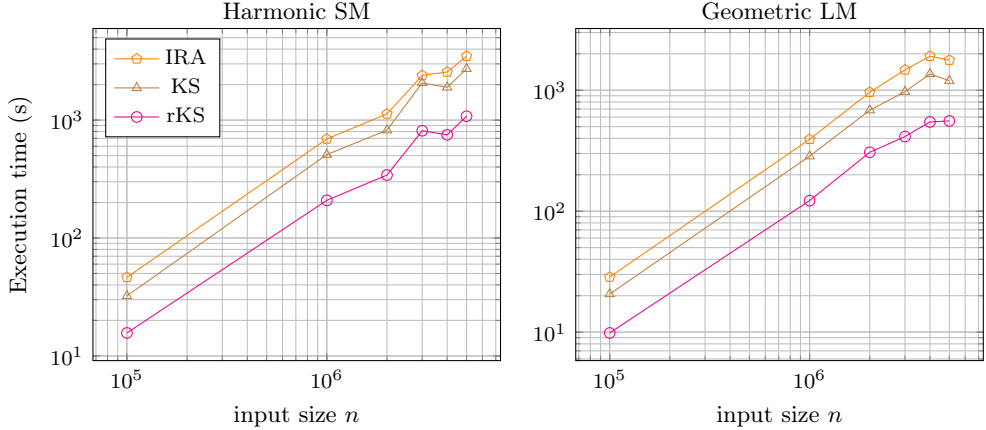
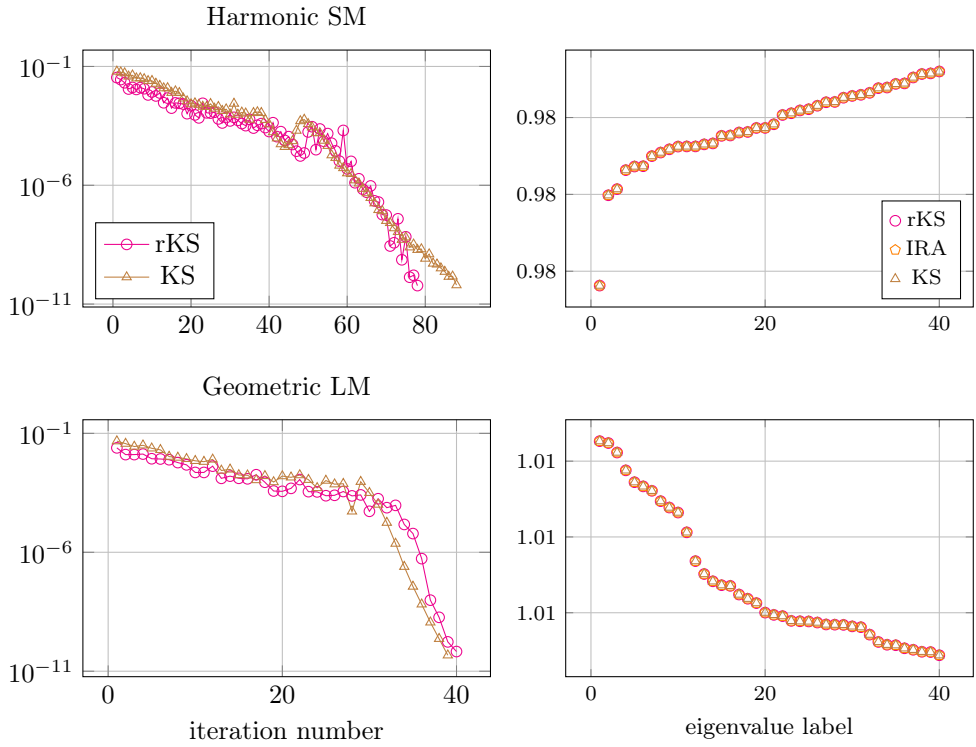
\begin{figure}
    \centering
    \input{figures/reserr_vp_noth.pgf}
    \caption{Left: maximum over $i$ of the residual errors $\norm{\Omega(A \tilde{\vec{x}}_i  - \tilde{\lambda} \tilde{\vec{x}}_i) }$, against the number of iterations of KS and rKS. Convergence is declared when this maximum reaches $10^{-10}$. Right: obtained Ritz eigenvalues for KS and rKS compared to the reference IRA. }
    \label{fig:rKS-accuracy}
\end{figure}

%% file: figures/mosaics_synth_timings.pgf
\begin{tikzpicture}
\begin{groupplot}[group style={group size={2 by 1}, ylabels at={edge left}, xlabels at={edge bottom}}, ylabel={Execution time (s)}, xlabel={input size $n$}]
       \nextgroupplot[footnotesize, title={Harmonic SM}, height={6cm}, width={7cm}, ticks={both}, grid={both}, legend pos={north west}, max space between ticks={45}, x tick scale label style={yshift=10pt,xshift=6pt}, y tick scale label style={yshift=0pt,xshift=-18pt}, ymode=log, xmode=log]
    \addplot[color={orange}, mark={pentagon}]
        coordinates {
            (100000,46.403687034)
            (1000000,691.71289998)
            (2000000,1128.963799867)
            (3000000,2393.482875992)
            (4000000,2559.684628228)
            (5000000,3485.233101315)
        }
        ;
        \addlegendentry {IRA}
    \addplot[color={brown}, mark={triangle}]
        coordinates {
            (100000,32.268412087)
            (1000000,509.460285317)
            (2000000,821.587366933)
            (3000000,2066.506102456)
            (4000000,1895.491830494)
            (5000000,2721.227717158)
        }
        ;
        \addlegendentry {KS}
    \addplot[color={magenta}, mark={o}]
        coordinates {
            (100000,15.655238056)
            (1000000,208.609111025)
            (2000000,341.857560481)
            (3000000,812.766002839)
            (4000000,750.328983843)
            (5000000,1077.835513801)
        }
        ;
        \addlegendentry {rKS}
    \nextgroupplot[footnotesize, title={Geometric LM}, height={6cm}, width={7cm}, ticks={both}, grid={both}, legend pos={north west}, max space between ticks={45}, x tick scale label style={yshift=10pt,xshift=6pt}, y tick scale label style={yshift=0pt,xshift=-18pt}, ymode=log, xmode=log]
    \addplot[color={orange}, mark={pentagon}]
        coordinates {
            (100000,28.544154381)
            (1000000,393.375947838)
            (2000000,962.641672298)
            (3000000,1473.260071596)
            (4000000,1913.051551491)
            (5000000,1774.659923633)
        }
        ;
    \addplot[color={brown}, mark={triangle}]
        coordinates {
            (100000,20.659718125)
            (1000000,284.728816501)
            (2000000,683.646882805)
            (3000000,972.756284825)
            (4000000,1363.70826308)
            (5000000,1199.054598439)
        }
        ;
    \addplot[color={magenta}, mark={o}]
        coordinates {
            (100000,9.848747412)
            (1000000,122.270385225)
            (2000000,306.905139667)
            (3000000,414.078517222)
            (4000000,546.949608183)
            (5000000,557.335632431)
        }
        ;
\end{groupplot}
\end{tikzpicture}

%% file: figures/reserr_vp_noth.pgf
\begin{tikzpicture}
\begin{groupplot}[group style={group size={2 by 2}, ylabels at={edge left}, xlabels at={edge bottom}, vertical sep=1.5cm}]
    \nextgroupplot[height={5 cm}, width={7cm}, xlabel={}, ticks={both}, grid={both}, ymode={log}, legend pos={south west}, max space between ticks={45}, title={Harmonic SM}]
    \addplot[color={magenta}, mark={o}]
        coordinates {
            (1,0.0340162201906722)
            (2,0.02795730412802913)
            (3,0.020537087190162583)
            (4,0.011224816033908803)
            (5,0.013055667339593626)
            (6,0.010774975673073424)
            (7,0.012868087485440561)
            (8,0.011545601270997595)
            (9,0.006323313021178962)
            (10,0.008441430036855496)
            (11,0.005735740991950446)
            (12,0.007650845076299058)
            (13,0.002959934789083297)
            (14,0.004672162627662099)
            (15,0.0017627121766644997)
            (16,0.002923713609567899)
            (17,0.002799750542278529)
            (18,0.002606908906794832)
            (19,0.0010286228987793447)
            (20,0.002218720190305197)
            (21,0.0009112550344209832)
            (22,0.0006845031608085136)
            (23,0.0028093942354131675)
            (24,0.0011179255561385429)
            (25,0.0010877063459724253)
            (26,0.001183061561001269)
            (27,0.0006189167671227625)
            (28,0.0004181957173507422)
            (29,0.0007090139761161241)
            (30,0.0005094244353944333)
            (31,0.0007215800267942863)
            (32,0.0005387972973967587)
            (33,0.00039462940289472243)
            (34,0.0003225831436079098)
            (35,0.0005948622059424955)
            (36,0.00025481762248451517)
            (37,0.0003342759691882804)
            (38,0.000392530142025825)
            (39,0.00026699902059898744)
            (40,0.00018226208083190327)
            (41,0.00042862198886448884)
            (42,0.0001111668806711639)
            (43,0.00018033472956253204)
            (44,7.742094185255854e-5)
            (45,0.0001111705199882921)
            (46,5.1986149771085176e-5)
            (47,2.7154048653339534e-5)
            (48,1.7102666729297236e-5)
            (49,2.275214775184082e-5)
            (50,0.0001755695499045682)
            (51,0.00029095376702990513)
            (52,3.1432496205849596e-5)
            (53,0.000233229649248737)
            (54,7.117163829060605e-5)
            (55,0.00014935338207824712)
            (56,5.881500001395755e-5)
            (57,2.7979270799278188e-5)
            (58,1.0381979773630552e-5)
            (59,0.00020266629186107643)
            (60,4.860423647610773e-6)
            (61,1.0118677698543046e-5)
            (62,1.2754382619897128e-6)
            (63,1.871873881353013e-6)
            (64,6.797911191305904e-7)
            (65,5.074841545120013e-7)
            (66,9.366033086261446e-7)
            (67,2.2027699529519126e-7)
            (68,1.9683642506025547e-7)
            (69,5.813990378336906e-8)
            (70,5.58000599701387e-8)
            (71,2.7754545159693398e-9)
            (72,3.793198378357818e-9)
            (73,3.8728421893268565e-8)
            (74,7.299880870312012e-10)
            (75,6.775515454488263e-9)
            (76,1.3483898278560234e-10)
            (77,1.6251199475866298e-10)
            (78,5.990973308024417e-11)
        }
        ;
    \addlegendentry {rKS}
    \addplot[color={brown}, mark={triangle}]
        coordinates {
            (1,0.06136770154482913)
            (2,0.05674041739359513)
            (3,0.052275545256912095)
            (4,0.03813117376130646)
            (5,0.042526407242902305)
            (6,0.03186819098468721)
            (7,0.03252762760631303)
            (8,0.029673025454652297)
            (9,0.02458912519339205)
            (10,0.024821488788046738)
            (11,0.01896372050764663)
            (12,0.015454630949145034)
            (13,0.01218507886978889)
            (14,0.012379917088772157)
            (15,0.007848735470236026)
            (16,0.0091233237192166)
            (17,0.007727801458646162)
            (18,0.004544707566915246)
            (19,0.0034587732674274006)
            (20,0.002593889085807545)
            (21,0.002848827415048873)
            (22,0.0024509320936556963)
            (23,0.002051433128306006)
            (24,0.00200728397051888)
            (25,0.002669911000491524)
            (26,0.001695338353617985)
            (27,0.002001775673784231)
            (28,0.001355636589203443)
            (29,0.0015905607116541414)
            (30,0.0009702196167123331)
            (31,0.0027371805860532694)
            (32,0.0009436731965289117)
            (33,0.0011636904884657735)
            (34,0.0010621450773026832)
            (35,0.0005839461348850524)
            (36,0.0011105090455749873)
            (37,0.001204357418585798)
            (38,0.0011155348615288002)
            (39,0.0005682310729945579)
            (40,0.00037030896311190866)
            (41,0.00017209294437645726)
            (42,9.534715544551234e-5)
            (43,5.403922795199724e-5)
            (44,4.1475728879880116e-5)
            (45,4.929855782017612e-5)
            (46,7.294519377606574e-5)
            (47,0.00020023240565342971)
            (48,0.0005048117146661512)
            (49,0.0005815459719320847)
            (50,0.0003796582495038269)
            (51,0.0002770878741922951)
            (52,0.00015856583732744127)
            (53,0.0001578038839052688)
            (54,7.585536064490219e-5)
            (55,4.5938463603427836e-5)
            (56,1.852111211589876e-5)
            (57,1.4360681287175309e-5)
            (58,6.797100710750998e-6)
            (59,5.341206647458003e-6)
            (60,3.2295265808663987e-6)
            (61,3.021977293318353e-6)
            (62,1.6306246325264247e-6)
            (63,1.20231929824806e-6)
            (64,7.153811222966284e-7)
            (65,4.6965668035893864e-7)
            (66,2.969462286512651e-7)
            (67,1.756753403227379e-7)
            (68,9.202254141129764e-8)
            (69,8.209423909686224e-8)
            (70,3.2080225945518615e-8)
            (71,2.6148658769017376e-8)
            (72,1.573049404422491e-8)
            (73,1.0939852510937597e-8)
            (74,5.139339890984581e-9)
            (75,6.467912461494462e-9)
            (76,2.4595058636661374e-9)
            (77,3.2628247759951665e-9)
            (78,1.9316649558016987e-9)
            (79,2.060325443904058e-9)
            (80,8.238307617906019e-10)
            (81,1.262726395394613e-9)
            (82,5.101380042968583e-10)
            (83,4.694047928716107e-10)
            (84,3.2397355714450346e-10)
            (85,2.310651206208008e-10)
            (86,1.427966481830239e-10)
            (87,1.4639853123729443e-10)
            (88,6.293036381476066e-11)
        }
        ;
    \addlegendentry {KS}
    \nextgroupplot[footnotesize, xlabel={}, grid={both}, ticks={both}, height={5 cm}, width={7cm}, max space between ticks={45}, title={}, legend pos={south east}]
    \addplot[only marks, mark size={2}, color={magenta}, mark={o}]
        coordinates {
            (1,0.9756260544441594)
            (2,0.9779767516150925)
            (3,0.9781369284396734)
            (4,0.9786311111986147)
            (5,0.9787247646576656)
            (6,0.9787328420811524)
            (7,0.9789960626176576)
            (8,0.9790882702029721)
            (9,0.979172283042637)
            (10,0.9792439706314184)
            (11,0.9792477884569933)
            (12,0.9792481965127936)
            (13,0.9793012611659677)
            (14,0.9793258056606747)
            (15,0.9795234682025056)
            (16,0.9795325443476881)
            (17,0.9796097105695902)
            (18,0.9796248689116119)
            (19,0.9797288379629657)
            (20,0.9797296033077316)
            (21,0.9798242984932768)
            (22,0.9800660236477498)
            (23,0.9801111391279759)
            (24,0.9801838966233472)
            (25,0.9802155672953831)
            (26,0.9803069829969608)
            (27,0.9803972960197569)
            (28,0.9804058370325509)
            (29,0.9805155649957943)
            (30,0.980562274129648)
            (31,0.9805913161428353)
            (32,0.9806426706867767)
            (33,0.9807642202343823)
            (34,0.9807806256089454)
            (35,0.9808753394597283)
            (36,0.980887828400518)
            (37,0.981042140192285)
            (38,0.9811327218579604)
            (39,0.9811485282132161)
            (40,0.9811987296183057)
        }
        ;
    \addlegendentry {rKS}
    \addplot[only marks, mark size={2}, color={orange}, mark={pentagon}]
        coordinates {
            (1,0.9756260544441105)
            (2,0.9779767516151975)
            (3,0.9781369284396689)
            (4,0.9786311111986171)
            (5,0.9787247646576192)
            (6,0.9787328420811133)
            (7,0.9789960626176314)
            (8,0.9790882702029565)
            (9,0.9791722830425867)
            (10,0.979243970631401)
            (11,0.9792477884569473)
            (12,0.9792481965127875)
            (13,0.9793012611659718)
            (14,0.9793258056606601)
            (15,0.979523468202474)
            (16,0.9795325443476821)
            (17,0.979609710569575)
            (18,0.9796248689115954)
            (19,0.9797288379629288)
            (20,0.9797296033077254)
            (21,0.9798242984932959)
            (22,0.980066023647727)
            (23,0.9801111391279617)
            (24,0.980183896623334)
            (25,0.9802155672953548)
            (26,0.9803069829969351)
            (27,0.9803972960197238)
            (28,0.9804058370325438)
            (29,0.9805155649957724)
            (30,0.9805622741296345)
            (31,0.9805913161428257)
            (32,0.9806426706867438)
            (33,0.9807642202343632)
            (34,0.9807806256089258)
            (35,0.9808753394597381)
            (36,0.9808878284005105)
            (37,0.9810421401922326)
            (38,0.9811327218579777)
            (39,0.9811485282154214)
            (40,0.9811987296116611)
        }
        ;
    \addlegendentry {IRA}
    \addplot[only marks, mark size={2}, color={brown}, mark={triangle}]
        coordinates {
            (1,0.975626054444152)
            (2,0.9779767516151852)
            (3,0.9781369284396576)
            (4,0.9786311111986112)
            (5,0.9787247646576412)
            (6,0.9787328420811165)
            (7,0.978996062617639)
            (8,0.9790882702030091)
            (9,0.9791722830426417)
            (10,0.9792439706314101)
            (11,0.9792477884569485)
            (12,0.9792481965128128)
            (13,0.9793012611659949)
            (14,0.9793258056606566)
            (15,0.9795234682025077)
            (16,0.9795325443477024)
            (17,0.9796097105696172)
            (18,0.9796248689116387)
            (19,0.9797288379629443)
            (20,0.9797296033077785)
            (21,0.9798242984933254)
            (22,0.980066023647739)
            (23,0.9801111391279753)
            (24,0.9801838966233822)
            (25,0.9802155672953728)
            (26,0.9803069829969795)
            (27,0.9803972960197317)
            (28,0.9804058370325424)
            (29,0.9805155649958099)
            (30,0.9805622741296529)
            (31,0.9805913161428619)
            (32,0.9806426706868059)
            (33,0.9807642202343785)
            (34,0.9807806256089695)
            (35,0.9808753394597886)
            (36,0.9808878284005355)
            (37,0.9810421401922609)
            (38,0.9811327218580118)
            (39,0.9811485282154276)
            (40,0.9811987296116702)
        }
        ;
    \addlegendentry {KS}
    \nextgroupplot[height={5 cm}, width={7cm}, ticks={both}, grid={both}, ymode={log}, legend pos={south west}, max space between ticks={45}, title={Geometric LM}, xlabel={iteration number}]
    \addplot[color={magenta}, mark={o}]
        coordinates {
            (1,0.024513394702039303)
            (2,0.012768202281442022)
            (3,0.012770128893382695)
            (4,0.013449090514589243)
            (5,0.008539406244246523)
            (6,0.00833741305537824)
            (7,0.007776463866245341)
            (8,0.005889469602897761)
            (9,0.004735932352748177)
            (10,0.002266955566628555)
            (11,0.0022709960794555883)
            (12,0.004175677464535368)
            (13,0.0012812680838677276)
            (14,0.0015880237518184257)
            (15,0.0012870134280426539)
            (16,0.0012230200238188946)
            (17,0.0018196549234665489)
            (18,0.0008786336990592305)
            (19,0.00037584777864340383)
            (20,0.00036412777685900226)
            (21,0.0004859821885859078)
            (22,0.0011419024908145034)
            (23,0.0003532331198888488)
            (24,0.0003454875051023967)
            (25,0.00023830847992993256)
            (26,0.00024802065624233234)
            (27,0.0003586768728993849)
            (28,0.0002341451344287863)
            (29,0.0002575806044582467)
            (30,5.276665235693446e-5)
            (31,0.00017803544606162766)
            (32,7.490402553371163e-5)
            (33,9.372491328746932e-5)
            (34,1.4775199994440828e-5)
            (35,6.107066002484496e-6)
            (36,5.402504040264094e-7)
            (37,9.701144097111121e-9)
            (38,1.8418612329122482e-9)
            (39,1.7088517193619793e-10)
            (40,6.704625512996894e-11)
        }
        ;
    \addlegendentry {rKS}
    \addplot[color={brown}, mark={triangle}]
        coordinates {
            (1,0.04688297207033754)
            (2,0.036923686120499115)
            (3,0.027865103987588243)
            (4,0.03136341797903586)
            (5,0.023469333872958523)
            (6,0.020192954077607763)
            (7,0.010511373964902246)
            (8,0.009603628008054784)
            (9,0.008440228903626243)
            (10,0.006994952618308636)
            (11,0.006427046089970969)
            (12,0.00799060782057968)
            (13,0.002747198385011122)
            (14,0.0030428290620244983)
            (15,0.0017747805310532813)
            (16,0.0017253902576020025)
            (17,0.0010901168795593242)
            (18,0.001621321109002949)
            (19,0.0008971276129527372)
            (20,0.0015961083707230668)
            (21,0.0014159886620990286)
            (22,0.0017579663969909004)
            (23,0.0011192338363681272)
            (24,0.000506515810996662)
            (25,0.0010332191713026246)
            (26,0.0007275417567278224)
            (27,0.0007615066089773671)
            (28,5.253633999502169e-5)
            (29,0.0009240577494762045)
            (30,0.0003163039118205201)
            (31,0.00010160728659397039)
            (32,1.759659541019747e-5)
            (33,2.3235474617915675e-6)
            (34,2.476804178325276e-7)
            (35,3.625229094348817e-8)
            (36,6.5028980660614025e-9)
            (37,1.1369037360377875e-9)
            (38,2.3115287660444538e-10)
            (39,4.67784346080992e-11)
        }
        ;
    \addlegendentry {KS}
    \legend{}
    \nextgroupplot[footnotesize, grid={both}, ticks={both}, height={5 cm}, width={7cm}, max space between ticks={45}, title={}, xlabel={eigenvalue label}]
    \addplot[only marks, mark size={2}, color={magenta}, mark={o}]
        coordinates {
            (1,1.0125274848431645)
            (2,1.0124774775074012)
            (3,1.0122270587382032)
            (4,1.0117568563365111)
            (5,1.011444239286269)
            (6,1.0113338165373273)
            (7,1.0112174928295077)
            (8,1.0109444904587277)
            (9,1.0107757660721721)
            (10,1.010643455557724)
            (11,1.0101211934271044)
            (12,1.0093621929718375)
            (13,1.0090219450852098)
            (14,1.0088291488420578)
            (15,1.0087266034138525)
            (16,1.0087093372188638)
            (17,1.0084804326477579)
            (18,1.0083699520789655)
            (19,1.0082546005040718)
            (20,1.0080034552188732)
            (21,1.0079436348544912)
            (22,1.0079111674601173)
            (23,1.0077858073294257)
            (24,1.0077748360633676)
            (25,1.0077699694046278)
            (26,1.007741539321359)
            (27,1.007690369349395)
            (28,1.0076852765892563)
            (29,1.0076796398853805)
            (30,1.0076400126238525)
            (31,1.0076186300755063)
            (32,1.0074226524861252)
            (33,1.0072290245375557)
            (34,1.0071555474245617)
            (35,1.007142748301617)
            (36,1.0070619496067377)
            (37,1.0070230263416262)
            (38,1.0069702153329827)
            (39,1.0069626925552786)
            (40,1.0068782454706364)
        }
        ;
    \addlegendentry {rKS}
    \addplot[only marks, mark size={2}, color={orange}, mark={pentagon}]
        coordinates {
            (1,1.012527484843167)
            (2,1.0124774775073855)
            (3,1.0122270587381943)
            (4,1.0117568563365076)
            (5,1.0114442392862828)
            (6,1.0113338165373582)
            (7,1.011217492829493)
            (8,1.0109444904587157)
            (9,1.0107757660721672)
            (10,1.010643455557698)
            (11,1.010121193427087)
            (12,1.009362192971866)
            (13,1.0090219450852222)
            (14,1.0088291488420549)
            (15,1.0087266034138433)
            (16,1.0087093372188554)
            (17,1.0084804326477663)
            (18,1.008369952078974)
            (19,1.0082546005040776)
            (20,1.0080034552188666)
            (21,1.00794363485448)
            (22,1.007911167460109)
            (23,1.007785807329434)
            (24,1.0077748360633594)
            (25,1.007769969404629)
            (26,1.007741539321351)
            (27,1.0076903693493948)
            (28,1.0076852765892552)
            (29,1.0076796398853816)
            (30,1.007640012623857)
            (31,1.007618630075493)
            (32,1.0074226524861196)
            (33,1.0072290245375495)
            (34,1.0071555474245693)
            (35,1.007142748301649)
            (36,1.0070619496066218)
            (37,1.007023026340867)
            (38,1.00697021532322)
            (39,1.0069626925651851)
            (40,1.006878245447922)
        }
        ;
    \addlegendentry {IRA}
    \addplot[only marks, mark size={2}, color={brown}, mark={triangle}]
        coordinates {
            (1,1.01252748484319)
            (2,1.0124774775074044)
            (3,1.0122270587381683)
            (4,1.0117568563365131)
            (5,1.0114442392862868)
            (6,1.0113338165373642)
            (7,1.0112174928295177)
            (8,1.0109444904587248)
            (9,1.010775766072163)
            (10,1.010643455557725)
            (11,1.0101211934270995)
            (12,1.0093621929718648)
            (13,1.0090219450852256)
            (14,1.0088291488420496)
            (15,1.0087266034138673)
            (16,1.0087093372188427)
            (17,1.0084804326477732)
            (18,1.0083699520789786)
            (19,1.0082546005040636)
            (20,1.0080034552188604)
            (21,1.0079436348544968)
            (22,1.0079111674601389)
            (23,1.007785807329422)
            (24,1.007774836063359)
            (25,1.0077699694046252)
            (26,1.0077415393213587)
            (27,1.0076903693493926)
            (28,1.0076852765892563)
            (29,1.0076796398853756)
            (30,1.0076400126238543)
            (31,1.007618630075496)
            (32,1.0074226524861198)
            (33,1.0072290245375577)
            (34,1.0071555474245637)
            (35,1.0071427483016504)
            (36,1.007061949606636)
            (37,1.0070230263408722)
            (38,1.0069702153232234)
            (39,1.0069626925651762)
            (40,1.0068782454479366)
        }
        ;
    \addlegendentry {KS}
    \legend{}
\end{groupplot}
\end{tikzpicture}

%% file: matrix-functions.tex
\section{Evaluation of matrix functions}
\label{sec:matrix-functions}

Given a matrix $A \in \R^{n \times n}$ and a function $f$ that is analytic on and inside a contour $\Gamma \subset \C$ which encloses the spectrum of $A$, the matrix function $f(A)$ can be defined as
\begin{equation*}
	f(A) = \frac{1}{2 \pi i}\int_\Gamma f(t) (tI - A)^{-1} \, dt.
\end{equation*}
We refer to \cite{Higham08} for other equivalent definitions for $f(A)$.  
The computation of matrix functions arises in many areas, such as the solution of partial differential equations \cite{BotchevKnizhnerman20}, network analysis \cite{BenziBoito20, EstradaHigham10}, and electronic structure computations \cite{BenziBoitoRazouk13}. In these applications, one is often interested in the computation of $f(A) \vec b$ for a given vector $\vec b \in \R^n$, rather than the full matrix function $f(A)$. When $A$ is large and sparse, the computation of $f(A)$ through methods such as Schur-Parlett \cite{DaviesHigham03} is usually infeasible as it has a computational cost of $O(n^3)$, especially if the (generally dense) matrix $f(A)$ is too large to store explicitly. In this setting, Krylov subspace methods are the most popular methods for the approximation of $f(A) \vec b$ \cite{DruskinKnizhnerman89, Saad92}. 
In this section, we discuss the computation of $f(A) \vec b$ with randomized Krylov methods, which has been recently investigated in \cite{CKN24, GuettelSchweitzer23, PSS25mf}. We are going to present the different approximations that have been proposed in the literature, emphasizing the links and equivalences between the different approaches.

A well-established way to approximate $f(A) \vec b$ is to use a projection onto the Krylov subspace $\K_m(A, \vec b)$, which yields the Arnoldi approximation \cite{Saad92,DruskinKnizhnerman89}
\begin{equation}
	\label{eqn:fAb-arnoldi-approx-orth}
	\vec f_m := \widetilde{\beta} Q_m f(H_m) \vec e_1,
\end{equation} 
where $Q_m$ and $H_m$ satisfy the (orthogonal) Arnoldi relation \cref{eqn:arnoldi-relation} and $\vec b = \widetilde{\beta} Q_m \vec e_1$. This approximation is exact when $f$ is a polynomial of degree up to $m-1$, and it is equivalent to computing $p_{m-1}(A) \vec b$, where $p_{m-1}$ is a polynomial which interpolates $f$ at the eigenvalues of $H_m$, see for instance \cite[Theorem~3.3]{Saad92}. 

Given an arbitrary basis $W_m$ of $\K_m(A, \vec b)$ and the associated Hessenberg matrix $\underline{L}_m$ satisfying the Arnoldi-like relation \cref{eqn:arnoldi-like-relation}, the approximation \cref{eqn:fAb-arnoldi-approx-orth} satisfies 
\begin{equation}
	\label{eqn:fAb-arnoldi-approx-general}
	\vec f_m = W_m f(W_m^+ A W_m) W_m^+ \vec b = \widetilde{\beta} W_m f(L_m + W_m^+ \vec w_{m+1} \vec \ell_m^T) \vec e_1,
\end{equation}
where we used \cref{eqn:arnoldi-projection-identity--general-basis} and $W_m^+ \vec b = W_m^+ (\widetilde{\beta} \vec w_1) = \widetilde{\beta} \vec e_1$ to rewrite the identity. See, for instance, \cite[Lemma~3.1]{CKN24} for a proof of the equivalence between \cref{eqn:fAb-arnoldi-approx-orth} and \cref{eqn:fAb-arnoldi-approx-general}. 
When the basis $W_m$ is not orthonormal, the downside of the approximation \cref{eqn:fAb-arnoldi-approx-general} is that computing $W_m^+ \vec w_{m+1}$ requires the solution of the least squares problem 
\begin{equation*}
	W_m^+ \vec w_{m+1} = \argmin_{\vec h \in \R^m} \norm{W_m \vec h - \vec w_{m+1}}.
\end{equation*}
When $W_m$ is sketch-orthonormal it follows from \cref{eqn:singular-values-eps-embedding} that $\cond(W_m) \le \sqrt{(1+\varepsilon)/(1-\varepsilon)}$, so $W_m$ is well-conditioned and $W_m^+ \vec w_{m+1}$ can be computed by solving a least squares problem with the LSQR algorithm. This approach is proposed in \cite[Algorithm~3.1]{CKN24}. For a general basis $W_m$, in \cite[Algorithm~3.2]{CKN24} it is proposed to solve the least squares problem for $W_m^+ \vec w_{m+1}$ by using Blendenkip \cite{AMT10Blendenpik}. 

A further alternative that has been explored in \cite[Algorithm~3.3]{CKN24} is to completely ignore the rank-one update in \cref{eqn:fAb-arnoldi-approx-general} and approximate $f(A) \vec b$ with
\begin{equation*}
	\widetilde{\vec f}_m = \widetilde{\beta} W_m f(L_m) \vec e_1,
\end{equation*}
but for a general basis $W_m$, the matrix $f(L_m)$ can be quite far from $f(L_m + W_m^+ \vec w_{m+1} \vec \ell_m^T)$, so $\widetilde{\vec f}_m$ may be an inaccurate approximation of $f(A) \vec b$.  
However, with a sketch-orthonormal basis of $\K_m(A, \vec b)$, we show below that an approximation of this form has a natural link with \cref{eqn:fAb-arnoldi-approx-orth}. Assume that we have a sketch-orthonormal basis $V_m$ and a corresponding Hessenberg matrix $\underline{G}_m$ which satisfy the randomized Arnoldi relation \cref{eqn:arnoldi-like-relation--sketched-basis}, and let $\vec b = \beta V_m \vec e_1$. Then we have
\begin{equation}
	\label{eqn:fAb-arnoldi-approx-sketched}
	\vec f_m^\Omega := \beta V_m f(G_m) \vec e_1 = V_m f\left((\Omega V_m)^T \Omega A V_m\right) (\Omega V_m)^T \Omega \vec b,
\end{equation}
where for the last equality we used the sketch-orthonormality of $V_m$ and the identity $G_m = (\Omega V_m)^T \Omega A V_m$. 
This approximation is introduced in \cite{GuettelSchweitzer23} with the name \emph{sketched FOM} (sFOM). To describe its connection with \cref{eqn:fAb-arnoldi-approx-orth}, we follow the presentation in \cite[Section~2]{GuettelSchweitzer23} and consider a function that admits an integral representation of the form
\begin{equation*}
	f(z) = \int_\Gamma (t + z)^{-1} \, d\mu(t),
\end{equation*}
which includes as special cases both Stieltjes functions \cite{Berg08} and the Cauchy integral representation for analytic functions. The integral expression for $f$ translates to the matrix function representation
\begin{equation}
    \label{eqn:fAb-integral-expression}
	f(A) \vec b = \int_\Gamma (t I + A)^{-1} \vec b \, d \mu(t) = \int_\Gamma \vec x(t) \, d\mu(t), \qquad \text{where} \quad (tI + A) \vec x(t) = \vec b.
\end{equation}  
For this class of functions, the Arnoldi approximation \cref{eqn:fAb-arnoldi-approx-orth} can be interpreted as follows. For each $t \in \Gamma$, let $\vec x_m(t) := \widetilde{\beta} Q_m (tI + H_m)^{-1} \vec e_1$ be the approximate solution to the shifted linear system $(tI + A) \vec x(t) = \vec b$ after $m$ iterations of FOM. Then, we have
\begin{equation*}
	\vec f_m = \widetilde{\beta} Q_m f(H_m) \vec e_1 = \int_\Gamma \widetilde{\beta} Q_m (tI + H_m)^{-1} \vec e_1 \, d\mu(t) = \int_\Gamma \vec x_m(t) \, d\mu(t).
\end{equation*}
The residuals $\vec r_m(t) = \vec b - (tI + A) \vec x_m(t)$ are orthogonal to $\K_m(A, \vec b) = \range(Q_m)$. 

The sFOM approximation \cref{eqn:fAb-arnoldi-approx-sketched} stems from the following observation: instead of imposing the orthogonality condition on $\vec r_m(t)$ exactly, we can alternatively solve the shifted linear systems using randomized FOM, i.e., impose that the sketched residuals are orthogonal to the sketch of the Krylov subspace. In other words, for each $t \in \Gamma$ we look for an approximate solution $\vec x_m^\Omega(t) \in \K_m(A, \vec b)$ such that the residual $\vec r_m^\Omega(t) = \vec b - (tI + A) \vec x_m^\Omega(t)$ satisfies the sketched Galerkin condition $\Omega \vec r_m^\Omega(t) \perp \Omega \K_m(A, \vec b)$. This implies that $\vec x_m^\Omega(t) = V_m \vec y_m^\Omega(t)$ with
\begin{equation*}
	(\Omega V_m)^T (\Omega\vec b - \Omega (tI + A)V_m \vec y_m^\Omega(t)) = \vec 0,
\end{equation*}
and hence
\begin{equation}
	\label{eqn:sketched-fom-shifted-linsys}
	\vec x_m^\Omega(t) = \beta V_m (t I + (\Omega V_m)^T \Omega A V_m)^{-1} \vec e_1 = \beta V_m (t I + G_m)^{-1} \vec e_1,
\end{equation}
where we used the identities $\vec b = \beta \vec v_1$, $(\Omega V_m)^T \Omega V_m = I$ and $(\Omega V_m)^T \Omega A V_m = G_m$. 
It then follows that 
\begin{equation*}
	\int_\Gamma \vec x_m^\Omega(t) d \mu(t) = \int_\Gamma \beta V_m (t I + G_m)^{-1} \vec e_1 = \beta V_m f(G_m) \vec e_1 = \vec f_m^\Omega,
\end{equation*} 
showing that the approximation $\vec f_m^\Omega$ can be obtained by imposing a sketched Galerkin condition on the residuals of the shifted linear systems in the integral expression \cref{eqn:fAb-integral-expression} for $f(A) \vec b$. The relation between \cref{eqn:fAb-arnoldi-approx-sketched} and \cref{eqn:fAb-arnoldi-approx-orth} then mimics the one between the approximate solutions for linear systems obtained with randomized FOM~\cref{eqn:linear-system-sketched-fom-sol} and FOM~\cref{eqn:linear-system-fom-sol}.
We refer to \cite[Section~2]{GuettelSchweitzer23} for further details on the sFOM approximation for $f(A) \vec b$.

The authors of \cite{GuettelSchweitzer23} also consider a sketched GMRES approximation, in which the shifted linear systems $(t I + A) \vec x(t) = \vec b$ are solved by using randomized GMRES instead of randomized FOM. However, this approximation to $f(A) \vec b$ has no simple closed-form expression, so it requires using a quadrature rule to evaluate the integral expression of $f$. We refer to \cite[Section~3]{GuettelSchweitzer23} for additional details.

Both in \cite{CKN24} and in \cite{GuettelSchweitzer23}, it is proposed to compute the approximation \cref{eqn:fAb-arnoldi-approx-sketched} by using a sketch-orthonormal basis $V_m$ obtained implicitly through the whitening approach desribed in \cref{subsubsec:whitening}. Specifically, a non-orthonormal basis $W_m$ for $\K_m(A, \vec b)$ is constructed using the $k$-truncated Arnoldi process, and then a sketch-orthonormal basis $V_m = W_m R_m^{-1}$ is obtained, where $\Omega W_m = S_m R_m$ is a QR factorization. In this setting, we can compute the approximation \cref{eqn:fAb-arnoldi-approx-sketched} without ever forming the basis $V_m$ explicitly. Indeed, we have \cite[Section~2.1]{GuettelSchweitzer23}
\begin{equation*}
	\vec f_m^\Omega = V_m f(G_m) \beta \vec e_1 = W_m R_m^{-1} f\left(S_m^T \Omega A W_m R_m^{-1}\right) S_m^T \Omega \vec b,
\end{equation*}
and in this expression we only need to apply $R_m^{-1}$ to the right of $S_m^T \Omega A V_m \in \R^{m \times m}$ and to the left of the vector $f(S_m^T \Omega A V_m R_m^{-1}) S_m^T \Omega \vec b \in \R^m$, for a cost of $\mathcal{O}(m^3)$. Recall that, on the other hand, explicitly computing $V_m R_m^{-1}$ would cost $\mathcal{O}(n m^2)$, which would eliminate the computational adavantage of the whitening approach.

The approximation \cref{eqn:fAb-arnoldi-approx-sketched} is also studied in \cite{PSS25mf}, where the authors also employ the $k$-truncated Arnoldi process combined with whitening in order to implicitly construct the sketch-orthonormal basis $V_m$. In particular, they derive the identities \cref{eqn:whitened-sketched-arnoldi-relation} and \cref{eqn:whitened-hessenberg-matrix} and use them to rewrite \cref{eqn:fAb-arnoldi-approx-sketched} as
\begin{equation}
	\label{eqn:fAb-arnoldi-approx-whitened-sketched}
	\vec f^\Omega_m = W_m R_m^{-1} f(\widehat{L}_m + \rho_m^{-1} \vec z_m \vec \ell_m^T) \beta \vec e_1,
\end{equation}
see \cite[Algorithm~1]{PSS25mf} for the implementation details. in \cite[Section~7]{PSS25mf}, the authors argue that the approximation \cref{eqn:fAb-arnoldi-approx-whitened-sketched} is quite robust in finite-precision arithmetic despite the potentially ill-conditioned matrix $R_m$.

We also mention that in the recent preprint \cite{GMAM25}, sketch-orthogonalization is proposed as a means to reduce the orthogonalization costs of restarted Krylov methods for matrix functions \cite{EiermannErnst06}. Although the authors do not provide a rigorous convergence analysis, the method they present exhibits competitive performance in their numerical tests, occasionally even converging in fewer iterations compared to the deterministic restarted method; see, e.g., \cite[Figure~3]{GMAM25}.

%% file: matrix-equations.tex
\section{Solution of matrix equations}
\label{sec:matrix-equations}

Matrix Sylvester equations appear in numerous applications, for instance in model order reduction and in the discretization of certain partial differential equations; we refer to \cite{BennerSaak13,Simoncini16} for additional details. In several applications, the right-hand side is low-rank and the matrix Sylvester equation can be written in the form
\begin{equation}
	\label{eqn:sylvester-equation}
	A X + X B = C_1 C_2^T,
\end{equation}
with $A, B \in \R^{n \times n}$ and $C_1, C_2 \in \R^{n \times r}$, with $r \ll n$. In this setting, efficient approaches for the solution of \cref{eqn:sylvester-equation} are often based on projection on the polynomial block Krylov subspaces $\K_m(A, C_1)$ and $\K_m(B^T, C_2)$, or on extended and rational block Krylov subspaces; we refer to the review paper \cite{Simoncini16} for further details and references.

In \cite{PSS25me}, it was proposed to use randomized sketching to reduce the cost of orthogonalization within Krylov methods for the solution of \cref{eqn:sylvester-equation}. In this section, we briefly describe the approach presented in \cite{PSS25me}.
Assume that we generate two non-orthonormal bases $W_m^A$ and $W_m^B \in \R^{n \times mr}$ of $\K_m(A, C_1)$ and $\K_m(B^T, C_2)$, using a truncated block Arnoldi procedure, where orthogonalization is only performed against the previous $k$ blocks, which leads to the Arnoldi relations
\begin{equation*}
	A W_m^A = W_{m+1}^A \underline{H}_m^A \qquad \text{and} \qquad B^T W_m^B = W_{m+1}^B \underline{H}_m^B,
\end{equation*}
where $\underline{H}_m^A$ and $\underline{H}_m^B \in \R^{(m+1) r \times mr}$ are block upper Hessenberg matrices with upper bandwidth $kr$. 
Assume that we have two $\varepsilon$-subspace embeddings $\Omega_A$ and $\Omega_B$ for $\K_{m+1}(A, C_1)$ and $\K_{m+1}(B^T, C_2)$, respectively. Then, we can use the basis whitening approach to cheaply compute sketch-orthonormal bases of the two block Krylov subspaces: given the QR factorizations $\Omega^A W_m^A = Q_m^A T_m^A$ and $\Omega^B W_m^B = Q_m^B T_m^B$, the whitened bases 
\begin{equation*}
	V_m^A = W_m^A (T_m^A)^{-1} \qquad \text{and} \qquad V_m^B = W_m^B (T_m^B)^{-1}
\end{equation*}     
are sketch-orthogonal due to the $\epsilon$-embedding property. Associated to the bases $V_m^A$ and $V_m^B$ are the modified block upper Hessenberg matrices $G_m^A = T_m^A H_m^A (T_m^A)^{-1}$ and $G_m^B = T_m^B H_m^B (T_m^B)^{-1}$. For further details we refer to the whitened-sketched Arnoldi relations presented in \cite[Section~2.3]{PSS25me}, which generalize \cref{eqn:whitened-sketched-arnoldi-relation} to the block case.

Projection methods for the Sylvester matrix equation \cref{eqn:sylvester-equation} usually look for an approximate solution of the form $X_m = Q_m^A Y_m (Q_m^B)^T$, with corresponding residual $R_m = A X_m + X_m B - C_1 C_2^T$, where $Q_m^A$ and $Q_m^B$ are orthonormal bases of the block Krylov subspaces $\K_m(A, C_1)$ and $\K_m(B^T, C_2)$, and $Y_m \in \R^{mr \times mr}$ is computed by solving a smaller projected problem. The authors of \cite{PSS25me} propose a sketched-and-truncated method that searches for a solution of the form $X_m = V_m^A Y_m (V_m^B)^T$, and imposes on the associated residual matrix $R_m$ the following \emph{sketched Galerkin condition}
\begin{equation}
	(\Omega^A V_m^A)^T (\Omega^A R_m (\Omega^B)^T) (\Omega^B V_m^B) = 0,
\end{equation} 
which can be satisfied by taking as $Y_m$ the solution of the projected equation  
\begin{equation}
	(G_m^A + \widehat{G}^A E_m^T) Y_m + Y_m (G_m^B + \widehat{G}^B E_m^T)^T = E_1 \beta_1 \beta_2^T E_1^T,
\end{equation}
where $\widehat{G}^A$ and $\widehat{G}^B \in \R^{n \times r}$ are suitable block vectors which can be obtained when computing the whitened bases $V_m^A$ and $V_m^B$, and the block scalars $\beta_1$, $\beta_2 \in \R^{r \times r}$ can be determined from the expressions $C_1 = W_m^A E_1 \beta_1$, $C_2 = W_m^B E_1 \beta_2$, where $E_i \in \R^{mr \times r}$ denotes the $i$-th block column of a $mr \times mr$ identity matrix. 
We refer to \cite[Section 3]{PSS25me} for a detailed description of their algorithm and explicit expressions for $\widehat{G}^A$ and $\widehat{G}^B$, and to \cite[Algorithm~1]{PSS25me} for a pseudocode.